\documentclass[11pt,reqno,a4paper]{amsart}

\usepackage[T1]{fontenc}
\usepackage[utf8]{inputenc}

\usepackage{libertinus}
\usepackage{microtype}
\usepackage{geometry}
\usepackage{amsmath,amssymb,amsthm,amsfonts}
\usepackage{mathtools}
\usepackage{mathrsfs}
\usepackage{bbm}

\numberwithin{equation}{section}

\usepackage{graphicx}
\usepackage{xcolor}
\usepackage{enumitem}
\setlist{itemsep=2pt,topsep=4pt}

\usepackage[
  colorlinks=true,
  linkcolor=blue,
  citecolor=blue,
  urlcolor=blue,
  pagebackref=true
]{hyperref}

\renewcommand*{\backrefalt}[4]{%
  \ifcase #1\relax
  \or
    {\footnotesize Cited on page #2.}%
  \else
    {\footnotesize Cited on pages #2.}%
  \fi
}

\usepackage[nameinlink,noabbrev]{cleveref}
\crefname{section}{Section}{Sections}
\Crefname{section}{Section}{Sections}
\crefname{subsection}{Subsection}{Subsections}
\Crefname{subsection}{Subsection}{Subsections}

\usepackage{titlesec}

\titleformat{\section}
  {\normalfont\large\bfseries\centering}
  {\thesection.}
  {0.75em}
  {}

\titlespacing*{\section}
  {0pt}
  {3.2ex plus 1ex minus .2ex}
  {1.6ex plus .3ex}

\titleformat{\subsection}
  {\normalfont\normalsize\bfseries}
  {\thesubsection.}
  {0.75em}
  {}

\titlespacing*{\subsection}
  {0pt}
  {2.4ex plus .8ex minus .2ex}
  {1.0ex plus .2ex}

\titleformat{\subsubsection}
  {\normalfont\normalsize\bfseries\itshape}
  {\thesubsubsection.}
  {0.75em}
  {}

\titlespacing*{\subsubsection}
  {0pt}
  {2.0ex plus .6ex minus .2ex}
  {0.7ex plus .2ex}

\usepackage{aliascnt}
\usepackage{etoolbox}

\newtheoremstyle{compact}
  {6pt}
  {1pt}
  {\itshape}
  {}
  {\bfseries}
  {.}
  {0.5em}
  {}

\theoremstyle{compact}
\newtheorem{theorem}{Theorem}[section]

\newaliascnt{proposition}{theorem}
\newtheorem{proposition}[proposition]{Proposition}
\aliascntresetthe{proposition}

\newaliascnt{lemma}{theorem}
\newtheorem{lemma}[lemma]{Lemma}
\aliascntresetthe{lemma}

\newaliascnt{corollary}{theorem}
\newtheorem{corollary}[corollary]{Corollary}
\aliascntresetthe{corollary}

\newtheorem{maintheorem}{Theorem}

\theoremstyle{definition}
\newaliascnt{definition}{theorem}

\aliascntresetthe{definition}

\theoremstyle{definition}
\newtheorem{remark}[theorem]{Remark}
\newtheorem*{remark*}{Remark}

\crefname{theorem}{Theorem}{Theorems}
\Crefname{theorem}{Theorem}{Theorems}
\crefname{proposition}{Proposition}{Propositions}
\Crefname{proposition}{Proposition}{Propositions}
\crefname{lemma}{Lemma}{Lemmas}
\Crefname{lemma}{Lemma}{Lemmas}
\crefname{corollary}{Corollary}{Corollaries}
\Crefname{corollary}{Corollary}{Corollaries}
\crefname{maintheorem}{Theorem}{Theorems}
\Crefname{maintheorem}{Theorem}{Theorems}
\crefname{maincorollary}{Corollary}{Corollaries}
\Crefname{maincorollary}{Corollary}{Corollaries}

\AtBeginEnvironment{proof}{\vspace{6pt}}

\makeatletter
\AtBeginEnvironment{thebibliography}{%
  \def\@mklab#1{#1\hfil}%
}
\makeatother

\newcommand{\Ad}{\mathbf{A}}
\newcommand{\PP}{\mathbf{PP}}
\newcommand{\EPP}{\mathbf{EPP}}
\newcommand{\RR}{\mathbf{R}}
\newcommand{\Hyp}{\mathbf{Hyp}}
\newcommand{\IHyp}{\mathbf{IHyp}}
\newcommand{\FHyp}{\mathbf{FHyp}}
\newcommand{\PNC}{\mathbf{PNC}}
\newcommand{\PerPP}{\mathbf{Per}}
\newcommand{\MRig}{\mathbf{MRig}}
\newcommand{\LMRig}{\mathbf{LMRig}}
\newcommand{\Rev}{\mathbf{Rev}}
\newcommand{\MFVH}{\mathbf{MFVH}}

\newcommand{\Var}{\operatorname{Var}}
\newcommand{\Cov}{\operatorname{Cov}}
\newcommand{\Vol}{\operatorname{Vol}}

\begin{document}

\title{Generic Point Processes Are Intermittently Hyperuniform and Maximally Rigid}

\author{Michael Bj\"orklund}
\address{Department of Mathematics, Chalmers University of Technology and University of Gothenburg, Gothenburg, Sweden}
\email{micbjo@chalmers.se}

\begin{abstract}
We study the prevalence of hyperuniformity and spatial rigidity among
translation-invariant ergodic simple point processes of fixed intensity, in a topology
controlling local second moments. Hyperuniform processes are dense but
meagre. We prove that intermittent hyperuniformity is residual: the number variance
exhibits hyperuniform behavior along some sequences of radii and
hyperfluctuating behavior along others. Maximal rigidity is also residual. Thus generic processes
are maximally rigid but not hyperuniform. A local-to-global principle
reduces the proofs to cross-section constructions for an arbitrary fixed
essentially free ergodic $\mathbb R^d$-action. We also prove that processes admitting revival
sequences form a residual set and exhibit near-quadratic variance along
suitable van Hove sequences.
\end{abstract}

\maketitle

\section{Introduction and main results}
\label{sec:introduction}

Hyperuniformity means that the variance of the number of points in a
large Euclidean ball grows more slowly than its volume. It occurs in
periodic point configurations, certain quasicrystals, and some disordered
systems \cite{TS03,Tor18}. We call a point process \emph{number rigid} if
its configuration outside any bounded region determines the number of
points inside, and \emph{maximally rigid} if it determines their
positions as well. Examples arise from determinantal kernels, zeros of
Gaussian analytic functions, and spectral gaps
\cite{GP17,GL17rig,GL18,GK21,LR26}.

Ghosh and Lebowitz \cite{GL17rig} proved number rigidity for
hyperuniform processes in dimensions one and two under additional decay
assumptions on truncated pair correlations. They also raised the question
whether hyperuniformity, possibly with correlation decay, is necessary
for number rigidity. In \cite{GL18}, they established maximal rigidity
for processes whose structure function vanishes on an open set, and also
under superpolynomial vanishing at the origin. Klatt and Last \cite{KL22}
constructed hyperfluctuating point processes for which the entire
configuration in every bounded region can be recovered from a bounded
exterior stopping set. Their \emph{strong rigidity} implies maximal
rigidity as defined here. Conversely, in dimensions at least three,
sufficiently large independent Gaussian displacements of lattices give
hyperuniform processes that are not number rigid \cite{PS14,BB26}.
Thus hyperuniformity does not imply number rigidity in dimensions at
least three, while the hyperfluctuating examples of Klatt and Last show
that number rigidity does not imply hyperuniformity.

These results concern particular families of point processes. In this paper,
we study the Baire-category prevalence of hyperuniformity and rigidity
among translation-invariant ergodic simple point processes of fixed intensity
and finite local second moments. We prove, for the local $L^2$ topology
defined below, that hyperuniformity is dense but meagre, whereas linear
maximal rigidity is residual.

We identify a locally finite point set with its counting measure. Let
\(\operatorname{Conf}(\mathbb R^d)\) be the space of locally finite simple
counting measures, equipped with the vague topology; see
\cite[Appendix~A2.6]{DVJ03} and \cite{Pug02}.

A point process is a probability measure \(\eta\) on
\(\operatorname{Conf}(\mathbb R^d)\). For \(t\in\mathbb R^d\), set
\[
(\theta_t\omega)(A):=\omega(A+t),
\]
and call \(\eta\) \emph{translation-invariant} if
\((\theta_t)_*\eta=\eta\) for every \(t\in\mathbb R^d\). For a bounded Borel
set \(A\subset\mathbb R^d\), let
\[
N_A(\omega):=\omega(A).
\]
Such a process has \emph{intensity} \(\rho>0\) if
\[
\int N_A\,d\eta=\rho\Vol_d(A)
\]
for every bounded Borel set \(A\subset\mathbb R^d\). Write
\(\PP_{\rho,2}\) for the translation-invariant point processes of intensity
\(\rho\) whose counts also satisfy
\[
N_A\in L^2(\eta)
\qquad
\text{for every bounded Borel }A\subset\mathbb R^d.
\]
Finite intensity alone does not imply local square-integrability.

We equip \(\PP_{\rho,2}\) with the \emph{local \(L^2\) topology}.
In this topology, convergence means weak convergence of point-process laws
together with convergence of their count second moments in each cube
\(Q_m:=[-m/2,m/2)^d\), \(m\geq1\).
We prove the precise convergence criterion in
\cref{prop:local-L2-characterization} and show in
\cref{prop:global-polish} that \(\PP_{\rho,2}\) is Polish and its
subspace \(\EPP_{\rho,2}\) of ergodic processes is a
\(G_\delta\)-subset, hence also Polish.

A set is \emph{nowhere dense} if its closure has empty interior, and
\emph{meagre} if it is a countable union of nowhere-dense sets. A set is
\emph{residual} if its complement is meagre. In a Polish space this is
equivalent to containing a dense \(G_\delta\)-subset.

\subsection{Hyperuniformity}
\label{subsec:intro-hyperuniformity}

For \(R>0\), let \(B_R\) be the closed Euclidean ball of radius \(R\)
centred at the origin. We write
\[
V_\eta(R):=\Var_\eta(N_{B_R}).
\]
We call \(\eta\in\PP_{\rho,2}\) \emph{hyperuniform} if its ball-count
variance grows sublinearly in volume, that is, if
\[
\frac{V_\eta(R)}{\Vol_d(B_R)}\longrightarrow0
\qquad(R\to\infty).
\]
Let \(\Hyp_\rho\subset\PP_{\rho,2}\) denote the class of hyperuniform processes.

For $d=1$, let \(\Vol_0\) denote counting measure. Beck's discrepancy
theorem \cite[Theorem~2A]{Bec87} and the spectral lower bound of
\cite[Theorem~5.1]{BB26} imply that every
\(\eta\in\PP_{\rho,2}\), without an ergodicity assumption, satisfies
\begin{equation}\label{eq:variance-boundary-limsup}
\limsup_{R\to\infty}
\frac{V_\eta(R)}{\Vol_{d-1}(\partial B_R)}>0.
\end{equation}
Equivalently, there are
\(c=c(\eta)>0\) and radii \(R_k\to\infty\) such that
\[
V_\eta(R_k)\geq c\,\Vol_d(B_{R_k})^{(d-1)/d}.
\]
At the other end of the scale, every ergodic
\(\eta\in\PP_{\rho,2}\) satisfies
\[
V_\eta(R)=o\!\left(\Vol_d(B_R)^2\right)
\qquad(R\to\infty),
\]
by the mean ergodic theorem; see \cref{prop:van-hove-subquadratic}.

We call \(\eta\in\smash{\EPP_{\rho,2}}\) \emph{intermittently hyperuniform}
if, for every
\(\beta>(d-1)/d\),
\begin{equation}\label{eq:ihyp-lower}
\liminf_{R\to\infty}
\frac{V_\eta(R)}{\Vol_d(B_R)^\beta}=0,
\end{equation}
and, for every \(\alpha<2\),
\begin{equation}\label{eq:ihyp-upper}
\limsup_{R\to\infty}
\frac{V_\eta(R)}{\Vol_d(B_R)^\alpha}=+\infty.
\end{equation}
Write \(\IHyp_\rho\) for this class.

\begin{maintheorem}[Hyperuniformity]
\label{thm:main-hyperuniformity}
The set \(\IHyp_\rho\) is residual in \(\smash{\EPP_{\rho,2}}\), while
\(\Hyp_\rho\cap\smash{\EPP_{\rho,2}}\) is dense and meagre in
\(\smash{\EPP_{\rho,2}}\).
\end{maintheorem}

Hyperuniform examples include determinantal processes with projection kernels,
such as the sine and Ginibre processes \cite{GL17}, the two-dimensional
one-component plasma \cite{Leb26}, certain cut-and-project quasicrystals
\cite{BH24}, and several transport constructions \cite{KLLY25}.
Stable matching with a lattice has also produced hyperuniform and
number-rigid processes \cite{KLY20}.

\begin{remark}
\label{rem:KLLY-density}
The fair-partition construction of \cite[Example~4.7]{KLLY25} produces,
from every ergodic finite-intensity point process, a hyperuniform process
satisfying
\(\Cov_\eta(N_A,N_B)\leq0\) for disjoint bounded Borel sets \(A,B\).
By \cref{cor:pnc-nowhere-dense,prop:periodic-density}, this class is
nowhere dense in \(\EPP_{\rho,2}\); the fair-partition construction
therefore does not prove density of hyperuniform processes.
\end{remark}

\subsection{Maximal rigidity}
\label{subsec:intro-rigidity}

For a Borel set \(B\subset\mathbb R^d\), write
\(\omega|_B:=\omega(\,\cdot\cap B)\). For bounded Borel \(A\), let
\[
\mathcal F_{A^c}:=\sigma(N_E:E\subset A^c\text{ bounded Borel}).
\]
The process \(\eta\in\PP_{\rho,2}\) is \emph{number rigid} if, for every
bounded Borel \(A\), the count \(N_A\) is
\(\mathcal F_{A^c}\)-measurable modulo \(\eta\)-null sets. It is
\emph{maximally rigid} if, for every such \(A\), the restriction
\(\omega|_A\) is \(\mathcal F_{A^c}\)-measurable modulo \(\eta\)-null sets.
Equivalently, every \(N_E\) with \(E\subset A\) bounded Borel is
\(\mathcal F_{A^c}\)-measurable. Write \(\MRig_\rho\subset\PP_{\rho,2}\)
for the class of maximally rigid processes.

We use real $L^2$ for rigidity and complex $L^2$ for spectral arguments.

We also consider approximation by exterior linear statistics.
For a bounded, compactly supported Borel function \(f\) and
\(\omega\in\operatorname{Conf}(\mathbb R^d)\), set
\[
S_f(\omega):=\int_{\mathbb R^d}f\,d\omega,
\qquad
Z_f^\eta:=S_f-\rho\int_{\mathbb R^d}f.
\]
For each bounded Borel set $A$, let $\mathcal H_{A^c}^\eta$ be the
closed span in real $L^2(\eta)$ of $Z_g^\eta$, over all bounded,
compactly supported Borel functions $g$ vanishing on $A$.
We call $\eta$ \emph{linearly maximally rigid} if
\[
Z_f^\eta\in\mathcal H_{A^c}^\eta
\]
for every bounded Borel $A$ and every bounded Borel function $f$
vanishing outside $A$. Write \(\LMRig_\rho\subset\PP_{\rho,2}\) for
this class. Since
\(\mathcal H_{A^c}^\eta\subset L^2(\mathcal F_{A^c},\eta)\), linear maximal
rigidity implies maximal rigidity. Indeed, for every bounded Borel
\(E\subset A\), the random variable \(N_E-\rho\Vol_d(E)\) is
\(\mathcal F_{A^c}\)-measurable.

\begin{maintheorem}[Rigidity]
\label{thm:main-rigidity}
The set
\[
\LMRig_\rho\cap\smash{\EPP_{\rho,2}}
\]
is residual in \(\smash{\EPP_{\rho,2}}\).  In particular,
\(\MRig_\rho\cap\smash{\EPP_{\rho,2}}\) is residual.
\end{maintheorem}

\begin{corollary}
\label{cor:main-intersection}
The sets
\[
\IHyp_\rho\cap\LMRig_\rho\cap\smash{\EPP_{\rho,2}}
\]
and
\[
(\MRig_\rho\setminus\Hyp_\rho)\cap\smash{\EPP_{\rho,2}}
\]
are residual in \(\smash{\EPP_{\rho,2}}\).
\end{corollary}

In Appendix~\ref{app:rotation-examples}, we construct a one-dimensional
family of linearly maximally rigid point processes with uniformly discrete
and relatively dense configurations. Their ball-count variance is bounded
along a sequence of radii. The processes are hyperuniform for
Lebesgue-almost every parameter, but intermittently hyperuniform for
a dense $G_\delta$-set of Liouville parameters. Here the
full-measure and residual statements refer to the parameter interval.

The distinction between number rigidity and maximal rigidity is
visible in the Ginibre process and the zero process of the planar Gaussian
analytic function: their exterior determines certain interior statistics,
but not all interior point positions \cite{GP17}. Through
uniqueness-pair methods, maximal rigidity also holds for further
examples including quasicrystals \cite{LR26}. Spectral criteria for finite-order rigidity provide another
approach \cite{LRrig24}. Deletion and insertion tolerance provide
related distinctions \cite{HS13}.

Baire-category methods for measure-preserving transformations and actions
originate in the work of Halmos and Rokhlin and are developed, among other
places, in \cite{Hal44,Rok48,GK98,GW19,FRW11,Kec10}.  We use
related approximation methods for translation-invariant point processes,
with intensity and local second moments fixed.

\subsection{Fixed actions and the local-to-global principle}
\label{subsec:intro-fixed-actions}

Let \((X,\mu)\) be a nonatomic standard probability space and
let
\[
a:\mathbb R^d\times X\longrightarrow X
\]
be an essentially free ergodic probability-measure-preserving (p.m.p.)
action. Write \(t.x:=a(t,x)\). A \emph{realization} of \(a\) is a measurable
map \(\kappa:X\to\operatorname{Conf}(\mathbb R^d)\) satisfying
\(\kappa(t.x)=\theta_t\kappa(x)\) almost everywhere. It is \emph{generating}
if the \(\sigma\)-algebra generated by \(\kappa\) agrees with the full
\(\sigma\)-algebra modulo null sets.

We denote by \(\RR^a_{\rho,2}(X,\mu)\) the Polish space of realizations of
intensity \(\rho\) with local second moments, in the rooted local
\(L^2\) topology, and by \(\RR^{a,\mathrm{iso}}_{\rho,2}(X,\mu)\) its
\(G_\delta\)-subspace of generating realizations. The rooted topology is
defined in \cref{subsec:configuration-maps,subsec:rooted}. The push-forward map
\[
\RR^a_{\rho,2}(X,\mu)\longrightarrow\PP_{\rho,2},
\qquad
\kappa\longmapsto\kappa_*\mu,
\]
is continuous by \cref{prop:global-polish}.

We first prove a fixed-action result from which
\cref{thm:main-hyperuniformity,thm:main-rigidity} follow.

\begin{proposition}[Fixed-action residuality]
\label{prop:main-fixed-action}
For every essentially free ergodic p.m.p. action
$a:\mathbb R^d\times X\to X$, the set
\[
\left\{
\kappa\in\RR^{a,\mathrm{iso}}_{\rho,2}(X,\mu):
\kappa_*\mu\in\IHyp_\rho\cap\LMRig_\rho
\right\}
\]
is residual in \(\smash{\RR^{a,\mathrm{iso}}_{\rho,2}(X,\mu)}\).
\end{proposition}

The essentially free processes form a dense \(G_\delta\)-subset
\(\EPP^{\mathrm{fr}}_{\rho,2}\) of \(\EPP_{\rho,2}\), by
\cref{prop:free-dense}. The local-to-global principle,
\cref{prop:local-global}, says that a \(G_\delta\)-class
\(\mathcal A\subset\EPP_{\rho,2}\) is residual whenever its pullback is
dense in \(\RR^{a,\mathrm{iso}}_{\rho,2}(X,\mu)\) for every essentially free
ergodic action \(a\). With \cref{prop:main-fixed-action}, this yields
\cref{thm:main-hyperuniformity,thm:main-rigidity,cor:main-intersection}.

The constructions in \cite{Bjo26b} give maximally rigid Delone
(uniformly discrete and relatively dense) realizations of arbitrary
essentially free ergodic actions at sufficiently large prescribed intensity.

\subsection{Fluctuations along van Hove sequences}
\label{subsec:intro-van-hove}

Let \((F_n)\) be bounded Borel sets satisfying
\(0<\Vol_d(F_n)<\infty\). We call \((F_n)\) a \emph{F\o lner sequence} if,
for every \(t\in\mathbb R^d\),
\[
\frac{\Vol_d((F_n+t)\mathbin{\triangle}F_n)}{\Vol_d(F_n)}
\longrightarrow0\qquad(n\to\infty).
\]
For $r>0$, define
\[
\partial_r F:=\{x\in\mathbb R^d:\operatorname{dist}(x,\partial F)<r\}.
\]
The sequence is \emph{van Hove} if
\[
\frac{\Vol_d(\partial_rF_n)}{\Vol_d(F_n)}\longrightarrow0
\qquad(n\to\infty,\ \text{for every }r>0).
\]
If \(K\subset\mathbb R^d\) is a compact convex body with nonempty
interior and \(R_n\to\infty\), then \((R_nK)\) is van Hove. This includes
Euclidean balls and closed \(\ell^p\)-balls, \(1\leq p\leq\infty\).
Every van Hove sequence is F\o lner. For every
\(\eta\in\smash{\EPP_{\rho,2}}\) and every van Hove sequence \((F_n)\),
\[
\Var_\eta(N_{F_n})=o\!\left(\Vol_d(F_n)^2\right);
\]
see \cref{prop:van-hove-subquadratic}. The bound can fail for F\o lner
sequences that are not van Hove; see \cref{rem:folner-versus-van-hove}.

For \(\eta\in\smash{\EPP_{\rho,2}}\), we call a sequence
\((t_n)\subset\mathbb R^d\) a \emph{revival sequence} if \(|t_n|\to\infty\)
and
\[
\bigl\|N_{A+t_n}-N_A\bigr\|_{L^2(\eta)}\longrightarrow0
\qquad(n\to\infty)
\]
for every bounded Borel set \(A\subset\mathbb R^d\) with
\(\Vol_d(\partial A)=0\). Let \(\Rev_\rho\) be the class of processes
admitting a revival sequence. Every subsequence of a revival sequence is
again a revival sequence.

For \(\eta\in\smash{\EPP_{\rho,2}}\), the property of
\emph{near-quadratic van Hove fluctuations} means that
\[
\begin{gathered}
\text{for every }\alpha_n>0\text{ with }\alpha_n\to0
\text{ and every }M_n>0\text{ with }M_n\to\infty,\\
\text{there is a van Hove sequence }(F_n)\text{ with }
\Var_\eta(N_{F_n})\geq
M_n\Vol_d(F_n)^{2-\alpha_n}\quad(n\geq1).
\end{gathered}
\]
Let \(\MFVH_\rho\) denote the class of processes with this property.

\begin{maintheorem}[Revival and near-quadratic fluctuations]
\label{thm:main-mfvh}
The set \(\Rev_\rho\) is residual in \(\smash{\EPP_{\rho,2}}\) and
\(\Rev_\rho\subseteq\MFVH_\rho\).  Consequently \(\MFVH_\rho\) is
residual in \(\smash{\EPP_{\rho,2}}\).
\end{maintheorem}

By \cref{prop:revival-koopman,prop:rigidity-implies-linear},
every process in $\Rev_\rho$ is linearly maximally rigid.

\subsection{Organization}
\label{subsec:organization}

Section~\ref{sec:spaces} establishes the Polish spaces and the
local-to-global principle.  Section~\ref{sec:intermittent-hyperuniformity}
proves the residuality and meagreness assertions of
\cref{thm:main-hyperuniformity}; its density assertion and the final
proof appear in \cref{subsec:periodic-approximation}.
The proofs of \cref{prop:main-fixed-action,thm:main-rigidity,cor:main-intersection}
appear at the end of \cref{subsec:linear-maximal-rigidity}.
\Cref{thm:main-mfvh} follows from
\cref{prop:revival-residual,prop:revival-mfvh} in
Section~\ref{sec:revival}.  Appendix~\ref{subsec:generic-rigidity}
treats rigidity sequences for actions, and
Appendix~\ref{app:rotation-examples} gives intermittently hyperuniform point processes.

\section{Point-process spaces and fixed actions}\label{sec:spaces}

\subsection{Configuration maps and cross-sections}\label{subsec:configuration-maps}

Let $\mathcal M_+(\mathbb R^d)$ denote the space of locally finite positive
Radon measures on $\mathbb R^d$.  For $\varphi\in C_c(\mathbb R^d)$ and
$\omega\in\mathcal M_+(\mathbb R^d)$, set
\[
S_\varphi(\omega):=\int_{\mathbb R^d}\varphi\,d\omega.
\]
The \emph{vague topology} on $\mathcal M_+(\mathbb R^d)$ is the
weakest topology making all maps $S_\varphi$,
$\varphi\in C_c(\mathbb R^d)$, continuous.

Let $\mathcal N(\mathbb R^d)$ be the subspace of locally finite counting
measures, allowing multiplicities. For a locally finite set
$P\subset\mathbb R^d$, write
\[
\delta_P:=\sum_{p\in P}\delta_p.
\]
The space of locally finite simple counting measures is
\[
\operatorname{Conf}(\mathbb R^d)
:=\{\delta_P:P\subset\mathbb R^d\text{ is locally finite}\}.
\]  With the vague
topology, $\mathcal N(\mathbb R^d)$ is Polish and its Borel
$\sigma$-algebra is generated by the evaluation maps
\[
N_A(\omega):=\omega(A),
\qquad A\subset\mathbb R^d\text{ bounded and Borel};
\]
see \cite[Proposition~9.1.IV]{DVJ08}.  The subspace
$\operatorname{Conf}(\mathbb R^d)$ is a $G_\delta$ subset of
$\mathcal N(\mathbb R^d)$, hence is itself Polish; see
\cite{Pug02}. Weak convergence of point-process laws will always refer
to this vague topology unless stated otherwise.

We use \(\omega\) for the counting measure and
\(\operatorname{supp}\omega\) for its underlying point set. Thus sums
indexed by points of a configuration are written as
\(\sum_{u\in\operatorname{supp}\omega}\).

A point-process law \(\eta\) on \(\operatorname{Conf}(\mathbb R^d)\)
has \emph{finite local second moments} if
\begin{equation}\label{eq:finite-local-second-moments}
\int N_A(\omega)^2\,d\eta(\omega)<\infty
\qquad\text{for every bounded Borel }A\subset\mathbb R^d.
\end{equation}
Equivalently, it suffices to check the cubes
\(Q_m:=[-m/2,m/2)^d\), \(m\geq1\).

Recall that, for $t\in\mathbb R^d$,
\[
(\theta_t\omega)(A):=\omega(A+t).
\]
Then $\operatorname{supp}(\theta_t\omega)=\operatorname{supp}(\omega)-t$.
The map $(t,\omega)\mapsto\theta_t\omega$ is continuous. If
$t_n\to t$ and $\omega_n\to\omega$ vaguely, then
$\varphi(\,\cdot-t_n)\to\varphi(\,\cdot-t)$ uniformly for
every $\varphi\in C_c(\mathbb R^d)$. These functions have support
in a common compact set, whose $\omega_n$-masses are bounded.
Consequently,
$S_\varphi(\theta_{t_n}\omega_n)\to
S_\varphi(\theta_t\omega)$.

Throughout, $(X,\mu)$ is a nonatomic standard probability space.  Let
$\Ad_d(X,\mu)$ denote the essentially free ergodic p.m.p. actions
\[
a:\mathbb R^d\times X\longrightarrow X.
\]
For $a\in\Ad_d(X,\mu)$, set $t.x:=a(t,x)$ and denote by
\[
L^0\!\left(X,\operatorname{Conf}(\mathbb R^d),\mu\right)
\]
the space of equivalence classes, modulo $\mu$-almost everywhere equality,
of measurable maps from $X$ to $\operatorname{Conf}(\mathbb R^d)$.  A class
$\kappa$ is \emph{equivariant} if, for every $t\in\mathbb R^d$,
\[
\kappa(t.x)=\theta_t\kappa(x)
\qquad\text{for $\mu$-almost every }x\in X.
\]
The exceptional null set may depend on $t$. For each fixed $t$,
equivariance and invariance of $\mu$ give
\[
(\theta_t)_*(\kappa_*\mu)
=(\kappa(t.\,\cdot))_*\mu
=\kappa_*\mu.
\]
\begin{lemma}\label{lem:concrete-version}
Let $\kappa\in L^0(X,\operatorname{Conf}(\mathbb R^d),\mu)$ be equivariant.
There exist an invariant conull Borel set $X_\kappa\subset X$ and a Borel
representative
\[
\kappa^\circ:X_\kappa\longrightarrow\operatorname{Conf}(\mathbb R^d)
\]
of $\kappa$ such that
\[
\kappa^\circ(t.x)=\theta_t\kappa^\circ(x)
\qquad
(t,x)\in\mathbb R^d\times X_\kappa.
\]
If $\kappa(x)$ is nonempty almost everywhere, then $X_\kappa$ may be chosen
so that $\kappa^\circ(x)$ is nonempty for every $x\in X_\kappa$.
\end{lemma}

\begin{proof}
Apply \cite[Proposition~5.7]{ABC25} to the measurable factor $\kappa$.
\end{proof}

Let $X_0\subset X$ be an invariant conull Borel set.  A Borel set
$Y\subset X_0$ is a \emph{cross-section} if, for every $x\in X_0$, the
return-time set
\[
Y_x:=\{t\in\mathbb R^d:t.x\in Y\}
\]
is nonempty and locally finite.

\begin{proposition}\label{prop:configuration-cross-section}
Let $\kappa\in L^0(X,\operatorname{Conf}(\mathbb R^d),\mu)$ be equivariant
and nonempty almost everywhere.  Then there exist an invariant conull Borel
set $X_\kappa\subset X$ and a nonempty Borel representative
$\kappa^\circ:X_\kappa\to\operatorname{Conf}(\mathbb R^d)$ satisfying
\[
\kappa^\circ(t.x)=\theta_t\kappa^\circ(x)
\qquad(t,x)\in\mathbb R^d\times X_\kappa,
\]
such that
\[
Y_\kappa
:=
\{x\in X_\kappa:N_{\{0\}}(\kappa^\circ(x))=1\}
\]
is a cross-section and
\[
(Y_\kappa)_x=\operatorname{supp}(\kappa^\circ(x)),
\qquad x\in X_\kappa.
\]
\end{proposition}

\begin{proof}
Choose $X_\kappa$ and $\kappa^\circ$ from \cref{lem:concrete-version}.
The set $Y_\kappa$ is Borel because $N_{\{0\}}$ is Borel on
$\operatorname{Conf}(\mathbb R^d)$.  For $t\in\mathbb R^d$ and
$x\in X_\kappa$, pointwise equivariance gives
\[
t\in(Y_\kappa)_x
\iff
0\in\operatorname{supp}(\kappa^\circ(t.x))
\iff
0\in\operatorname{supp}(\kappa^\circ(x))-t
\iff
t\in\operatorname{supp}(\kappa^\circ(x)).
\]
Since $\operatorname{supp}(\kappa^\circ(x))$ is nonempty and locally
finite, $Y_\kappa$ is a cross-section.
\end{proof}

Conversely, a cross-section $Y\subset X_0$ gives the equivariant Borel map
\[
\kappa_Y(x):=\sum_{t\in Y_x}\delta_t.
\]
The relation $\{(x,t):x\in X_0,\ t.x\in Y\}$ is Borel,
and its sections $Y_x$ are countable by local finiteness.
Lusin--Novikov \cite[Theorem~18.10]{Kec95} expresses this
relation as a union of graphs of partial Borel maps $t_j$.
For each bounded Borel $D$, the function
$x\mapsto N_D(\kappa_Y(x))=\sum_{j:x\in\operatorname{dom}(t_j)}
\mathbf1_D(t_j(x))$ is Borel.  The identity
$Y_{t.x}=Y_x-t$ gives $\kappa_Y(t.x)=\theta_t\kappa_Y(x)$.

On an invariant conull Borel set $X_0$ where $a$ is free,
\cite[Theorem~A.3]{CGMTD26}, applied with the compact set
$\overline B_{r/2}$, gives a cocompact Borel cross-section $C_r$
for every $r>0$ such that each return set $(C_r)_x$ is
$r$-separated. Its configuration map has uniformly bounded local
counts, but need not have intensity $\rho$.

\smallskip\noindent\emph{Voronoi cells of a cross-section.}
Fix an $r$-separated cocompact Borel cross-section
$C\subset X_0$. Cocompactness gives $R<\infty$ such that
$C_x$ is $R$-dense in $\mathbb R^d$ for every $x\in X_0$.
Choose a Borel injection $\ell:C\to[0,1]$, which exists by
\cite[Theorem~15.6]{Kec95}. For $c\in C$, define
\begin{equation}\label{eq:labelled-voronoi}
W_c:=\left\{t\in\mathbb R^d:
 (\|t\|,\ell(c))\leq_{\mathrm{lex}}
 (\|t-v\|,\ell(v.c))
 \text{ for every }v\in C_c\right\},
\end{equation}
where $\leq_{\mathrm{lex}}$ denotes the lexicographic order.
Thus the cell $W_c$ consists of the points for which $0$ is a nearest
return vector in $C_c$, with equal distances resolved by $\ell$.
In particular, $W_c$ depends on the chosen cross-section and its
Borel injection; we retain the notation $W_c$ when both are fixed.

The relation $\{(c,v)\in C\times\mathbb R^d:v.c\in C\}$ is
Borel and has countable sections, since $C_c$ is $r$-separated.
By Lusin--Novikov \cite[Theorem~18.10]{Kec95}, it is the union
of graphs of partial Borel maps $v_j:D_j\to\mathbb R^d$.
Applying the comparisons in \eqref{eq:labelled-voronoi} to these
maps proves that $\{(c,t):c\in C,\ t\in W_c\}$ is Borel.
Each comparison defines an open or closed half-space, so $W_c$
is convex. Separation and cocompactness give
\[
B_{r/2}^{\circ}\subset W_c\subset\overline B_R.
\]
For every $x\in X_0$, the sets
$u+W_{u.x}$, $u\in C_x$, partition $\mathbb R^d$:
each $t$ has a nearest return vector $u\in C_x$, and $\ell$
resolves ties between the finitely many nearest returns.
The partition is equivariant under translation and locally finite,
since all its cells lie within distance $R$ of their return vectors.

\subsection{Rooted realization spaces}\label{subsec:rooted}

Fix a complete compatible metric $d_{\mathrm{Conf}}$ on
$\operatorname{Conf}(\mathbb R^d)$, bounded by $1$.
Choose a countable dense family $(\omega_j)_{j\geq1}$ in this space,
and let $\mathcal D$ be the countable unital $\mathbb Q$-algebra
generated by the functions
$\omega\mapsto d_{\mathrm{Conf}}(\omega,\omega_j)$.
Every element of $\mathcal D$ is bounded and continuous.
For $\kappa,\lambda\in L^0(X,\operatorname{Conf}(\mathbb R^d),\mu)$, set
\[
d_0(\kappa,\lambda)
:=
\int_X d_{\mathrm{Conf}}(\kappa(x),\lambda(x))\,d\mu(x).
\]
Convergence in $d_0$ is convergence in measure, since for every
$\varepsilon>0$,
\[
\varepsilon\,\mu\{d_{\mathrm{Conf}}(\kappa,\lambda)>\varepsilon\}
\leq d_0(\kappa,\lambda)
\leq
\varepsilon+
\mu\{d_{\mathrm{Conf}}(\kappa,\lambda)>\varepsilon\}.
\]
Recall the cubes \(Q_m\) from \eqref{eq:finite-local-second-moments}.
An equivariant \(\kappa\in L^0(X,\operatorname{Conf}(\mathbb R^d),\mu)\)
has \emph{finite local second moments} if
\(N_{Q_m}\circ\kappa\in L^2(X,\mu)\) for every \(m\geq1\);
equivalently, its law \(\kappa_*\mu\) satisfies
\eqref{eq:finite-local-second-moments}.
It has \emph{intensity $\rho$} if
\[
\int_X N_A(\kappa(x))\,d\mu(x)
=
\rho\,\Vol_d(A)
\]
for every bounded Borel set $A\subset\mathbb R^d$.

Let $\smash{\RR^a_{\rho,2}(X,\mu)}$ be the set of $a$-equivariant classes
of intensity $\rho$ with finite local second moments,
equipped with the metric
\[
d_{\mathrm{root},2}(\kappa,\lambda)
:=
d_0(\kappa,\lambda)
+
\sum_{m=1}^\infty 2^{-m}
\min\!\left\{1,
\left\|N_{Q_m}\circ\kappa-N_{Q_m}\circ\lambda\right\|_2
\right\}.
\]
This metric defines the \emph{rooted local $L^2$ topology}.
A sequence $(\kappa_n)$ converges to $\kappa$ in
$\smash{\RR^a_{\rho,2}(X,\mu)}$ if and only if
$\kappa_n\to\kappa$ in measure and
\[
N_{Q_m}\circ\kappa_n\longrightarrow N_{Q_m}\circ\kappa
\quad\text{in }L^2(X,\mu)\qquad(n\to\infty,\ m\geq1).
\]

For $\kappa\in\smash{\RR^a_{\rho,2}(X,\mu)}$, let $\mathcal F_\kappa$ be the
$\mu$-completion of the $\sigma$-algebra generated by $\kappa$.  An element
$\kappa$ is \emph{generating} if $\mathcal F_\kappa=\mathcal B_X^\mu$.  Define
\[
\RR^{a,\mathrm{iso}}_{\rho,2}(X,\mu)
:=
\{\kappa\in\RR^a_{\rho,2}(X,\mu):
  \mathcal F_\kappa=\mathcal B_X^\mu\}.
\]
For such $\kappa$, the factor map $\kappa$ is a measurable isomorphism,
modulo null sets, from
$\mathbb R^d\overset{a}{\curvearrowright}(X,\mu)$ to the translation
action on $(\operatorname{Conf}(\mathbb R^d),\kappa_*\mu)$: on standard
Borel probability spaces, a generating factor admits a measurable inverse
modulo null sets by \cite[Theorem~15.1]{Kec95}.

\begin{lemma}\label{lem:null-window}
Let $\eta$ be a translation-invariant probability measure on
$\operatorname{Conf}(\mathbb R^d)$.  If $A\subset\mathbb R^d$ is bounded and Borel, with
$\Vol_d(A)=0$, then $N_A=0$ for $\eta$-almost every configuration.
\end{lemma}

\begin{proof}
Fix a bounded Borel set $B\subset\mathbb R^d$ with positive volume.  Tonelli's
theorem gives
\[
\int_B\int N_{A+t}(\omega)\,d\eta(\omega)\,dt
=
\int\sum_{p\in\operatorname{supp}(\omega)}
\Vol_d\bigl(B\cap(p-A)\bigr)\,d\eta(\omega)
=0.
\]
Translation invariance makes the inner expectation on the left independent of $t$, so
$\smash{\int} N_A\,d\eta=0$.
\end{proof}

\begin{lemma}[Vitali convergence]\label{lem:vitali-lp}
Let \((f_n)\) be measurable functions on a probability space and let
\(1\leq p<\infty\).  If \(f_n\to f\) in measure and
\(\{|f_n|^p:n\geq1\}\) is uniformly integrable, then
\(f\in L^p\) and \(f_n\to f\) in \(L^p\).
\end{lemma}

\begin{proof}
Every subsequence has a further subsequence converging almost everywhere to
\(f\).  Fatou gives \(f\in L^p\), while
\(|f_n-f|^p\leq2^{p-1}(|f_n|^p+|f|^p)\) makes the differences uniformly
integrable.  Vitali's theorem gives \(L^p\)-convergence along the further
subsequence, hence along the original sequence.
\end{proof}

\begin{proposition}[Polish structure of realization spaces]\label{prop:rooted-polish}
The space $\smash{\RR^a_{\rho,2}(X,\mu)}$ is Polish.  Moreover:
\begin{enumerate}
\item if $A\subset\mathbb R^d$ is bounded and Borel, with $\Vol_d(\partial A)=0$, then
$\kappa\mapsto N_A\circ\kappa$ is continuous into $L^2(X,\mu)$;
\item if $\varphi\in C_c(\mathbb R^d)$, then
$\kappa\mapsto S_\varphi\circ\kappa$ is continuous into $L^2(X,\mu)$;
\item $\smash{\RR^{a,\mathrm{iso}}_{\rho,2}(X,\mu)}$ is a $G_\delta$ subset of
$\smash{\RR^a_{\rho,2}(X,\mu)}$, hence Polish.
\end{enumerate}
\end{proposition}

\begin{proof}
The ambient space
$L^0(X,\operatorname{Conf}(\mathbb R^d),\mu)$ is Polish for convergence in
measure.  A $d_0$-Cauchy sequence has a subsequence whose
successive $d_0$-distances are summable; Tonelli and completeness of
$d_{\mathrm{Conf}}$ give an almost-everywhere pointwise limit, and hence a
$d_0$-limit.  Separability follows after identifying $(X,\mu)$ with
$([0,1],\operatorname{Leb})$ modulo null sets
\cite[Theorem~17.41]{Kec95}: step functions with rational intervals and values
in a fixed countable dense subset of $\operatorname{Conf}(\mathbb R^d)$ are
dense.

Consider
\[
J:\RR^a_{\rho,2}(X,\mu)
\longrightarrow
L^0(X,\operatorname{Conf}(\mathbb R^d),\mu)
\times\prod_{m\geq1}L^2(X,\mu),
\qquad
J(\kappa):=
\bigl(\kappa,(N_{Q_m}\circ\kappa)_{m\geq1}\bigr).
\]
On the product we use $d_0$ on the first coordinate and
\(\smash{\sum_{m\geq1}}2^{-m}\min\{1,\|g_m-h_m\|_2\}\) on the $L^2$ coordinates.
With these metrics, $J$ is isometric. Let
$J(\kappa_n)\to(\kappa,(g_m)_{m\geq1})$ in the product.  For each fixed
$t\in\mathbb R^d$, composition with $x\mapsto t.x$ is an isometry of
$L^0$, while composition with $\theta_t$ is continuous in measure because
$\theta_t$ is continuous on $\operatorname{Conf}(\mathbb R^d)$.  Passing to
the limit in
$\kappa_n(t.x)=\theta_t\kappa_n(x)$ shows that $\kappa$ is equivariant.
Hence $\kappa_*\mu$ is translation-invariant.

After passing to a subsequence, still denoted $(\kappa_n)$, we may assume
$\smash{\sum_n} d_0(\kappa_n,\kappa)<\infty$.  Tonelli's theorem then gives
$\kappa_n(x)\to\kappa(x)$ vaguely for almost every $x$.  By \cref{lem:null-window}, $\kappa(x)(\partial Q_m)=0$ for every $m$ outside
one null set.  For such an $x$, the continuity-set criterion for vague convergence
\cite[Appendix~A2.6]{DVJ03} gives
\[
N_{Q_m}(\kappa_n(x))\longrightarrow N_{Q_m}(\kappa(x))
\qquad(n\to\infty).
\]
Since
$N_{Q_m}\circ\kappa_n\to g_m$ in $L^2$ for every $m$, a diagonal
subsequence gives $g_m=N_{Q_m}\circ\kappa$ almost everywhere for all $m$.
Hence $N_{Q_m}\circ\kappa\in L^2(X,\mu)$ for every $m$.

Let $A\subset\mathbb R^d$ be a bounded Borel set with
$\Vol_d(\partial A)=0$, and choose $m$ with $A\subset Q_m$.  By \cref{lem:null-window} and the continuity-set criterion for vague
convergence, $N_A\circ\kappa_n\to N_A\circ\kappa$ in measure.  Since
$N_{Q_m}\circ\kappa_n\to N_{Q_m}\circ\kappa$ in $L^2$, the
squares of these local counts are uniformly integrable.  Since $N_A\leq N_{Q_m}$, the variables $N_A\circ\kappa_n$ are
uniformly integrable.  Hence
\[
\int_X N_A(\kappa(x))\,d\mu(x)
=
\lim_{n\to\infty}\int_X N_A(\kappa_n(x))\,d\mu(x)
=
\rho\,\Vol_d(A).
\]
Uniqueness of locally finite measures extends the identity from rational
half-open rectangles to all bounded Borel sets.  Hence $\kappa$ has intensity $\rho$.

Thus $J(\RR^a_{\rho,2}(X,\mu))$ is closed in a Polish space,
and $\RR^a_{\rho,2}(X,\mu)$ is Polish.

\smallskip\noindent\emph{Continuity of local statistics.}
Let $\kappa_n\to\kappa$ in $\smash{\RR^a_{\rho,2}(X,\mu)}$, let
$A\subset\mathbb R^d$ be bounded and Borel, with $\Vol_d(\partial A)=0$, and choose
$m$ with $A\subset Q_m$.  By \cref{lem:null-window}, $\kappa(x)(\partial A)=0$ for almost every $x$.
Vague convergence in measure therefore gives
$N_A\circ\kappa_n\to N_A\circ\kappa$ in measure, while
\[
|N_A\circ\kappa_n-N_A\circ\kappa|^2
\leq
2(N_{Q_m}\circ\kappa_n)^2+2(N_{Q_m}\circ\kappa)^2.
\]
Since $N_{Q_m}\circ\kappa_n\to N_{Q_m}\circ\kappa$ in $L^2$, the squares on
the right are uniformly integrable.  Hence
\cref{lem:vitali-lp} gives
$N_A\circ\kappa_n\to N_A\circ\kappa$ in $L^2$.

Let $\varphi\in C_c(\mathbb R^d)$ and choose $m$ with
$\operatorname{supp}\varphi\subset Q_m$.  Vague convergence in measure gives
$S_\varphi\circ\kappa_n\to S_\varphi\circ\kappa$ in measure, while
\[
|S_\varphi\circ\kappa_n-S_\varphi\circ\kappa|^2
\leq
2\|\varphi\|_\infty^2
\bigl((N_{Q_m}\circ\kappa_n)^2+(N_{Q_m}\circ\kappa)^2\bigr).
\]
By \cref{lem:vitali-lp}, the linear statistics converge in $L^2$.

\smallskip\noindent\emph{Generating realizations.}
The algebra $\mathcal D$ defined above separates points and generates
the Borel $\sigma$-algebra. Hence the functional monotone-class theorem
shows that $\mathcal D$ is dense in real $L^2(\eta)$ for every
probability measure $\eta$ on $\operatorname{Conf}(\mathbb R^d)$.

Let $(u_j)_{j\geq1}$ be a dense family of bounded real-valued
functions in real $L^2(X,\mu)$.  The map $\kappa$ is generating if and only if, for every
$j,n\geq1$, there exists $h\in\mathcal D$ such that
\[
\|u_j-h\circ\kappa\|_2<\frac1n.
\]
For a fixed bounded continuous $h$, convergence in $d_0$ gives
$h\circ\kappa_n\to h\circ\kappa$ in measure and, by boundedness,
in $L^2$. Thus each set
\[
\left\{\kappa\in\RR^a_{\rho,2}(X,\mu):
\|u_j-h\circ\kappa\|_2<\frac1n\right\}
\]
is open.  Therefore
\[
\RR^{a,\mathrm{iso}}_{\rho,2}(X,\mu)
=
\bigcap_{j,n\geq1}
\bigcup_{h\in\mathcal D}
\left\{\kappa:
\|u_j-h\circ\kappa\|_2<\frac1n\right\}
\]
is $G_\delta$.
\end{proof}

\subsection{Global point-process spaces}\label{subsec:global}

For a Polish space $Y$ and a continuous function
$\Psi:Y\to[1,\infty)$, let $\mathscr P^\Psi(Y)$ be the space of
probability measures $\eta$ on $Y$ with $\smash{\int}\Psi\,d\eta<\infty$.
We equip it with the \emph{$\Psi$-weak topology}, defined by convergence
of the integrals of all continuous functions $F:Y\to\mathbb R$ with
$|F|\leq C\Psi$ for some $C<\infty$.
The map
\[
\eta\longmapsto\frac{\Psi\eta}{\int\Psi\,d\eta}
\]
is a homeomorphism onto the space of probability measures on $Y$ with
the weak topology; its inverse sends $\nu$ to
$\Psi^{-1}\nu/\int\Psi^{-1}\,d\nu$. Thus $\mathscr P^\Psi(Y)$
is Polish. Equivalently, $\eta_n\to\eta$ in the $\Psi$-weak topology
if and only if $\eta_n\Rightarrow\eta$ and
$\smash{\int}\Psi\,d\eta_n\to\smash{\int}\Psi\,d\eta$
as $n\to\infty$.

For each $m\geq1$, choose $\varphi_m\in C_c(\mathbb R^d)$ such that
$0\leq\varphi_m\leq1$, $\varphi_m=1$ on $Q_m$, and
$\operatorname{supp}\varphi_m\subset Q_{m+1}$. Define
\[
\Psi_m(\omega):=1+S_{\varphi_m}(\omega)^2
\]
and abbreviate
$\mathscr P^{\Psi_m}:=\mathscr P^{\Psi_m}(\operatorname{Conf}(\mathbb R^d))$.
Each $\Psi_m$ is continuous. Let
\[
\mathscr P_{\mathrm{loc},2}:=\bigcap_{m\geq1}\mathscr P^{\Psi_m},
\]
with the topology induced by the diagonal map into
$\smash{\prod_{m\geq1}}\mathscr P^{\Psi_m}$.
Thus $\eta_n\to\eta$ in $\mathscr P_{\mathrm{loc},2}$ if and only if
$\eta_n\Rightarrow\eta$ and
\[
\int S_{\varphi_m}^2\,d\eta_n
\longrightarrow\int S_{\varphi_m}^2\,d\eta
\qquad(n\to\infty,\ m\geq1).
\]
Since $N_{Q_m}\leq S_{\varphi_m}\leq N_{Q_{m+1}}$, the condition $\eta\in\mathscr P_{\mathrm{loc},2}$ is equivalent to
$N_{Q_m}\in L^2(\eta)$ for every $m\geq1$.

Define $\PP_{\rho,2}$ to be the set of translation-invariant Borel
probability measures $\eta$ on $\operatorname{Conf}(\mathbb R^d)$
with intensity $\rho$ and finite local second moments, in the sense of
\eqref{eq:finite-local-second-moments}. Equivalently,
$\PP_{\rho,2}$ is the translation-invariant intensity-$\rho$ subspace of
$\mathscr P_{\mathrm{loc},2}$.  We equip it with the subspace topology,
called the \emph{local $L^2$ topology}.  Let
$\smash{\EPP_{\rho,2}}\subset\PP_{\rho,2}$ denote the ergodic elements.

\begin{proposition}\label{prop:global-polish}
The space $\PP_{\rho,2}$ is Polish for the local $L^2$ topology, and
$\EPP_{\rho,2}$ is a $G_\delta$ subset, hence Polish.  Moreover:
\begin{enumerate}
\item if $A\subset\mathbb R^d$ is bounded and Borel, with $\Vol_d(\partial A)=0$, then
$\eta\mapsto\smash{\int} N_A^2\,d\eta$ and $\eta\mapsto\Var_\eta(N_A)$ are
continuous on $\PP_{\rho,2}$;
\item for every $a\in\Ad_d(X,\mu)$, the map
$\kappa\mapsto\kappa_*\mu$ from $\smash{\RR^a_{\rho,2}(X,\mu)}$ to
$\PP_{\rho,2}$ is continuous.
\end{enumerate}
\end{proposition}

\begin{proof}
The diagonal embedding of $\mathscr P_{\mathrm{loc},2}$ into
$\prod_{m\geq1}\mathscr P^{\Psi_m}$ has closed image: each
$\Psi_m$-weak topology refines weak convergence, so all coordinates of
a limit agree. Thus $\mathscr P_{\mathrm{loc},2}$ is Polish.

The set of translation-invariant probability measures is closed: if
$\eta_n\to\eta$ in
$\mathscr P_{\mathrm{loc},2}$ and each $\eta_n$ is translation-invariant, then for
$t\in\mathbb R^d$ and $F\in C_b(\operatorname{Conf}(\mathbb R^d))$,
continuity of $\theta_t$ gives
\[
\int F\circ\theta_t\,d\eta
=
\lim_{n\to\infty}\int F\circ\theta_t\,d\eta_n
=
\lim_{n\to\infty}\int F\,d\eta_n
=
\int F\,d\eta.
\]

Suppose in addition that every $\eta_n$ has intensity $\rho$.  Fix a bounded
Borel set $A$ with $\Vol_d(\partial A)=0$ and choose $m$ with $A\subset Q_m$.
By translation invariance of $\eta$ and \cref{lem:null-window},
$\eta\{\omega:N_{\partial A}(\omega)=0\}=1$.  The continuity-set criterion
for vague convergence \cite[Appendix~A2.6]{DVJ03} therefore shows that the
distribution of $N_A$ under $\eta_n$ converges weakly on $\mathbb R$ to its distribution under
$\eta$.
We also have $\smash{\sup_n}\smash{\int} S_{\varphi_m}^2\,d\eta_n<\infty$.
By Cauchy--Schwarz and the bound
$\sup_{n\geq1}\int S_{\varphi_m}^2\,d\eta_n<\infty$,
\[
\lim_{R\to\infty}\sup_{n\geq1}
\int S_{\varphi_m}\mathbf1_{\{S_{\varphi_m}>R\}}\,d\eta_n=0.
\]
Since $N_A\leq S_{\varphi_m}$, the same tail bound holds for $N_A$.
Weak convergence of the distributions of $N_A$ and this bound give
\[
\int N_A\,d\eta
=
\lim_{n\to\infty}\int N_A\,d\eta_n
=
\rho\,\Vol_d(A).
\]
On each $Q_k$, the two first-moment measures are finite and agree on
rational half-open rectangles.  The monotone-class theorem extends the
identity to all bounded Borel sets.  Thus
$\PP_{\rho,2}$ is closed in $\mathscr P_{\mathrm{loc},2}$ and is Polish.

Enumerate the countable algebra $\mathcal D$ defined in
\cref{subsec:rooted} as $(h_j)_{j\geq1}$.  For $n\geq1$, put
\[
A_nh(\omega):=\frac{1}{\Vol_d(Q_n)}
\int_{Q_n}h(\theta_t\omega)\,dt,
\qquad
F_{j,n}(\eta):=\int\left|A_nh_j-\int h_j\,d\eta\right|^2d\eta.
\]
For fixed $j,n$, the function $A_nh_j$ is bounded and continuous, so
$F_{j,n}$ is weakly continuous in $\eta$.  For invariant $\eta$, let
$P_\eta$ be the orthogonal projection onto the invariant functions in
$L^2(\eta)$.  The mean ergodic theorem gives
\[
\lim_{n\to\infty}F_{j,n}(\eta)
=\left\|P_\eta h_j-\int h_j\,d\eta\right\|_{L^2(\eta)}^2.
\]  Since $\mathcal D$ is dense in $L^2(\eta)$ for every
$\eta$, this limit vanishes for all $j$ exactly when $\eta$ is ergodic.
Consequently, among invariant probability measures,
\[
\{\eta:\eta\text{ is ergodic}\}
=\bigcap_{j,k,N\geq1}\bigcup_{n\geq N}
\{\eta:F_{j,n}(\eta)<k^{-1}\}.
\]
The set of invariant laws is closed, and the displayed sets are open
in the weak topology.  Hence the ergodic laws form a $G_\delta$ in
the invariant-measure space and in the finer local $L^2$ topology on
$\PP_{\rho,2}$.

Let $\eta_n\to\eta$ in $\PP_{\rho,2}$, let $A\subset\mathbb R^d$ be
bounded and Borel, with $\Vol_d(\partial A)=0$, and choose $m$ with $A\subset Q_m$.
By \cref{lem:null-window} and the continuity-set criterion, the distribution of
$N_A$ under $\eta_n$ converges weakly to its distribution under $\eta$. Moreover,
\[
\int S_{\varphi_m}^2\,d\eta_n
\longrightarrow
\int S_{\varphi_m}^2\,d\eta.
\]
Since $S_{\varphi_m}^2$ is continuous and nonnegative, weak convergence
and convergence of its expectations imply, by applying weak convergence
to $\min\{S_{\varphi_m}^2,R\}$ for each fixed $R>0$, that
\[
\lim_{R\to\infty}\sup_n
\int (S_{\varphi_m}^2-R)_+\,d\eta_n=0.
\]
Since $u\mathbf1_{\{u>2R\}}\leq 2(u-R)_+$ for $u\geq0$, this implies
\[
\lim_{R\to\infty}\sup_n
\int S_{\varphi_m}^2
\mathbf1_{\{S_{\varphi_m}^2>R\}}\,d\eta_n=0.
\]
Since $N_A^2\leq S_{\varphi_m}^2$, its corresponding tail integrals
also tend to zero uniformly in $n$. Weak convergence of the distributions
of $N_A$ therefore yields
\[
\int N_A^2\,d\eta_n\longrightarrow\int N_A^2\,d\eta.
\]
Since every process in $\PP_{\rho,2}$ satisfies
$\smash{\int} N_A\,d\eta=\rho\Vol_d(A)$, we also have
\[
\Var_{\eta_n}(N_A)
=\int N_A^2\,d\eta_n-\rho^2\Vol_d(A)^2
\longrightarrow
\Var_{\eta}(N_A).
\]

Let $\kappa_n\to\kappa$ in $\smash{\RR^a_{\rho,2}(X,\mu)}$.  For every
$F\in C_b(\operatorname{Conf}(\mathbb R^d))$,
$F\circ\kappa_n\to F\circ\kappa$ in measure and, because $F$ is
bounded, in $L^1(X,\mu)$, so
$(\kappa_n)_*\mu\Rightarrow\kappa_*\mu$.  By
\cref{prop:rooted-polish},
$S_{\varphi_m}\circ\kappa_n\to S_{\varphi_m}\circ\kappa$ in $L^2$ for
every $m$.  Therefore
\[
\int S_{\varphi_m}^2\,d((\kappa_n)_*\mu)
\longrightarrow
\int S_{\varphi_m}^2\,d(\kappa_*\mu).
\]
This is convergence of $(\kappa_n)_*\mu$ to $\kappa_*\mu$
in the local $L^2$ topology.
\end{proof}

\begin{proposition}[Characterization of local $L^2$ convergence]
\label{prop:local-L2-characterization}
Let $\eta,\eta_1,\eta_2,\ldots\in\PP_{\rho,2}$.
\[
\eta_n\longrightarrow\eta
\quad\Longleftrightarrow\quad
\eta_n\Rightarrow\eta
\ \text{and}\ 
\int N_{Q_m}^2\,d\eta_n
\longrightarrow
\int N_{Q_m}^2\,d\eta
\quad(m\geq1).
\]
In particular, the topology does not depend on the cutoffs $(\varphi_m)$.
\end{proposition}

\begin{proof}
Assume first that $\eta_n\to\eta$ in the local $L^2$ topology.  Fix $m$.
Weak convergence and translation invariance imply that the distribution of
$N_{Q_m}$ under $\eta_n$ converges weakly on $\mathbb R$ to its distribution under
$\eta$, since $\Vol_d(\partial Q_m)=0$.  Choose $j$ with
$Q_m\subset Q_j$ and $\varphi_j\equiv1$ on $Q_m$.  Convergence of the
$S_{\varphi_j}$-second moments gives
\[
\lim_{R\to\infty}\sup_n
\int S_{\varphi_j}^2
\mathbf1_{\{S_{\varphi_j}^2>R\}}\,d\eta_n=0.
\]
Since $N_{Q_m}^2\leq S_{\varphi_j}^2$, the same tail bound holds for
$N_{Q_m}^2$.  Hence
\[
\int N_{Q_m}^2\,d\eta_n
\longrightarrow
\int N_{Q_m}^2\,d\eta.
\]

Conversely, suppose that $\eta_n\Rightarrow\eta$ and that the displayed
second-moment convergence holds for every $m$.  Since
$\operatorname{supp}(\varphi_m)\subset Q_{m+1}$,
\[
S_{\varphi_m}^2\leq N_{Q_{m+1}}^2.
\]
Convergence of the $Q_{m+1}$-second moments and weak convergence give
\[
\lim_{R\to\infty}\sup_n
\int N_{Q_{m+1}}^2
\mathbf1_{\{N_{Q_{m+1}}^2>R\}}\,d\eta_n=0.
\]
Since $S_{\varphi_m}^2\leq N_{Q_{m+1}}^2$ and $S_{\varphi_m}$ is
continuous, weak convergence and this tail bound give
\[
\int S_{\varphi_m}^2\,d\eta_n
\longrightarrow
\int S_{\varphi_m}^2\,d\eta.
\]
This is convergence in the local $L^2$ topology.
\end{proof}

For \(\varphi\in L^1(\mathbb R^d)\), set
\[
\widehat\varphi(\xi)
:=
\int_{\mathbb R^d}
\varphi(x)e^{-2\pi i\langle x,\xi\rangle}\,dx.
\]
If \(\eta\in\PP_{\rho,2}\) and \(\varphi\in\mathcal S(\mathbb R^d)\), then
\(S_\varphi\in L^2(\eta)\).  Indeed, the unit-cube decomposition gives, for
every \(N>d\),
\[
\|S_\varphi\|_{L^2(\eta)}
\leq
\|N_{Q_1}\|_{L^2(\eta)}
\sum_{k\in\mathbb Z^d}\sup_{Q_1+k}|\varphi|
\leq
C_N\sup_x(1+\|x\|)^N|\varphi(x)|.
\]
Define
\[
Z_\varphi^\eta
:=
S_\varphi-\rho\int_{\mathbb R^d}\varphi(x)\,dx.
\]

\begin{proposition}[Bartlett spectrum and continuity]
\label{prop:bartlett-spectrum}\label{prop:bartlett-continuity}
For every \(\eta\in\PP_{\rho,2}\), there is a unique positive symmetric
translation-bounded Radon measure \(\sigma_\eta\) on \(\mathbb R^d\) such
that
\[
\int Z_\varphi^\eta\,\overline{Z_\psi^\eta}\,d\eta
=
\int_{\mathbb R^d}
\widehat\varphi(\xi)\overline{\widehat\psi(\xi)}\,d\sigma_\eta(\xi),
\qquad
\varphi,\psi\in\mathcal S(\mathbb R^d).
\]
The measure \(\sigma_\eta\) is the \emph{Bartlett spectrum} of \(\eta\).
If \(\eta_n\to\eta\) in \(\PP_{\rho,2}\), then
\[
\sigma_{\eta_n}\longrightarrow\sigma_\eta\qquad(n\to\infty)
\]
vaguely on \(\mathbb R^d\).
\end{proposition}

\begin{proof}
The spectral representation follows from
\cite[Proposition~8.2.I and Definition~8.2.II]{DVJ03}, with angular
frequencies $\omega=2\pi\xi$ in the notation of that reference.
Uniqueness follows because every
$h\in C_c^\infty(\mathbb R^d)$ has the form
$\widehat\varphi\,\overline{\widehat\psi}$ for some
$\varphi,\psi\in\mathcal S(\mathbb R^d)$: take
$\widehat\varphi=h$ and $\widehat\psi=1$ on
$\operatorname{supp}h$.
Since the covariance of real statistics is symmetric, so is $\sigma_\eta$.
For translation boundedness, let \(\chi\in C_c^\infty(\mathbb R^d)\) satisfy
\(|\widehat\chi|\geq c>0\) on a neighbourhood \(Q\) of the origin.  Then
\[
c^2\sigma_\eta(a+Q)
\leq
\|Z_{e^{2\pi i\langle a,\cdot\rangle}\chi}^\eta\|_2^2
\leq
\int S_{|\chi|}^2\,d\eta,
\qquad a\in\mathbb R^d.
\]
Fix a compact set \(K\subset\mathbb R^d\).  Choose a real
\(\chi\in C_c^\infty(\mathbb R^d)\) such that \(\widehat\chi\) has no zero
on a neighbourhood of \(K\), for instance a sufficiently
concentrated smooth approximate identity.  Set
\[
\nu_n:=|\widehat\chi|^2\sigma_{\eta_n},
\qquad
\nu:=|\widehat\chi|^2\sigma_\eta.
\]
For \(t\in\mathbb R^d\), the Bartlett identity gives
\[
\int e^{-2\pi i\langle t,\xi\rangle}\,d\nu_n(\xi)
=
\int
Z_\chi^{\eta_n}\,
\overline{Z_{\chi(\,\cdot+t\,)}^{\eta_n}}\,d\eta_n.
\]
Choose $M$ so that
\(\operatorname{supp}\chi\cup(\operatorname{supp}\chi-t)\subset Q_M\).
The integrand is continuous in the vague topology, and its absolute value
is at most \(C(1+N_{Q_M}^2)\).
By \cref{prop:local-L2-characterization}, the laws and second moments of
$N_{Q_M}$ converge. The resulting uniform tail estimate for
$N_{Q_M}^2$ permits passage to the limit in the covariance integral.  At \(t=0\) this also gives
\(\nu_n(\mathbb R^d)\to\nu(\mathbb R^d)\).  L\'evy's continuity theorem
\cite[Theorem~5.3]{Kal02}, after normalization when the limiting mass is
nonzero, yields weak convergence of the finite measures $\nu_n$ to $\nu$ on
$\mathbb R^d$.

If \(h\in C_c(\mathbb R^d)\) is supported in \(K\), then
\[
\int h\,d\sigma_{\eta_n}
=
\int\frac{h}{|\widehat\chi|^2}\,d\nu_n
\longrightarrow
\int\frac{h}{|\widehat\chi|^2}\,d\nu
=
\int h\,d\sigma_\eta.
\]
Hence \(\sigma_{\eta_n}\to\sigma_\eta\) vaguely.
\end{proof}

\begin{lemma}[Bartlett identity for regular windows]
\label{lem:bartlett-windows}
Let \(\eta\in\PP_{\rho,2}\), and let
\(A\subset\mathbb R^d\) be bounded and Borel, with
\(\Vol_d(\partial A)=0\).  Then
\[
\Var_\eta(N_A)
=
\int_{\mathbb R^d}
|\widehat{\mathbf1}_A(\xi)|^2\,d\sigma_\eta(\xi).
\]
\end{lemma}

\begin{proof}
Let \(\varphi_n\in C_c^\infty(\mathbb R^d)\) be uniformly bounded and
supported in a fixed compact set, with
\(\varphi_n\to\mathbf1_A\) pointwise off \(\partial A\) and in \(L^1\).
By \cref{lem:null-window}, \(S_{\varphi_n}\to N_A\) almost surely under $\eta$;
domination by a fixed local count gives convergence in \(L^2(\eta)\).
By the Bartlett identity,
\((\widehat\varphi_n)\) is Cauchy in \(L^2(\sigma_\eta)\).  Since
\(\widehat\varphi_n(\xi)\to\widehat{\mathbf1}_A(\xi)\) for every \(\xi\), its
\(L^2(\sigma_\eta)\)-limit is
\(\widehat{\mathbf1}_A\).  Passing to the limit in the variance identity
gives the result.
\end{proof}

For every $\eta\in\PP_{\rho,2}$, the measure $\sigma_\eta$
has positive mass on some compact annulus. Otherwise
$Z_{\partial_j\varphi}^\eta=0$ almost surely for every
$\varphi\in C_c^\infty(\mathbb R^d)$ and $1\leq j\leq d$.
Taking a countable dense family gives
$\partial_j\omega=0$ distributionally for almost every $\omega$.
Thus $\omega$ is a constant multiple of Lebesgue measure,
contrary to simplicity and positive intensity.
For $0<a<b$, the Bessel asymptotic gives, uniformly for
$a\leq\|\xi\|\leq b$,
\[
\frac1T\int_T^{2T}R^{1-d}
|\widehat{\mathbf1}_{B_R}(\xi)|^2\,dR
=\frac{1}{2\pi^2\|\xi\|^{d+1}}+O_{a,b,d}(T^{-1}).
\]
Choose such an annulus of positive $\sigma_\eta$-mass.
The Bartlett identity and Tonelli's theorem give
\eqref{eq:variance-boundary-limsup}; cf.
\cite[Theorem~5.1]{BB26}.

\subsection{The local-to-global principle}\label{subsec:local-global}

Let $\smash{\EPP^{\mathrm{fr}}_{\rho,2}}$ consist of laws
$\eta\in\smash{\EPP_{\rho,2}}$ for which the translation action
$\mathbb R^d\curvearrowright
(\operatorname{Conf}(\mathbb R^d),\eta)$ is essentially free.

\begin{proposition}\label{prop:free-dense}
The set $\smash{\EPP^{\mathrm{fr}}_{\rho,2}}$ is a dense $G_\delta$ subset of
$\smash{\EPP_{\rho,2}}$.
\end{proposition}

\begin{proof}
For $m\geq1$, let
\[
K_m:=\{t\in\mathbb R^d:m^{-1}\leq \|t\|\leq m\},
\qquad
h_m(\omega):=
\min_{t\in K_m}d_{\mathrm{Conf}}(\theta_t\omega,\omega).
\]
Joint continuity of the translation action and compactness of $K_m$
make $h_m$ continuous, and
$\omega$ has trivial stabilizer exactly when $h_m(\omega)>0$ for every $m$.
Thus, with $C_m:=\{h_m=0\}$,
\[
\smash{\EPP^{\mathrm{fr}}_{\rho,2}}
=
\bigcap_{m,k\geq1}
\{\eta\in\smash{\EPP_{\rho,2}}:\eta(C_m)<k^{-1}\}.
\]
The sets $C_m$ are closed, so $\eta\mapsto\eta(C_m)$ is upper
semicontinuous for weak convergence.  Each strict sublevel set in the
preceding intersection is therefore open.

Fix $\eta\in\smash{\EPP_{\rho,2}}$ and a probability measure
$\upsilon\ll\Vol_d$ supported on $B_1$.  The relation
$\{(\omega,p):p\in\operatorname{supp}\omega\}$ is Borel:
point incidence is equivalent to
$N_{B_{1/n}}(\theta_p\omega)>0$ for all $n\geq1$.
These conditions are Borel by joint continuity of translation
and Borel measurability of $N_{B_{1/n}}$.  Its sections
are countable by local finiteness.  Lusin--Novikov
\cite[Theorem~18.10]{Kec95} gives partial Borel maps with
disjoint graphs covering this relation.  For each $\omega$,
list their defined values in increasing index order as
$p_1(\omega),p_2(\omega),\ldots$.  The domains and values
of these reindexed maps are Borel; the indices range from $1$
to $\#\operatorname{supp}\omega$ when this number is finite.
On the product space
$\operatorname{Conf}(\mathbb R^d)\times B_1^{\mathbb N}$ with measure
$\eta\otimes\upsilon^{\otimes\mathbb N}$, define
\[
\Phi_\varepsilon(\omega,\mathbf u)
:=\sum_{j\leq\#\operatorname{supp}\omega}
\delta_{p_j(\omega)+\varepsilon u_j},
\qquad
\eta_\varepsilon:=(\Phi_\varepsilon)_*
(\eta\otimes\upsilon^{\otimes\mathbb N}).
\]
This independent marking construction is Borel; see
\cite[Definition~5.3]{LP17}. The resulting measure
$\eta_\varepsilon$ does not depend on the chosen enumeration. Translating the original configuration permutes its points
and therefore their independent identically distributed displacements.
Hence $\eta_\varepsilon$ is translation-invariant.
For $\varphi\in C_c(\mathbb R^d)$, the intensity identity gives
\[
\int S_\varphi\,d\eta_\varepsilon
=\int\sum_{p\in\operatorname{supp}\omega}
\int_{B_1}\varphi(p+\varepsilon u)\,d\upsilon(u)\,d\eta(\omega)
=\rho\int\varphi\,d\Vol_d.
\]

Independent marking preserves ergodicity. Let $F,G$ be bounded
observables of the displaced configuration on bounded regions. For
$H=F,G$, set
\[
\overline H(\omega):=
\int H(\Phi_\varepsilon(\omega,\mathbf u))
\,d\upsilon^{\otimes\mathbb N}(\mathbf u).
\]
For sufficiently large $|t|$, the marks affecting $F$ and
$G\circ\theta_t$ belong to disjoint sets of original points.
Conditional independence gives
\[
\int F\,(G\circ\theta_t)\,d\eta_\varepsilon
=\int\overline F\,(\overline G\circ\theta_t)\,d\eta.
\]
Averaging over $t\in Q_T$ and applying the mean ergodic theorem to
$\eta$ shows that the averages of these correlations converge, as
$T\to\infty$, to $(\int F\,d\eta_\varepsilon)
(\int G\,d\eta_\varepsilon)$. Since bounded local functions
generate the configuration $\sigma$-algebra, $\eta_\varepsilon$ is ergodic.

For distinct $p,q$ in any fixed configuration, absolute continuity of
$\upsilon$ gives
\[
\int_{B_1^2}
\mathbf1_{\{p+\varepsilon u=q+\varepsilon v\}}
\,d\upsilon(u)\,d\upsilon(v)=0.
\]
For three distinct points $p,q,r$, likewise,
\[
\int_{B_1^3}
\mathbf1_{\{p+\varepsilon u+r+\varepsilon w=2q+2\varepsilon v\}}
\,d\upsilon(u)\,d\upsilon(v)\,d\upsilon(w)=0.
\]
There are only countably many such pairs and triples. Thus
$\Phi_\varepsilon(\omega,\mathbf u)$ is almost surely simple
and has no three distinct points in arithmetic progression.
The empty configuration is invariant under every translation. Since
$\eta_\varepsilon$ is ergodic and has positive intensity, it gives
zero mass to the empty configuration. If a nonempty set $P$ has period
$t\ne0$, then $x-t,x,x+t\in P$ for any $x\in P$, contradicting
the absence of three-term progressions. Hence $\eta_\varepsilon$
is essentially free.

For $0<\varepsilon<1/2$ and every $m\geq1$,
\[
S_{\varphi_m}(\Phi_\varepsilon(\omega,\mathbf u))^2
\leq N_{Q_{m+2}}(\omega)^2.
\]
The right-hand side is $\eta$-integrable, so
$\eta_\varepsilon\in\smash{\EPP^{\mathrm{fr}}_{\rho,2}}$.
For each $\omega$ and $\mathbf u$, the configurations
$\Phi_\varepsilon(\omega,\mathbf u)$ converge vaguely to
$\omega$ as $\varepsilon\to0$. Dominated convergence gives
$\eta_\varepsilon\Rightarrow\eta$ and convergence of the
$S_{\varphi_m}$-second moments for every $m$.
Thus $\eta_\varepsilon\to\eta$ in the local $L^2$ topology.
\end{proof}

\begin{proposition}[Local-to-global principle]\label{prop:local-global}
Let $\mathcal A\subset\smash{\EPP_{\rho,2}}$ be $G_\delta$.  Suppose that for every
$a\in\Ad_d(X,\mu)$ the set
\[
\{\kappa\in\RR^{a,\mathrm{iso}}_{\rho,2}(X,\mu):
  \kappa_*\mu\in\mathcal A\}
\]
is dense in $\smash{\RR^{a,\mathrm{iso}}_{\rho,2}(X,\mu)}$.  Then $\mathcal A$ is
residual in $\smash{\EPP_{\rho,2}}$.
\end{proposition}

\begin{proof}
Let
$\mathcal U\subset\smash{\EPP_{\rho,2}}$ be nonempty and open.  By
\cref{prop:free-dense}, $\mathcal U$ contains some
$\eta\in\smash{\EPP^{\mathrm{fr}}_{\rho,2}}$.  The probability space
$(\operatorname{Conf}(\mathbb R^d),\eta)$ is nonatomic.  Indeed, since
the action is essentially free, the translates of any atom would give
uncountably many distinct atoms of the same positive mass.
By \cite[Theorem~17.41]{Kec95}, choose a measure-preserving Borel
isomorphism modulo null sets
\[
\psi:(X,\mu)\longrightarrow(\operatorname{Conf}(\mathbb R^d),\eta).
\]
By shrinking the conull Borel sets, we may assume that their
complements are uncountable Borel null sets. Extend $\psi$ to a Borel
bijection using a Borel isomorphism between these complements
\cite[Theorem~15.6]{Kec95}.
Define
\[
a(t,x):=\psi^{-1}(\theta_t\psi(x)),
\qquad(t,x)\in\mathbb R^d\times X.
\]
Since $\psi$ preserves measure and intertwines $a$ with $\theta$,
this action is p.m.p., ergodic, and essentially free.
The map $\kappa_0:=\psi$ is equivariant and generating, with
$(\kappa_0)_*\mu=\eta$.

Continuity of the push-forward map makes
\[
\{\kappa\in\RR^{a,\mathrm{iso}}_{\rho,2}(X,\mu):
  \kappa_*\mu\in\mathcal U\}
\]
an open neighbourhood of $\kappa_0$.  By the density hypothesis it contains
some $\kappa$ with $\kappa_*\mu\in\mathcal A$.  Hence
$\mathcal U\cap\mathcal A\neq\varnothing$, so $\mathcal A$ is dense.
Since $\mathcal A$ is $G_\delta$ in the Polish space
$\smash{\EPP_{\rho,2}}$, it is residual.
\end{proof}

\section{Intermittent hyperuniformity}\label{sec:intermittent-hyperuniformity}

In this section we prove that, for every essentially free ergodic action,
realizations with small or large number variance at prescribed scales are
dense in the local $L^2$ topology. For small variance, we replace points
inside large cells so that their polynomial moments agree with Lebesgue
measure. For large variance, we move points into smaller subcubes while
preserving the number of points in each cell. We then obtain generating
realizations by encoding the underlying action with sparse point triples.

\subsection{Baire reduction}\label{subsec:ihyp-baire}

In this subsection we express intermittent hyperuniformity as a countable
intersection of open variance conditions, reducing its residuality to two
density assertions. Recall \(V_\eta(R)\), \(\Hyp_\rho\), and
\(\IHyp_\rho\) from \cref{subsec:intro-hyperuniformity}.  For
\(a\in\Ad_d(X,\mu)\) and
\(\kappa\in\smash{\RR^a_{\rho,2}(X,\mu)}\), set
\[
\mathcal V_\kappa(R):=V_{\kappa_*\mu}(R).
\]
For $\eta\in\EPP_{\rho,2}$ and $\beta>(d-1)/d$,
\eqref{eq:ihyp-lower} holds if and only if, for every
$\varepsilon>0$, there are arbitrarily large $R$ with
\[
V_\eta(R)<\varepsilon\Vol_d(B_R)^\beta.
\]
For $\alpha<2$, \eqref{eq:ihyp-upper} holds if and only if,
for every $M>0$, there are arbitrarily large $R$ with
\[
V_\eta(R)>M\Vol_d(B_R)^\alpha.
\]
Set
\[
\mathbb Q_-:=\mathbb Q\cap(-\infty,2),
\qquad
\mathbb Q_+:=\mathbb Q\cap\left(\frac{d-1}{d},\infty\right).
\]

\begin{proposition}\label{prop:ihyp-gdelta}
The set $\IHyp_\rho$ is $G_\delta$ in $\smash{\EPP_{\rho,2}}$.  For every
$a\in\Ad_d(X,\mu)$, the set
\[
\{\kappa\in\RR^a_{\rho,2}(X,\mu):\kappa_*\mu\in\IHyp_\rho\}
\]
is $G_\delta$ in $\smash{\RR^a_{\rho,2}(X,\mu)}$.
\end{proposition}

\begin{proof}
Fix $\eta\in\PP_{\rho,2}$.  The map $R\mapsto V_\eta(R)$ is continuous.
Indeed, if $R_n\to R$, then
$N_{B_{R_n}}\to N_{B_R}$ $\eta$-almost surely by
\cref{lem:null-window}.  For all sufficiently large $n$,
$N_{B_{R_n}}\leq N_{B_{R+1}}$. Since
$\eta\in\PP_{\rho,2}$, we have $N_{B_{R+1}}\in L^2(\eta)$,
so the convergence also holds in $L^2(\eta)$.

Hence rational radii suffice:
\[
\begin{split}
\IHyp_\rho
={}&
\bigcap_{\substack{\beta\in\mathbb Q_+\\ j,N\geq1}}
\ \bigcup_{\substack{R\in\mathbb Q\\ R>N}}
\left\{\eta\in\EPP_{\rho,2}:
V_\eta(R)<j^{-1}\Vol_d(B_R)^\beta\right\}\\
&\cap
\bigcap_{\substack{\alpha\in\mathbb Q_-\\ j,N\geq1}}
\ \bigcup_{\substack{R\in\mathbb Q\\ R>N}}
\left\{\eta\in\EPP_{\rho,2}:
V_\eta(R)>j\,\Vol_d(B_R)^\alpha\right\}.
\end{split}
\]
For fixed \(R>0\), the map \(\eta\mapsto V_\eta(R)\) is continuous by
\cref{prop:global-polish}.  This proves the first assertion, and continuity
of the push-forward map gives the second.
\end{proof}

For $\beta>(d-1)/d$, $\alpha<2$, and $j,N\geq1$, define
\[
\mathcal L^a_{\beta,j,N}
:=
\bigcup_{n\geq N}
\left\{
\kappa\in\RR^a_{\rho,2}(X,\mu):
\mathcal V_\kappa(n)<j^{-1}\Vol_d(B_n)^\beta
\right\},
\]
and
\[
\mathcal H^a_{\alpha,j,N}
:=
\bigcup_{n\geq N}
\left\{
\kappa\in\RR^a_{\rho,2}(X,\mu):
\mathcal V_\kappa(n)>j\,\Vol_d(B_n)^\alpha
\right\}.
\]
These sets are open by \cref{prop:rooted-polish}.

\begin{proposition}[Baire reduction]\label{prop:ihyp-baire-reduction}
Suppose that, for every $a\in\Ad_d(X,\mu)$,
$\beta>(d-1)/d$, $\alpha<2$, and $j,N\geq1$, the sets
$\smash{\mathcal L^a_{\beta,j,N}}$ and
$\smash{\mathcal H^a_{\alpha,j,N}}$ are dense in
$\smash{\RR^a_{\rho,2}(X,\mu)}$, and that
$\smash{\RR^{a,\mathrm{iso}}_{\rho,2}(X,\mu)}$ is dense there.  Then
\[
\{\kappa\in\RR^{a,\mathrm{iso}}_{\rho,2}(X,\mu):
  \kappa_*\mu\in\IHyp_\rho\}
\]
is residual in $\smash{\RR^{a,\mathrm{iso}}_{\rho,2}(X,\mu)}$ for every
$a\in\Ad_d(X,\mu)$, and $\IHyp_\rho$ is residual in
$\smash{\EPP_{\rho,2}}$.
\end{proposition}

\begin{proof}
Let $G_a$ be the intersection of
$\smash{\RR^{a,\mathrm{iso}}_{\rho,2}(X,\mu)}$ with the sets
$\smash{\mathcal L^a_{\beta,j,N}}$ and
$\smash{\mathcal H^a_{\alpha,j,N}}$ over
$j,N\geq1$, $\beta\in\mathbb Q_+$, and $\alpha\in\mathbb Q_-$.  The
hypotheses and \cref{prop:rooted-polish} show that $G_a$ is a dense
$G_\delta$ subset of $\smash{\RR^a_{\rho,2}(X,\mu)}$.

For $\kappa\in G_a$ and $\beta>(d-1)/d$, choose
$\beta'\in\mathbb Q_+$ with $\beta'<\beta$; for $\alpha<2$, choose
$\alpha'\in\mathbb Q_-$ with $\alpha<\alpha'$.  Since
$\Vol_d(B_n)\to\infty$ as $n\to\infty$, the corresponding integer-radius inequalities imply
\(\kappa_*\mu\in\IHyp_\rho\).  The pullback of \(\IHyp_\rho\) is therefore
dense in every generating realization space and is $G_\delta$ there by
\cref{prop:ihyp-gdelta}.  The local-to-global principle
\cref{prop:local-global} yields the global conclusion.
\end{proof}

\subsection{Low-variance approximation}\label{subsec:low-variance-surgery}

In this subsection we prove that low-variance realizations are dense.
The construction modifies a point configuration separately in the
interior of each of many small cubes. In each cube, we retain the
points outside a prescribed interior region and replace the points
inside that region by a finite set. The union of the retained and replacement points has
the prescribed count and polynomial moments on each small cube. These identities
will allow us to apply \cref{lem:uniform-periodic-variance}.

Write \(I:=[0,1)^d\). For $T>0$, a \emph{$T$-cell} is a half-open
cube $v+TI$, and a \emph{$T$-grid} is the partition
$\{v+Tk+TI:k\in\mathbb Z^d\}$ for a fixed $v\in\mathbb R^d$.
In particular, an $L$-cell has side length $L$; each $L$-grid
below has a specified translation.

\begin{proposition}[Low-variance density]\label{prop:low-variance-density}
For every \(a\in\Ad_d(X,\mu)\), \(\beta>(d-1)/d\), and \(j,N\geq1\),
the set \(\smash{\mathcal L^a_{\beta,j,N}}\) is dense in
\(\smash{\RR^a_{\rho,2}(X,\mu)}\).
\end{proposition}

\subsubsection{Ergodic approximation of cell moments}

We first estimate the number of points in large sets. Ergodicity gives
$L^2$ convergence of scaled point counts, also when the points are
weighted by a continuous function. We also need to control the expected
number of points in a surrounding cube on the event that a prescribed
weighted count deviates from its limit. This second estimate will be
used to bound the points in cells where the replacement construction
cannot be applied.

\begin{lemma}\label{lem:weighted-cell-ergodic}
Let \(a\in\Ad_d(X,\mu)\) and
\(\kappa\in\RR^a_{\rho,2}(X,\mu)\). Let
\(A\subset\mathbb R^d\) be a bounded Borel set satisfying
\(\Vol_d(\partial_r A)=O(r)\) as \(r\downarrow0\), and let
\(g\in C(\overline A)\). For \(s>0\), put
\begin{equation}\label{eq:weighted-cell-statistic}
\Sigma_s(x):=s^{-d}\sum_{u\in\operatorname{supp}\kappa(x)\cap sA}g(u/s).
\end{equation}
Then
\begin{equation}\label{eq:weighted-cell-L2}
\Sigma_s\longrightarrow\rho\int_A g(u)\,du
\qquad\text{in }L^2(X,\mu)\quad(s\to\infty).
\end{equation}
Moreover, for every bounded cube \(B\supset A\) and every
\(\varepsilon>0\),
\begin{equation}\label{eq:weighted-cell-exceptional}
 s^{-d}\int_{\{x:\,|\Sigma_s(x)-\rho\int_Ag|>\varepsilon\}}
 N_{sB}(\kappa(x))\,d\mu(x)\longrightarrow0
 \qquad(s\to\infty).
\end{equation}
\end{lemma}

\begin{proof}
Fix $a$ and $\kappa$ as in the statement. We first consider a
half-open cube $Q\subset\mathbb R^d$ and the unweighted count
$N_{sQ}\circ\kappa$. Choose a nonnegative
$\varphi\in C_c(\mathbb R^d)$ with integral one, and let $r>0$
satisfy $\operatorname{supp}\varphi\subset B_r$. Set
\(H(x):=\sum_{u\in\operatorname{supp}\kappa(x)}\varphi(u)\).
Equivariance gives
\begin{equation}\label{eq:weighted-cell-equivariance}
H(t.x)=\sum_{u\in\operatorname{supp}\kappa(x)}\varphi(u-t).
\end{equation}
The mean ergodic theorem and $\int H\,d\mu=\rho$ give
\begin{equation}\label{eq:weighted-cell-ergodic-mean}
\frac1{\Vol_d(sQ)}\int_{sQ}H(t.x)\,dt
\longrightarrow\rho
\qquad\text{in }L^2(X,\mu).
\end{equation}
On the other hand, Fubini and \eqref{eq:weighted-cell-equivariance} give
\begin{equation}\label{eq:weighted-cell-convolution}
\int_{sQ}H(t.x)\,dt
=\sum_{u\in\operatorname{supp}\kappa(x)}
  (\mathbf1_{sQ}*\varphi)(u).
\end{equation}
The functions $\mathbf1_{sQ}*\varphi$ and $\mathbf1_{sQ}$
coincide outside $\partial_r(sQ)$, and their difference has absolute
value at most one. For $s\geq1$, this boundary neighbourhood intersects at
most $C_Qs^{d-1}$ cubes of the unit grid. Since $\kappa_*\mu$ is translation-invariant and has finite local
second moments, the triangle inequality gives
\begin{equation}\label{eq:weighted-cell-boundary-error}
\left\|\int_{sQ}H(t.\cdot)\,dt-N_{sQ}\circ\kappa\right\|_{L^2(\mu)}
\leq C_Qs^{d-1}
\|N_{[0,1)^d}\circ\kappa\|_{L^2(\mu)}.
\end{equation}
Dividing \eqref{eq:weighted-cell-boundary-error} by $s^d$ and
using \eqref{eq:weighted-cell-ergodic-mean} yields
\begin{equation}\label{eq:weighted-cell-cube-limit}
s^{-d}N_{sQ}\circ\kappa\longrightarrow
\rho\Vol_d(Q)\qquad\text{in }L^2(\mu).
\end{equation}

To pass from cube counts to the weighted statistic $\Sigma_s$,
approximate $g\mathbf1_A$ by step functions on a grid. We regard
$g\mathbf1_A$ as a function on $\mathbb R^d$, equal to zero
outside $A$. For $h>0$, let $\mathscr J_h$ consist of the half-open cubes
$J=h(k+I)$, $k\in\mathbb Z^d$, that are contained in $A$.
This family is finite. Choose $u_J\in J$ for each $J\in\mathscr J_h$,
and put
\[
\psi_h:=\sum_{J\in\mathscr J_h}g(u_J)\mathbf1_J,
\qquad
\Sigma_{s,h}(x):=s^{-d}\sum_{u\in\operatorname{supp}\kappa(x)}
\psi_h(u/s).
\]
Our aim is to compare $\Sigma_s$ in
\eqref{eq:weighted-cell-statistic} with $\Sigma_{s,h}$, apply
\eqref{eq:weighted-cell-cube-limit} with $h$ fixed, and then let
$h\downarrow0$.

Define
\[
\omega_g(\delta):=\sup\{|g(u)-g(v)|:
 u,v\in\overline A,\ \|u-v\|\leq\delta\},
\qquad U_h:=A\cap\partial_{\sqrt d h}A.
\]
Each cube of side length $h$ has diameter $\sqrt d\,h$.
Consequently, if $u\in A\setminus U_h$, its grid cube is contained
in $A$: otherwise the segment from $u$ to a point of the cube outside
$A$ would intersect $\partial A$ at distance at most $\sqrt d\,h$ from
$u$. If $J$ is this cube, then $\|u-u_J\|\leq\sqrt d\,h$, so
$|g(u)-\psi_h(u)|\leq\omega_g(\sqrt d\,h)$.
For $u\in U_h$ we have
$|g(u)-\psi_h(u)|\leq2\|g\|_\infty$, and both functions vanish
outside $A$.

Set $M_2:=\|N_I\circ\kappa\|_{L^2(\mu)}$.
If $E\subset\mathbb R^d$ is bounded and Borel, the unit cubes intersecting
$E$ are disjoint and lie in $E+[-1,1]^d$. Translation invariance of
$\kappa_*\mu$ and the triangle inequality thus give
\begin{equation}\label{eq:weighted-cell-local-count-bound}
\|N_E\circ\kappa\|_{L^2(\mu)}
\leq M_2\Vol_d(E+[-1,1]^d).
\end{equation}
For fixed $h>0$ and $s\to\infty$, this bound applied to $sA$ and
$sU_h$ yields
\[
\begin{aligned}
\limsup_{s\to\infty}s^{-d}\|N_{sA}\circ\kappa\|_2
 &\leq M_2\Vol_d(A),\\
\limsup_{s\to\infty}s^{-d}\|N_{sU_h}\circ\kappa\|_2
 &\leq M_2\Vol_d(\partial_{2\sqrt d h}A).
\end{aligned}
\]
Indeed, for all sufficiently large $s$ (depending on $h$),
$U_h+[-1/s,1/s]^d\subset\partial_{2\sqrt d h}A$.
The pointwise bounds above therefore imply
\begin{equation}\label{eq:weighted-cell-step-error}
\begin{aligned}
\limsup_{s\to\infty}\|\Sigma_s-\Sigma_{s,h}\|_{L^2(\mu)}
\leq M_2\Bigl(&\omega_g(\sqrt d h)\Vol_d(A)\\
 &+2\|g\|_\infty\Vol_d(\partial_{2\sqrt d h}A)\Bigr).
\end{aligned}
\end{equation}
For every fixed $h>0$, $\psi_h$ is a finite linear combination of
indicators of half-open cubes. Hence
\eqref{eq:weighted-cell-cube-limit} gives
\begin{equation}\label{eq:weighted-cell-step-limit}
\Sigma_{s,h}\longrightarrow\rho\int_{\mathbb R^d}\psi_h(u)\,du
\qquad\text{in }L^2(\mu)\quad(s\to\infty).
\end{equation}
Combining the preceding two equations with the triangle inequality
produces, for every fixed $h>0$,
\begin{equation}\label{eq:weighted-cell-limit-comparison}
\begin{aligned}
\limsup_{s\to\infty}
\left\|\Sigma_s-\rho\int_A g\right\|_{L^2(\mu)}
\leq{}&M_2\Bigl(\omega_g(\sqrt d h)\Vol_d(A)
 +2\|g\|_\infty\Vol_d(\partial_{2\sqrt d h}A)\Bigr)\\
 &+\rho\left|\int_{\mathbb R^d}\psi_h-\int_Ag\right|.
\end{aligned}
\end{equation}
The right-hand side tends to zero as $h\downarrow0$: uniform
continuity gives $\omega_g(\sqrt d h)\to0$, the hypothesis on
$\partial_rA$ controls the boundary term, and the pointwise bounds
above imply $\|\psi_h-g\mathbf1_A\|_{L^1}\to0$.
Thus the left-hand side of
\eqref{eq:weighted-cell-limit-comparison}, which does not depend on
$h$, is zero. This proves \eqref{eq:weighted-cell-L2}.

Fix a cube $B\supset A$ and $\varepsilon>0$, and set
\[
E_s:=\left\{x:\left|\Sigma_s(x)-\rho\int_Ag\right|>
\varepsilon\right\},\qquad
Y_s:=s^{-d}N_{sB}\circ\kappa.
\]
Convergence \eqref{eq:weighted-cell-L2} implies $\mu(E_s)\to0$.
By \eqref{eq:weighted-cell-cube-limit},
$Y_s\to\rho\Vol_d(B)$ in $L^1(\mu)$; if $B$ is not
half-open, its boundary has zero expected count, so the choice of
boundary convention does not affect $Y_s$ almost surely. Since
$Y_s\geq0$,
\[
\int_{E_s}Y_s\,d\mu
\leq\|Y_s-\rho\Vol_d(B)\|_{L^1(\mu)}
 +\rho\Vol_d(B)\mu(E_s)\longrightarrow0,
\]
which proves \eqref{eq:weighted-cell-exceptional}.
\end{proof}

For a finite set of pairs \((A_i,g_i)\), \(1\leq i\leq n\), each
satisfying the hypotheses of \cref{lem:weighted-cell-ergodic}, fix
\(\varepsilon_i>0\) and a bounded cube \(B\) containing all \(A_i\).
Define
\[
E_{s,i}:=\left\{x:\left|
 s^{-d}\sum_{u\in\operatorname{supp}\kappa(x)\cap sA_i}g_i(u/s)
 -\rho\int_{A_i}g_i\right|>\varepsilon_i\right\}.
\]
Applying \eqref{eq:weighted-cell-L2} and
\eqref{eq:weighted-cell-exceptional} to each of the $n$ pairs gives
$\mu(\bigcup_{i=1}^nE_{s,i})\to0$ and
\begin{equation}\label{eq:weighted-cell-finite-family}
s^{-d}\int_{\bigcup_{i=1}^nE_{s,i}}N_{sB}\circ\kappa\,d\mu
\leq\sum_{i=1}^n s^{-d}\int_{E_{s,i}}N_{sB}\circ\kappa\,d\mu
\longrightarrow0.
\end{equation}

\subsubsection{Completion of polynomial moments}\label{subsec:relative-moment-repair}

Fix a closed cube $D\subset(0,1)^d$. Given a finite set of points
$F\subset I\setminus\operatorname{int}D$, we seek a finite set
$A\subset\operatorname{int}D$ such that the union $F\cup A$ has a
prescribed number of points and its polynomial sums equal the
corresponding integrals over $I$. Thus the points of $F$ remain fixed;
only the points of $A$ are to be determined. The proposition gives a
Borel choice of $A$ when the count and polynomial sums of $F$ are
sufficiently close to their prescribed values. It also gives a lower
bound on the distance between points of $A$ and on their distance to
$\partial D$.

For a positive integer $p$, let $\mathcal P_{p-1}$ be the space of
real polynomials on $\mathbb R^d$ of degree less than $p$, and set
\[
\mathcal P_{p-1}^0:=\{P\in\mathcal P_{p-1}:P(0)=0\}.
\]
For a finite set $A\subset\mathbb R^d$, write
\[
\operatorname{sep}(A):=
\min\{\|u-v\|:u,v\in A,\ u\neq v\},
\]
with $\operatorname{sep}(A):=+\infty$ when $\#A\leq1$.

\begin{proposition}[Completion of polynomial moments]\label{prop:relative-moment-repair}
Fix a closed cube \(D\subset(0,1)^d\), a positive integer
\(p\), and a basis \(P_1,\ldots,P_m\) of
\(\mathcal P_{p-1}^0\). Put
\(\vartheta:=\Vol_d(D)\) and
\(\mathbf P:=(P_1,\ldots,P_m)\); for \(p=1\), the basis is empty.
There exist \(\varepsilon_0>0\) and \(q_0\in\mathbb N\), depending
on \(D,p,\mathbf P\), and \(\sigma_d>0\), depending only on
\(d\), with the following property.

Fix an integer \(q\geq q_0\) and
\(0<\varepsilon\leq\varepsilon_0\). Let
\(\mathscr F_{q,\varepsilon}\) consist of the finite sets
\(F\subset I\setminus\operatorname{int}D\) satisfying
\begin{equation}\label{eq:moment-exterior-hypotheses}
\bigl|\#F-(1-\vartheta)q\bigr|\leq\varepsilon q
\quad\text{and}\quad
\left|\sum_{u\in F}P_i(u)-q\int_{I\setminus D}P_i(u)\,du\right|
\leq\varepsilon q\quad(1\leq i\leq m).
\end{equation}
Then $\mathscr F_{q,\varepsilon}$ is Borel in
$\operatorname{Conf}(\mathbb R^d)$, and there is a Borel map
\[
F\longmapsto\delta_{A(q,F)},
\qquad \mathscr F_{q,\varepsilon}\longrightarrow
\operatorname{Conf}(\mathbb R^d),
\]
with $A(q,F)\subset\operatorname{int}D$, such that
\begin{equation}\label{eq:moment-completed-identity}
\sum_{u\in F\cup A(q,F)}P(u)
=q\int_I P(u)\,du
\qquad(P\in\mathcal P_{p-1}),
\end{equation}
and
\begin{equation}\label{eq:moment-completed-separation}
\operatorname{sep}(A(q,F))\geq\sigma_d q^{-1/d},
\qquad
\operatorname{dist}(A(q,F),\partial D)\geq\sigma_d q^{-1/d}.
\end{equation}
\end{proposition}

\begin{remark*}
For each fixed $0<\varepsilon\leq\varepsilon_0$, the family
$\mathscr F_{q,\varepsilon}$ is nonempty for all sufficiently large
$q$. Indeed, let $h=q^{-1/d}$ and take points of a translate of
$h\mathbb Z^d$ in $I\setminus D$ whose distance from $D$ is at least
$h/4$. Riemann-sum estimates on the cube boundaries give an error
$O(q^{1-1/d})$ in both the point count and each polynomial sum in
\eqref{eq:moment-exterior-hypotheses}; for fixed $\varepsilon$, this
is at most $\varepsilon q$ when $q$ is large. This nonemptiness
assertion does not hold uniformly for all
$0<\varepsilon\leq\varepsilon_0$ at a fixed $q$.
\end{remark*}

\begin{proof}
Fix $q\geq q_0$, $0<\varepsilon\leq\varepsilon_0$, and
$F\in\mathscr F_{q,\varepsilon}$. Put $k:=q-\#F$.
We seek $A\subset\operatorname{int}D$ with exactly $k$ points
and, when $m\geq1$, with
\begin{equation}\label{eq:moment-interior-target}
\sum_{u\in A}\mathbf P(u)
=q\int_I\mathbf P(u)\,du-\sum_{u\in F}\mathbf P(u).
\end{equation}
The count bound in \eqref{eq:moment-exterior-hypotheses} gives
\[
|k-\vartheta q|\leq\varepsilon q.
\]
We take $\varepsilon_0\leq\vartheta/2$, so that
$\vartheta q/2\leq k\leq3\vartheta q/2$. In particular,
$k\to\infty$ uniformly over the sets $F$ under consideration as
$q\to\infty$.

There are constants $c_d,C_d>0$ and $K_d\in\mathbb N$,
depending only on $d$, such that for $k\geq K_d$ the cube $D$
can be partitioned into $k$ axis-parallel rectangles of equal
volume $\vartheta/k$, with every side length between
$2c_d(\vartheta/k)^{1/d}$ and $C_d(\vartheta/k)^{1/d}$.
Increase $q_0$ so that $k\geq K_d$. We construct these
partitions inductively on $d$. For $d=1$,
use $k$ equal intervals. For $d\geq2$, divide the first coordinate
into $n:=\lceil k^{1/d}\rceil$ intervals. Choose integers
$k_1,\ldots,k_n$ in $\{\lfloor k/n\rfloor,\lceil k/n\rceil\}$
with sum $k$, in a fixed order. Write $\ell:=\vartheta^{1/d}$
for the side length of $D$. Partition its first coordinate interval
into consecutive intervals of lengths $\ell k_i/k$, $1\leq i\leq n$.
Over the $i$th interval, partition the $(d-1)$-dimensional cube
into $k_i$ equal-volume rectangles by induction. The stated side-length bounds
follow from $k_i\asymp k^{(d-1)/d}$, with constants depending only
on $d$. Fix one such partition for each $k$, and let $A_k^0$
consist of its $k$ rectangle centres. Then
\[
\operatorname{sep}(A_k^0),\quad
\operatorname{dist}(A_k^0,\partial D)
\geq c_d(\vartheta/k)^{1/d}
\geq c_d(2/3)^{1/d}q^{-1/d}.
\]
For every $f\in C^1(D)$, integrating the difference between
$f$ at the centre and at each point of its rectangle yields
\begin{equation}\label{eq:rectangle-centre-error}
\left|\frac1k\sum_{u\in A_k^0}f(u)
-\frac1\vartheta\int_D f(u)\,du\right|
\leq C_d(\vartheta/k)^{1/d}\|\nabla f\|_{L^\infty(D)}.
\end{equation}

If $p=1$, the only polynomial condition concerns constants.
Set $A(q,F):=A_k^0$. It has $k$ points and satisfies the required
separation bounds. The choice is Borel because the fixed set
$A_k^0$ depends only on the integer $k=q-\#F$.
For the rest of the proof, assume $p\geq2$.

Dividing \eqref{eq:moment-interior-target} by $k$, the average
of $\mathbf P$ over $A$ must equal
\[
y(F):=\frac1k\left(q\int_I\mathbf P(u)\,du
-\sum_{u\in F}\mathbf P(u)\right).
\]
By splitting the integral over $I$ into $D$ and
$I\setminus D$, we have
\[
\begin{aligned}
y(F)-\frac1\vartheta\int_D\mathbf P(u)\,du
={}&\left(\frac qk-\frac1\vartheta\right)
\int_D\mathbf P(u)\,du\\
&+\frac1k\left(q\int_{I\setminus D}\mathbf P(u)\,du
-\sum_{u\in F}\mathbf P(u)\right).
\end{aligned}
\]
The bounds in \eqref{eq:moment-exterior-hypotheses} and
$k\geq\vartheta q/2$ therefore imply
\[
\left\|y(F)-\frac1\vartheta\int_D\mathbf P(u)\,du
\right\|\leq C\varepsilon.
\]
Here and below, constants may depend on $D,p$, and the fixed
polynomial basis, but not on $q$, $\varepsilon$, or $F$.

To correct the polynomial sums of $A_k^0$, we vary its points
through an $m$-parameter family of maps and choose the directions
so that the derivative of the moment map is invertible.
Choose a nonnegative function
$\chi\in C_c^\infty(\operatorname{int}D)$ positive on a
nonempty open set. For $1\leq j\leq m$, put
$X_j(u):=\chi(u)\nabla P_j(u)$, and define
\[
\Phi_t(u):=u+\sum_{j=1}^m t_jX_j(u)
\qquad(t\in\mathbb R^m,\ u\in\mathbb R^d).
\]
We shall apply $\Phi_t$ to every point of $A_k^0$.
The derivative at $t=0$ of the vector
$\vartheta^{-1}\int_D\mathbf P(\Phi_t(u))\,du$
is the matrix
\[
J_{ij}:=\frac1\vartheta\int_D
\chi(u)\nabla P_i(u)\cdot\nabla P_j(u)\,du.
\]
For $\xi\in\mathbb R^m$, let
$P_\xi:=\sum_j\xi_jP_j$. Then
\[
\xi^{\mathsf T}J\xi
=\frac1\vartheta\int_D\chi(u)\|\nabla P_\xi(u)\|^2\,du.
\]
If this integral is zero, the polynomial $\nabla P_\xi$
vanishes on a nonempty open set. Hence $P_\xi$ is constant,
and $P_\xi(0)=0$ makes it zero. Since the $P_j$ form a basis,
$\xi=0$. Thus $J$ is positive definite.

To obtain the prescribed polynomial sums, define
\[
G_{k,F}(t):=
\frac1k\sum_{u\in A_k^0}\mathbf P(\Phi_t(u))-y(F),
\qquad J_k:=DG_{k,F}(0).
\]
The matrix $J_k$ depends on $k$ but not on $F$, since $y(F)$
is independent of $t$. Applying
\eqref{eq:rectangle-centre-error} to the coordinates of
$\mathbf P$ and to the functions
$\chi\nabla P_i\cdot\nabla P_j$ gives
\[
\|G_{k,F}(0)\|\leq C(\varepsilon+k^{-1/d}),
\qquad \|J_k-J\|\leq Ck^{-1/d}.
\]
After increasing $q_0$, all the matrices $J_k$ that occur here
are invertible, with uniformly bounded inverses.
The second derivatives of $G_{k,F}$ are uniformly bounded for
$\|t\|\leq\tau$, where $\tau>0$ is fixed; they too do not
depend on $F$.

Set $\delta_k:=C_1(\varepsilon+k^{-1/d})$ and
\[
\mathcal N_{k,F}(t):=t-J_k^{-1}G_{k,F}(t).
\]
Choose $C_1$ so that
$\|\mathcal N_{k,F}(0)\|\leq\delta_k/2$.
By decreasing $\varepsilon_0$ and increasing $q_0$, we may
assume that $\delta_k\leq\tau$ and, for all
$\|t\|\leq\delta_k$,
\[
\|D\mathcal N_{k,F}(t)\|
=\|J_k^{-1}(J_k-DG_{k,F}(t))\|
\leq C\delta_k\leq\tfrac12.
\]
It follows that $\mathcal N_{k,F}$ maps the closed ball
$\{\|t\|\leq\delta_k\}$ into itself and is a contraction
there. Let $t(F)$ be its unique fixed point. Then
$G_{k,F}(t(F))=0$ and $\|t(F)\|\leq\delta_k$.

The vector fields $X_j$ vanish outside
$\operatorname{int}D$. By further decreasing $\varepsilon_0$
and increasing $q_0$, we can ensure that
$\|D\Phi_t-I\|\leq1/2$ on $\mathbb R^d$ whenever
$\|t\|\leq\delta_k$. Hence
\[
\tfrac12\|u-v\|\leq\|\Phi_t(u)-\Phi_t(v)\|
\leq\tfrac32\|u-v\|
\qquad(u,v\in\mathbb R^d).
\]
Moreover, $\Phi_t$ is a bijection of $\mathbb R^d$:
for each $y$, the map $u\mapsto y-(\Phi_t(u)-u)$
is a contraction. Since $\Phi_t$ is the identity outside
$\operatorname{int}D$, it maps $D$ onto itself and fixes
$\partial D$. Therefore, with
\[
A(q,F):=\Phi_{t(F)}(A_k^0),
\]
we have $A(q,F)\subset\operatorname{int}D$ and
\[
\operatorname{sep}(A(q,F)),\quad
\operatorname{dist}(A(q,F),\partial D)
\geq\tfrac12c_d(2/3)^{1/d}q^{-1/d}.
\]
We may take $\sigma_d:=\tfrac12c_d(2/3)^{1/d}$. The identity $G_{k,F}(t(F))=0$ gives
\[
\sum_{u\in A(q,F)}\mathbf P(u)
=q\int_I\mathbf P(u)\,du-\sum_{u\in F}\mathbf P(u).
\]
The sets $F$ and $A(q,F)$ are disjoint, since the former
lies outside $\operatorname{int}D$ and the latter lies inside.
Thus their union has $\#F+k=q$ points. The polynomial identities above and the identity for constants give
\[
\sum_{u\in F\cup A(q,F)}P(u)=q\int_I P(u)\,du
\qquad(P\in\mathcal P_{p-1}).
\]

Finally, the set $\mathscr F_{q,\varepsilon}$ is Borel: the condition $F\subset I\setminus\operatorname{int}D$, the
finite-point condition, and the count and polynomial-sum inequalities
are Borel conditions on configurations. On this domain, $k=q-\#F$
and $y(F)$ are Borel. The sets $A_k^0$, matrices $J_k$,
and maps $\Phi_t$ fixed for each integer $k$ make
$(F,t)\mapsto\mathcal N_{k,F}(t)$
Borel in $F$ and continuous in $t$. Starting with $t_0(F)=0$,
define $t_{n+1}(F):=\mathcal N_{k,F}(t_n(F))$.
Every $t_n$ is Borel, and contraction gives
$t_n(F)\to t(F)$. Hence $F\mapsto t(F)$ and
$F\mapsto\delta_{A(q,F)}$ are Borel.
\end{proof}

Fix $q\geq q_0$ and $0<\varepsilon\leq\varepsilon_0$, and set
$L:=(q/\rho)^{1/d}$. For an $L$-cell $P=v+LI$, put
$D_P:=v+LD$. Given a locally finite point set $\Lambda$, define
\begin{equation}\label{eq:rescaled-exterior-points}
F_{L,P}(\Lambda):=
\left\{\frac{u-v}{L}:
 u\in\Lambda\cap(P\setminus\operatorname{int}D_P)\right\}.
\end{equation}
The subscript $L$ records the scale used to normalize the exterior
points, and $P$ specifies their location. Whenever
$F_{L,P}(\Lambda)\in\mathscr F_{q,\varepsilon}$, define
\[
\mathcal R_P(\delta_\Lambda):=
\delta_{(\Lambda\setminus\operatorname{int}D_P)
\cup(v+LA(q,F_{L,P}(\Lambda)))}.
\]
The condition $F_{L,P}(\Lambda)\in\mathscr F_{q,\varepsilon}$ is Borel
in $\delta_\Lambda$. By \cref{prop:relative-moment-repair}, the map
$\mathcal R_P$ is Borel on this domain, agrees with $\delta_\Lambda$
outside $\operatorname{int}D_P$, and satisfies
\begin{equation}\label{eq:rescaled-completed-moments}
\int_PH(u)\,d\mathcal R_P(\delta_\Lambda)(u)
=\rho\int_PH(u)\,du
\qquad(H\in\mathcal P_{p-1}).
\end{equation}
In particular, $\mathcal R_P(\delta_\Lambda)(P)=q$.
Indeed, \eqref{eq:rescaled-completed-moments} is
\eqref{eq:moment-completed-identity} after the change of
variables $u\mapsto v+Lu$.

\subsubsection{Variance from exact cell moments}

We estimate the number variance of periodic point sets whose
polynomial moments on each $L$-cell agree with those of $\rho\,du$.
The estimate is uniform in the period and in the point set; its
constant depends on $L$. In the construction below, we choose an
integer $q$ sufficiently large and set $L=(q/\rho)^{1/d}$. Thus
$q=\rho L^d$ is the prescribed number of points in each modified cell.

\begin{lemma}[Periodic variance from exact moments]\label{lem:uniform-periodic-variance}
Fix a positive integer \(p\) with \(2p>d+1\), and an integer
\(q\geq1\). Set \(L:=(q/\rho)^{1/d}\). There is \(K_L<\infty\) such that
the following holds for every integer $m\geq1$. Put $S:=mL$ and
\[
\mathscr C_{L,m}:=\{Lk+[0,L)^d:k\in\{0,\ldots,m-1\}^d\}.
\]
Suppose $\Lambda\subset[0,S)^d$ is finite and, for every
$Q\in\mathscr C_{L,m}$,
\[
\#(\Lambda\cap Q)=q,
\qquad
\sum_{u\in\Lambda\cap Q}P(u)
=\rho\int_QP(u)\,du
\qquad(P\in\mathcal P_{p-1}).
\]
Define the $S\mathbb Z^d$-periodic point set
\[
\Lambda^{\rm per}:=\bigcup_{v\in S\mathbb Z^d}(\Lambda+v).
\]
Then
\[
\frac1{S^d}\int_{[0,S)^d}
\left|
N_{B_R}(\delta_{\Lambda^{\rm per}-t})-\rho\Vol_d(B_R)
\right|^2\,dt
\leq K_LR^{d-1}
\qquad(R\geq1).
\]
\end{lemma}

\begin{proof}
For each $Q\in\mathscr C_{L,m}$, define the signed measure
\[
\Delta_Q:=\delta_{\Lambda\cap Q}-\rho\Vol_d|_Q.
\]
The hypotheses imply
\begin{equation}\label{eq:periodic-cell-vanishing}
\int_Q P(u)\,d\Delta_Q(u)=0
\quad(P\in\mathcal P_{p-1}),
\qquad \|\Delta_Q\|_{\rm TV}\leq2q.
\end{equation}
We use \eqref{eq:periodic-cell-vanishing} to express each
$\Delta_Q$ as distributional derivatives of order $p$.

Let $c_Q$ be the centre of $Q$. For
$\varphi\in C_c^\infty(\mathbb R^d)$, define its Taylor polynomial
at $c_Q$ of degree less than $p$ by
\[
P_{Q,\varphi}(u):=
\sum_{|\alpha|<p}\frac{\partial^\alpha\varphi(c_Q)}{\alpha!}
(u-c_Q)^\alpha.
\]
Here $\alpha\in\mathbb N_0^d$ is a multi-index,
$\alpha!:=\prod_i\alpha_i!$, and
$\partial^\alpha:=\prod_i\partial_i^{\alpha_i}$.
Taylor's integral remainder formula gives
\begin{equation}\label{eq:periodic-taylor-remainder}
\varphi(u)-P_{Q,\varphi}(u)
=p\sum_{|\alpha|=p}\frac{(u-c_Q)^\alpha}{\alpha!}
\int_0^1(1-s)^{p-1}
\partial^\alpha\varphi(c_Q+s(u-c_Q))\,ds.
\end{equation}
For each $|\alpha|=p$, define a signed measure
$\Theta_{Q,\alpha}$ by
\[
\int\psi\,d\Theta_{Q,\alpha}
:=\frac{(-1)^pp}{\alpha!}
\int_Q\int_0^1(1-s)^{p-1}(u-c_Q)^\alpha
\psi(c_Q+s(u-c_Q))\,ds\,d\Delta_Q(u)
\]
for $\psi\in C_c(\mathbb R^d)$. Since
$P_{Q,\varphi}\in\mathcal P_{p-1}$,
\eqref{eq:periodic-cell-vanishing} and
\eqref{eq:periodic-taylor-remainder} give, as distributions,
\begin{equation}\label{eq:periodic-distribution-decomposition}
\Delta_Q=\sum_{|\alpha|=p}\partial^\alpha\Theta_{Q,\alpha},
\qquad
\operatorname{supp}\Theta_{Q,\alpha}\subset\overline Q,
\qquad
\|\Theta_{Q,\alpha}\|_{\rm TV}\leq C_{d,p}qL^p.
\end{equation}

Extend $\sum_{Q\in\mathscr C_{L,m}}\Delta_Q$ and
$\sum_{Q\in\mathscr C_{L,m}}\Theta_{Q,\alpha}$ periodically
under $S\mathbb Z^d$, and denote the extensions by $\Delta$
and $\Theta_\alpha$. For an $S\mathbb Z^d$-periodic signed
measure $\nu$, use the Fourier coefficients
\[
a_\nu(k):=S^{-d}\int_{[0,S)^d}
 e^{-2\pi i\langle k,u\rangle/S}\,d\nu(u)
\qquad(k\in\mathbb Z^d).
\]
Taking Fourier coefficients in
\eqref{eq:periodic-distribution-decomposition} yields
\begin{equation}\label{eq:periodic-fourier-vanishing}
|a_\Delta(k)|^2\leq C_{d,p}\|k/S\|^{2p}
\sum_{|\alpha|=p}|a_{\Theta_\alpha}(k)|^2.
\end{equation}
Also $a_\Delta(0)=0$ by
\eqref{eq:periodic-cell-vanishing}.

We next bound the squared Fourier coefficients of $\Delta$ and
$\Theta_\alpha$ in balls of radius one, uniformly in $S$.
Choose $\chi\in C_c^\infty(\mathbb R^d)$ such that
$|\widehat\chi|\geq1$ on $B_1$.
By \eqref{eq:periodic-cell-vanishing} and
\eqref{eq:periodic-distribution-decomposition}, a translate of
$\operatorname{supp}\chi$ intersects at most $C_{\chi,L}$ cells, and
therefore
\[
\sup_{t\in\mathbb R^d}\int|\chi(u-t)|\,d|\nu|(u)
\leq C_L
\qquad(\nu=\Delta,\Theta_\alpha).
\]
For $a\in\mathbb R^d$, Parseval on
$\mathbb R^d/S\mathbb Z^d$, applied to the $S\mathbb Z^d$-periodic
function
\[
t\longmapsto \int\chi(u-t)e^{-2\pi i\langle a,u-t\rangle}
\,d\nu(u),
\]
gives
\[
\begin{aligned}
\sum_{k\in\mathbb Z^d}|a_\nu(k)|^2
 |\widehat\chi(a-k/S)|^2
&=S^{-d}\int_{[0,S)^d}
 \left|\int\chi(u-t)e^{-2\pi i\langle a,u-t\rangle}
 \,d\nu(u)\right|^2dt\\
&\leq C_L^2.
\end{aligned}
\]
Hence
\begin{equation}\label{eq:uniform-fourier-ball}
\sup_{a\in\mathbb R^d}
\sum_{k/S\in a+B_1}|a_\nu(k)|^2\leq C_L^2
\qquad(\nu=\Delta,\Theta_\alpha).
\end{equation}
The constant $C_L$ depends on $L,d,p$ and the fixed $\chi$,
but not on $m$, $S$, or $\Lambda$.

For $0<\|k/S\|<1$, the exponent $2p>d+1$ and
\eqref{eq:periodic-fourier-vanishing} give
\[
\sum_{0<\|k/S\|<1}
\frac{|a_\Delta(k)|^2}{\|k/S\|^{d+1}}
\leq C_{d,p}\sum_{|\alpha|=p}
 \sum_{\|k/S\|<1}|a_{\Theta_\alpha}(k)|^2
\leq C_L.
\]
For $n\geq1$, the annulus $n\leq\|\xi\|<n+1$
can be covered by $O_d(n^{d-1})$ balls of radius one.
Applying \eqref{eq:uniform-fourier-ball} with $\nu=\Delta$
therefore yields
\[
\sum_{\|k/S\|\geq1}
\frac{|a_\Delta(k)|^2}{\|k/S\|^{d+1}}
\leq C_L\sum_{n\geq1}n^{-2}<\infty.
\]
Combining the estimates, we obtain
\begin{equation}\label{eq:periodic-weighted-fourier-sum}
\sum_{k\ne0}\frac{|a_\Delta(k)|^2}
 {\|k/S\|^{d+1}}\leq K_L,
\end{equation}
where $K_L$ is independent of $m$ and $\Lambda$.

By the definition of $\Delta$, its periodic extension is
$\delta_{\Lambda^{\rm per}}-\rho\Vol_d$.
Parseval now gives, with $t$ as the translation variable,
\[
\begin{aligned}
&S^{-d}\int_{[0,S)^d}
\left|N_{B_R}(\delta_{\Lambda^{\rm per}-t})
-\rho\Vol_d(B_R)\right|^2\,dt\\
&\qquad=\sum_{k\ne0}|a_\Delta(k)|^2
 |\widehat{\mathbf1}_{B_R}(k/S)|^2.
\end{aligned}
\]
For $R\geq1$ and $\xi\ne0$, the Fourier transform of the
Euclidean ball satisfies
\[
|\widehat{\mathbf1}_{B_R}(\xi)|^2
\leq C_dR^{d-1}\|\xi\|^{-(d+1)};
\]
see \cite[Sec.~IV and Appendix~A]{TS03}.
The assertion follows from
\eqref{eq:periodic-weighted-fourier-sum}.
\end{proof}

\subsubsection{Construction of a low-variance realization}

Fix a realization $\kappa$. We construct a nearby realization
by changing its points inside disjoint $S$-cubes. In the core of each
such cube, every $L$-cell will have exactly $q$ points and satisfy
\eqref{eq:rescaled-completed-moments}; the points in the outer shell
will be chosen so that the total point count in the $S$-cube is
unchanged. We then apply \cref{lem:uniform-periodic-variance}
to the configurations on the cores and estimate the local $L^2$
change separately.

Let \(a\in\Ad_d(X,\mu)\) and
\(\kappa\in\smash{\RR^a_{\rho,2}(X,\mu)}\). By
\cref{lem:concrete-version} and essential freeness, take a pointwise
equivariant Borel representative on an invariant conull Borel set
\(X_0\) on which \(a\) is free. Put
\(\Lambda_x:=\operatorname{supp}\kappa(x)\) for \(x\in X_0\).

\smallskip\noindent\emph{Exact moment configurations in \(L\)-cells.}
Fix an integer \(p\geq1\) with \(2p>d+1\) and a basis
\(P_1,\ldots,P_m\) of \(\mathcal P_{p-1}^0\).
Choose a closed cube \(D\subset(0,1)^d\), and set
\(\vartheta:=\Vol_d(D)\). Let \(\varepsilon_0,q_0,\sigma_d\) be the
constants of \cref{prop:relative-moment-repair} for these choices, and
put \(\varepsilon:=\varepsilon_0/2\).

For each integer \(q\geq q_0\), set \(L:=(q/\rho)^{1/d}\).
For every \(v\in\mathbb R^d\), put \(P:=v+LI\) and
\(D_P:=v+LD\). Recall that \(F_{L,P}(\Lambda)\) is the normalized
exterior point set defined in \eqref{eq:rescaled-exterior-points}.
We call \(P\) \emph{admissible} for \(\Lambda\) if
\begin{equation}\label{eq:admissible-L-cell}
F_{L,P}(\Lambda)\in\mathscr F_{q,\varepsilon},
\qquad
\#(\Lambda\cap D_P)\leq2\rho\vartheta L^d.
\end{equation}
Define
\[
\mathcal E_L:=\{x\in X_0:LI\text{ is not admissible for }\Lambda_x\}.
\]
Apply \cref{lem:weighted-cell-ergodic} with \(B=I\) to
\(A=I\setminus D\) with weights \(1,P_1,\ldots,P_m\), and to
\(A=D\) with weight \(1\). Since the expected count on
\(\partial(LD)\) is zero, we may disregard points there.
With tolerance \(\rho\min\{\varepsilon,\vartheta\}\), any
\(x\in\mathcal E_L\) outside the boundary-null set violates at
least one of these weighted-count inequalities, by the conditions in
\eqref{eq:moment-exterior-hypotheses} and
\eqref{eq:admissible-L-cell}.
The finite-family estimate \eqref{eq:weighted-cell-finite-family} therefore gives
\[
\mu(\mathcal E_L)\longrightarrow0,
\qquad
L^{-d}\int_{\mathcal E_L}N_{LI}(\kappa(x))\,d\mu(x)
\longrightarrow0
\qquad(q\to\infty).
\]

Put \(h_q:=q^{-1/d}\), and choose a translate \(\Gamma_q\) of
\(h_q\mathbb Z^d\) whose points are disjoint from
\(\partial I\cup\partial D\). Define
\[
F_q^0:=\{u\in\Gamma_q\cap(I\setminus D):
             \operatorname{dist}(u,D)\geq h_q/4\}.
\]
For each \(u\in\Gamma_q\), put \(J_u:=u+[0,h_q)^d\).
These cubes partition \(\mathbb R^d\). If
\(\mathbf1_{F_q^0}(u)\ne\mathbf1_{I\setminus D}(y)\) for some
\(y\in J_u\), then \(J_u\) is contained in the
\((\sqrt d+1)h_q\)-neighbourhood of \(\partial I\cup\partial D\).
The union of these cubes has volume at most \(K_{d,D}h_q\),
since both boundaries consist of finitely many flat faces.
For each $u\in F_q^0$ for which $J_u\subset I\setminus D$,
the mean-value inequality and $\|y-u\|\leq\sqrt d\,h_q$
for $y\in J_u$ give, for $1\leq i\leq m$,
\[
\left|h_q^dP_i(u)-\int_{J_u}P_i(y)\,dy\right|
\leq\sqrt d\,h_q^{d+1}
\sup_{y\in[-1,2]^d}\|\nabla P_i(y)\|.
\]
Every cube not covered by this estimate for which either its
base point lies in $F_q^0$ or the cube intersects $I\setminus D$
lies in the $(2\sqrt d+2)h_q$-neighbourhood of
$\partial I\cup\partial D$. The total volume of these cubes
is $O_{d,D}(h_q)$, and their contributions to both the Riemann sum
and the integral are bounded by this volume times
$\sup_{[-1,2]^d}|P_i|$. Combining these bounds, there is a
constant $K_{\mathrm{lat}}$, depending on
\(d,D,P_1,\ldots,P_m\), such that
\[
\bigl|h_q^d\#F_q^0-\Vol_d(I\setminus D)\bigr|
\leq K_{\mathrm{lat}}h_q,
\]
and
\[
\left|h_q^d\sum_{u\in F_q^0}P_i(u)
 -\int_{I\setminus D}P_i(u)\,du\right|
\leq K_{\mathrm{lat}}h_q
\qquad(1\leq i\leq m).
\]
Multiplying by \(q=h_q^{-d}\) gives an error at most
\(K_{\mathrm{lat}}q^{1-1/d}\).
For every \(q\) sufficiently large that
\(K_{\mathrm{lat}}q^{-1/d}\leq\varepsilon\),
we have \(F_q^0\in\mathscr F_{q,\varepsilon}\) and
\[
\operatorname{sep}(F_q^0)\geq h_q,
\qquad
\operatorname{dist}(F_q^0,D)\geq h_q/4.
\]
Set
\[
\Lambda_L^0:=L\bigl(F_q^0\cup A(q,F_q^0)\bigr)\subset LI.
\]
The configuration \(\Lambda_L^0\) has \(q\) points and the exact
polynomial moments required in \cref{lem:uniform-periodic-variance}.
The estimates of \cref{prop:relative-moment-repair} and the two
preceding separation bounds imply
\[
\operatorname{sep}(\Lambda_L^0)
\geq\sigma_*\rho^{-1/d},\qquad
\sigma_*:=\min\{\sigma_d,1/4\}.
\]
Let \(K_L\) be the constant in
\cref{lem:uniform-periodic-variance}.

\smallskip\noindent\emph{Large cubes and their cores.}
For \(0<\zeta<1/4\) and integers \(\ell>\ell'\) with
\(\ell-\ell'\) even, put \(S:=\ell L\) and \(S':=\ell' L\).
Suppose that
\[
\zeta\leq1-(S'/S)^d\leq2\zeta.
\]
In particular, $S'\geq S/2$.
Write \(\mathbf1:=(1,\ldots,1)\).  For \(Q=u+SI\), define
\[
\begin{aligned}
Q^{\mathrm{in}}&:=u+\tfrac{S-S'}2\mathbf1+S'I,
&Q^{\mathrm{out}}&:=Q\setminus Q^{\mathrm{in}},\\
\mathscr L(Q)&:=\left\{u+\tfrac{S-S'}2\mathbf1+Lk+LI:
  k\in\{0,\ldots,\ell'-1\}^d\right\},
&V_{\mathrm{out}}&:=S^d-(S')^d.
\end{aligned}
\]
We call $Q$ an $S$-cube, $Q^{\mathrm{in}}$ its core,
$Q^{\mathrm{out}}$ its outer shell, and the members of
$\mathscr L(Q)$ the $L$-cells of its core.
Since \(\ell-\ell'\) is even,
the cells in \(\mathscr L(Q)\) belong to the tiling \(\{v+LI:v\in u+L\mathbb Z^d\}\).
They partition \(Q^{\mathrm{in}}\), and
\(M:=\#\mathscr L(Q)=(S'/L)^d\).

For \(r\geq S^2\), choose an $r$-separated cocompact Borel
cross-section $C_r\subset X_0$ and its Voronoi cells $W_c$
as defined in \eqref{eq:labelled-voronoi}.
The relation
$\{(x,u)\in X_0\times\mathbb R^d:u.x\in C_r\}$ is Borel
and has countable sections. Lusin--Novikov
\cite[Theorem~18.10]{Kec95} gives Borel sets
$D_j\subset X_0$ and Borel maps $u_j:D_j\to\mathbb R^d$
whose disjoint graphs cover the relation; these maps will be used
below to establish measurability of the modified realization.
Let $\nu_r$ be the transverse measure normalized by the
cross-section formula
\cite[Theorem~1.14]{BHK25},
\[
\int_XH\,d\mu
=
\int_{C_r}\int_{W_c}H(t.c)\,dt\,d\nu_r(c).
\]
The cross-section mass-transport formula \cite[Theorem~1.14]{BHK25}
gives
\[
\int_{C_r}N_{W_c}(\kappa(c))\,d\nu_r(c)=\rho.
\]
For $s>0$, set
\[
W_c^{-s}:=\{u\in W_c:\operatorname{dist}(u,\partial W_c)>s\}.
\]
The inclusion \(\operatorname{int}B_{r/2}\subset W_c\) and convexity give, for
\(0<s<r/2\),
\[
\frac{\Vol_d(W_c\setminus W_c^{-s})}{\Vol_d(W_c)}
\leq\frac{2ds}{r}.
\]
For every $z\in[0,S)^d$, each point not covered by the cubes
$u+SI\subset W_c$, $u\in z+S\mathbb Z^d$, lies within distance
$\sqrt d\,S$ of $\partial W_c$.  Consequently,
\begin{equation}\label{eq:voronoi-grid-cover}
\Vol_d\!\left(W_c\setminus
\bigcup_{\substack{u\in z+S\mathbb Z^d\\u+SI\subset W_c}}(u+SI)\right)
\leq C_d\frac Sr\Vol_d(W_c).
\end{equation}

For \(z\in[0,S)^d\) and \(c\in C_r\), let
\[
\mathcal Q_z(c):=\{u+SI:u\in z+S\mathbb Z^d,\ u+SI\subset W_c\}.
\]
For $Q=u+[0,S)^d$ and $w\in\mathbb R^d$,
\begin{equation}\label{eq:cube-linear-supremum}
\sup_{t\in Q}\langle t,w\rangle
=\langle u,w\rangle+S\sum_{i=1}^d\max\{w_i,0\}.
\end{equation}
This supremum is attained on $Q$ exactly when every $w_i\leq0$.
The inclusion of $Q$ in a strict or non-strict affine half-space
is therefore a Borel condition.  For each $j$ with $c\in D_j$,
the half-space defined by $u_j(c)$ depends Borelly on $c$;
intersecting these countably many conditions shows that
$(c,z)\mapsto\mathcal Q_z(c)$ is Borel.

For \(Q\in\mathcal Q_z(c)\), put
\[
n_Q:=\#(\Lambda_c\cap Q),\qquad k_Q:=n_Q-qM.
\]
We call \(Q\) \emph{count-compatible} if
\begin{equation}\label{eq:count-compatible-S-cube}
|n_Q-\rho S^d|\leq\frac{\rho}{4}V_{\rm out}.
\end{equation}
\smallskip\noindent\emph{Replacement in each large cube.}
For \(Q\in\mathcal Q_z(c)\), the definitions of \(n_Q\), \(k_Q\), and
\(M\) give
\[
k_Q=\rho V_{\rm out}+(n_Q-\rho S^d).
\]
If \(Q\) is count-compatible, then \eqref{eq:count-compatible-S-cube}
gives
\[
\frac34\rho V_{\rm out}
\leq k_Q\leq
\frac54\rho V_{\rm out}.
\]
Write \(Q=u+SI\), and define the finite set of lattice sites
\[
\mathcal G_Q:=\bigl(u+(2\rho)^{-1/d}\mathbb Z^d\bigr)
\cap Q^{\mathrm{out}}.
\]
It has \(2\rho V_{\rm out}+O(S^{d-1})\) points.  Since
\(\zeta S^d\leq V_{\rm out}\leq2\zeta S^d\), for all sufficiently large
\(S\),
\[
\#\mathcal G_Q\geq\tfrac54\rho V_{\rm out}\geq k_Q
\qquad(Q\text{ count-compatible}).
\]
For a count-compatible \(Q\), let \(\mathcal G_Q(k_Q)\) be the
first \(k_Q\) points of \(\mathcal G_Q\) in lexicographic order.

Recall that an $L$-cell $P$ is admissible for $\Lambda_c$
precisely when $F_{L,P}(\Lambda_c)\in\mathscr F_{q,\varepsilon}$
and $\#(\Lambda_c\cap D_P)\leq2\rho\vartheta L^d$.
For each \(P=v+LI\in\mathscr L(Q)\), define
\[
\mathcal T_P(\Lambda_c):=
\begin{cases}
\bigl(\Lambda_c\cap(P\setminus\operatorname{int}D_P)\bigr)
 \cup\bigl(v+L A(q,F_{L,P}(\Lambda_c))\bigr),
 &P\text{ admissible for }\Lambda_c,\\
v+\Lambda_L^0,&P\text{ not admissible for }\Lambda_c.
\end{cases}
\]
By \eqref{eq:rescaled-completed-moments}, the first case has the required
polynomial moments in \(P\).  The definition of \(\Lambda_L^0\) gives
the same identities in the second case.  Thus, for each \(P\in\mathscr L(Q)\),
\begin{equation}\label{eq:modified-core-cell-moments}
\sum_{u\in\mathcal T_P(\Lambda_c)}H(u)=\rho\int_P H(u)\,du
\qquad(H\in\mathcal P_{p-1}).
\end{equation}
In particular, \(\#\mathcal T_P(\Lambda_c)=q\).

For each \(z\in[0,S)^d\), define \(\Lambda_{c,z}'\) on every
\(Q\in\mathcal Q_z(c)\) by
\[
\Lambda_{c,z}'\cap Q:=
\begin{cases}
\displaystyle\bigcup_{P\in\mathscr L(Q)}\mathcal T_P(\Lambda_c)
 \ \cup\ \mathcal G_Q(k_Q), &Q\text{ count-compatible},\\
\Lambda_c\cap Q,&Q\text{ not count-compatible}.
\end{cases}
\]
Outside these \(S\)-cubes, put \(\Lambda_{c,z}'=\Lambda_c\).
The union defining \(\Lambda_{c,z}'\cap Q\) is disjoint.  For each
count-compatible \(Q\),
\[
\#(\Lambda_{c,z}'\cap Q)=qM+k_Q=n_Q.
\]
Put $\Lambda^{W}_{c,z}:=\Lambda'_{c,z}\cap W_c$.  The global
realization is the counting measure with support
\[
\operatorname{supp}\kappa_z(x)
:=\bigcup_{u\in(C_r)_x}\bigl(u+\Lambda^{W}_{u.x,z}\bigr).
\]
The sets $u+W_{u.x}$ partition $\mathbb R^d$ and are locally
finite, so the displayed union is locally finite. The maps
$c\mapsto\Lambda^W_{c,z}$ are Borel. The inclusion of cubes in
$\mathcal Q_z(c)$ is Borel by \eqref{eq:cube-linear-supremum}
and the calculation following it; count-compatibility is the
inequality \eqref{eq:count-compatible-S-cube}, and admissibility
is determined by the two conditions in its definition. The
moment-completion map of \cref{prop:relative-moment-repair} and the
lexicographic choice from $\mathcal G_Q$ are Borel.  Using the return maps $u_j$,
for each bounded Borel set $D$ we have
\[
N_D(\kappa_z(x))
=\sum_{j:x\in D_j}
  N_{D-u_j(x)}(\delta_{\Lambda^W_{u_j(x).x,z}}).
\]
The summands are Borel, so $x\mapsto\kappa_z(x)$ is Borel.
For $v\in\mathbb R^d$,
\[
(C_r)_{v.x}=(C_r)_x-v,
\qquad (u-v).(v.x)=u.x,
\]
which gives $\kappa_z(v.x)=\theta_v\kappa_z(x)$.
The equalities $\#(\Lambda'_{c,z}\cap Q)=\#(\Lambda_c\cap Q)$
for modified cubes and $\Lambda'_{c,z}=\Lambda_c$ elsewhere give
\[
\#\Lambda^W_{c,z}=\#(\Lambda_c\cap W_c),\qquad
\int_{C_r}\#\Lambda^W_{c,z}\,d\nu_r(c)=\rho.
\]
Thus $\kappa_z$ is a translation-equivariant simple counting-measure
realization of intensity $\rho$.

For each $z\in[0,S)^d$, define
\[
\kappa_z^+(x):=\delta_{\operatorname{supp}\kappa_z(x)
               \setminus\operatorname{supp}\kappa(x)}.
\]

\begin{lemma}[Uniform inserted-point bound]\label{lem:lowvar-inserted}
For every $n\geq1$ there is a constant $D_n<\infty$,
depending only on $d,\rho,n$, such that, for every choice of
$D,p,q,L,\vartheta,\gamma,\zeta,S,r$ satisfying the preceding hypotheses and every $z\in[0,S)^d$,
\begin{equation}\label{eq:low-variance-inserted-count}
N_{Q_n}(\kappa_z^+(x))\leq D_n\qquad(x\in X_0).
\end{equation}
In particular, $\kappa_z\in\RR^a_{\rho,2}(X,\mu)$.
\end{lemma}

\begin{proof}
An $L$-cell $P$ is admissible here if
$F_{L,P}(\Lambda_c)\in\mathscr F_{q,\varepsilon}$ and
$\#(\Lambda_c\cap D_P)\leq2\rho\vartheta L^d$.
The points inserted in an admissible $L$-cell, in an $L$-cell
where either inequality fails, and in an outer shell have
respective separation bounds
$\sigma_d\rho^{-1/d}$, $\sigma_*\rho^{-1/d}$, and $(2\rho)^{-1/d}$.
Set $s_0:=\rho^{-1/d}\min\{\sigma_*,2^{-1/d}\}$.
The core $L$-cells have disjoint interiors, as do the $S$-cubes
containing the outer shells.  Disjoint axis-parallel cubes of side
$T$ intersecting $Q_n$ are contained in its coordinatewise $T$-enlargement,
so their number is at most $(n/T+2)^d$.  Since
$S\geq L\geq\rho^{-1/d}$, at most
$M_n:=(\lceil n\rho^{1/d}\rceil+2)^d$ cubes of either type intersect $Q_n$.
Let $J_n$ bound the number of points in an $s_0$-separated subset
of $Q_n$.  Then \eqref{eq:low-variance-inserted-count} holds with
$D_n:=2M_nJ_n$.  The constants $\sigma_d,\sigma_*$ depend only on $d$.
Finally,
$N_{Q_n}(\kappa_z)\leq N_{Q_n}(\kappa)+D_n$.
Since $\kappa$ has local second moments and $\kappa_z$ has intensity
$\rho$, $\kappa_z\in\RR^a_{\rho,2}(X,\mu)$.
\end{proof}

For each \(z\in[0,S)^d\), define
\[
\begin{aligned}
e_1(z)&:=\int_{C_r}\sum_{\substack{Q\in\mathcal Q_z(c)\\
 Q\text{ not count-compatible}}}S^d\,d\nu_r(c),\\
e_2(z)&:=\int_{C_r}\sum_{Q\in\mathcal Q_z(c)}
 \sum_{\substack{P\in\mathscr L(Q)\\
 P\text{ not admissible for }\Lambda_c}}L^d\,d\nu_r(c),\\
e_3(z)&:=\int_{C_r}\sum_{Q\in\mathcal Q_z(c)}
 \sum_{\substack{P\in\mathscr L(Q)\\
 P\text{ not admissible for }\Lambda_c}}\#(\Lambda_c\cap P)\,d\nu_r(c),\\
e_4(z)&:=\int_{C_r}\sum_{Q\in\mathcal Q_z(c)}
 \#(\Lambda_c\cap Q^{\mathrm{out}})\,d\nu_r(c).
\end{aligned}
\]
Set
\[
b_{S,\zeta}:=
\mu\left\{x:
\left|N_{SI}(\kappa(x))-\rho S^d\right|>
\frac{\rho\zeta S^d}{4}\right\}.
\]
By \cref{lem:weighted-cell-ergodic}, \(b_{S,\zeta}\to0\) as \(S\to\infty\)
for each fixed \(\zeta>0\).

\begin{lemma}[Choice of grid shift]\label{lem:lowvar-shift}
Let $\gamma>0$ and assume
\[
\mu(\mathcal E_L)<\gamma,\qquad
L^{-d}\int_{\mathcal E_L}N_{LI}(\kappa(x))\,d\mu(x)<\gamma.
\]
Then there is $z\in[0,S)^d$ such that
\[
e_1(z)\leq4b_{S,\zeta},\qquad
e_2(z)+e_3(z)\leq8\gamma,\qquad
e_4(z)\leq8\rho\zeta.
\]
\end{lemma}

\begin{proof}
For \(T=S\) or \(L\), let
\(u_{T,z}(t)\in[0,T)^d\) be the position of
\(t\) relative to the lower corner of its cell in the tiling
\(\{v+[0,T)^d:v\in z+T\mathbb Z^d\}\).
For every nonnegative Borel \(H\) on \(X\), the cross-section formula,
translation invariance of \(\mu\), and averaging \(z\) modulo \(T\) give
\begin{equation}\label{eq:shift-average}
\frac1{T^d}\int_{[0,T)^d}\int_{C_r}\int_{W_c}
H((t-u_{T,z}(t)).c)\,dt\,d\nu_r(c)\,dz
=\int_XH\,d\mu.
\end{equation}
Let $\mathcal P_{T,z}$ be the grid of half-open $T$-cubes with
lower corners in $z+T\mathbb Z^d$, and let $a(P)$ denote the
lower corner of $P$. For every nonnegative Borel $H$,
\begin{equation}\label{eq:restricted-grid-average}
\sum_{\substack{P\in\mathcal P_{T,z}\\P\subset W_c}}
T^dH(a(P).c)
\leq\int_{W_c}H((t-u_{T,z}(t)).c)\,dt.
\end{equation}
The integrand equals $H(a(P).c)$ on $P$.
The inequality also holds when the sum is restricted to fewer cells.
For $e_2,e_3$, the lower corner of a core $L$-cell differs
from that of its containing $S$-cube by
$\frac{\ell-\ell'}2L\mathbf1+Lk$, $k\in\mathbb Z^d$.
Since $\ell-\ell'$ is even, core $L$-cells belong to the
$z+L\mathbb Z^d$ grid. Applying \eqref{eq:shift-average}
with $T=S$ and
\[
H(x):=\mathbf1_{\{|N_{SI}(\kappa(x))-\rho S^d|
>\rho\zeta S^d/4\}}
\]
gives the bound for $e_1$. Since $S=\ell L$, averaging $z$
over $[0,S)^d$ also averages its residue modulo $L$ uniformly.
Taking $T=L$ and then
\[
H(x):=\mathbf1_{\mathcal E_L}(x),\qquad
H(x):=L^{-d}N_{LI}(\kappa(x))\mathbf1_{\mathcal E_L}(x)
\]
gives the bounds for $e_2$ and $e_3$. Thus
\[
\begin{aligned}
\frac1{S^d}\int_{[0,S)^d}e_1(z)\,dz&\leq b_{S,\zeta},\\
\frac1{S^d}\int_{[0,S)^d}e_2(z)\,dz&\leq\mu(\mathcal E_L)<\gamma,\\
\frac1{S^d}\int_{[0,S)^d}e_3(z)\,dz
&\leq L^{-d}\int_{\mathcal E_L}N_{LI}(\kappa(x))\,d\mu(x)<\gamma.
\end{aligned}
\]
For fixed \(c\) and \(u\in\Lambda_c\cap W_c\), the proportion of shifts
\(z\in[0,S)^d\) for which \(u\) belongs to the complement of a
core in the full \(S\)-grid is \(1-(S'/S)^d\).
Integrating first in \(z\) and then in \(c\), the mass-transport formula gives
\[
\frac1{S^d}\int_{[0,S)^d}e_4(z)\,dz
\leq\rho\bigl(1-(S'/S)^d\bigr)\leq2\rho\zeta.
\]
Markov's inequality gives a shift \(z\in[0,S)^d\) satisfying
\[
e_1(z)\leq4b_{S,\zeta},\qquad
e_2(z)+e_3(z)\leq8\gamma,\qquad
e_4(z)\leq8\rho\zeta.
\]
Indeed, each of the three strict reverse inequalities holds on a set of
shifts of measure at most \(S^d/4\).
\end{proof}

\begin{lemma}[Rooted approximation]\label{lem:lowvar-rooted}
For every $\varepsilon_*>0$ there is
$\delta_*=\delta_*(\kappa,\varepsilon_*)>0$ with the following
property. Assume that $q\ge q_0$,
$K_{\mathrm{lat}}q^{-1/d}\leq\varepsilon$, $0<\zeta<1/4$, $\ell>\ell'$ with
$\ell-\ell'$ even, $\zeta\leq1-(S'/S)^d\leq2\zeta$,
$r\ge S^2$, and
$\#\mathcal G_Q\geq\tfrac54\rho V_{\mathrm{out}}$
for every $S$-cube $Q=u+SI$.
Suppose that the hypotheses of \cref{lem:lowvar-shift}
hold and $z$ satisfies its three inequalities. If
\begin{equation}\label{eq:low-variance-rooted-threshold}
\vartheta+\gamma+\zeta+\frac Sr<\delta_*,
\end{equation}
then $d_{\mathrm{root},2}(\kappa_z,\kappa)<\varepsilon_*$.
The choice of $\delta_*$ is independent of
$p,D,q,L,S,r,z,\vartheta,\gamma,\zeta$, the cross-section $C_r$,
and its Borel order.
\end{lemma}

\begin{proof}
Put $\kappa':=\kappa_z$, $\Lambda_c':=\Lambda'_{c,z}$,
$\kappa^+:=\kappa_z^+$, and
\[
\kappa^-(x):=\delta_{\operatorname{supp}\kappa(x)
                         \setminus\operatorname{supp}\kappa'(x)}.
\]
Write $N_n^\pm:=N_{Q_n}\circ\kappa^\pm$ and let
$\rho_\pm$ be the intensities of $\kappa^\pm$.  The numbers of points
removed and inserted in an admissible core cell $P$ satisfy
\[
\#((\Lambda_c\setminus\Lambda_c')\cap P)
\leq2\rho\vartheta L^d,
\qquad
\#((\Lambda_c'\setminus\Lambda_c)\cap P)
\leq\tfrac32\rho\vartheta L^d.
\]
For a core $L$-cell $P$ that fails either of the admissibility
requirements $F_{L,P}(\Lambda_c)\in\mathscr F_{q,\varepsilon}$ or
$\#(\Lambda_c\cap D_P)\leq2\rho\vartheta L^d$,
the total number of removed and inserted points is at most
$\#(\Lambda_c\cap P)+\rho L^d$. Summing over such cells gives
$e_3(z)+\rho e_2(z)$.  A modified $S$-cube contributes at most
$\#(\Lambda_c\cap Q^{\mathrm{out}})+5\rho V_{\mathrm{out}}/4$
from its outer shell. Since $V_{\mathrm{out}}/S^d\leq2\zeta$,
\cref{lem:lowvar-shift} and the cross-section mass-transport identity give
\[
\begin{aligned}
\rho_-+\rho_+
&\leq\tfrac72\rho\vartheta+e_3(z)+\rho e_2(z)
+e_4(z)+\tfrac52\rho\zeta\\
&\leq\tfrac72\rho\vartheta
+8\max\{1,\rho\}\gamma+\tfrac{21}{2}\rho\zeta.
\end{aligned}
\]
In particular, for every fixed $n$,
\[
\int_X(N_n^-+N_n^+)\,d\mu
\leq C_n(\vartheta+\gamma+\zeta).
\]

By \cref{lem:lowvar-inserted},
$N_{Q_n}(\kappa')\leq N_{Q_n}(\kappa)+D_n$.

It follows that
\[
(N_n^-)^2\leq N_{Q_n}(\kappa)^2\mathbf1_{\{N_n^->0\}},
\qquad (N_n^+)^2\leq D_nN_n^+.
\]
As $\vartheta+\gamma+\zeta\to0$ over the allowed parameter choices,
the first squared $L^2$-norm tends to zero by absolute continuity,
since
$\mu\{N_n^->0\}\leq\int_XN_n^-\,d\mu\to0$.
The second squared norm tends to zero by the preceding $L^1$ bound.
Thus $\|N_n^-\|_2+\|N_n^+\|_2\to0$ for each fixed $n$.
For
$\varphi\in C_c(\mathbb R^d)$ vanishing outside $Q_n$,
\[
|S_\varphi(\kappa')-S_\varphi(\kappa)|
\leq\|\varphi\|_\infty(N_n^-+N_n^+).
\]
The same estimate with $\varphi=\mathbf1_{Q_n}$ gives
$L^2$-convergence of each local count.

For any sequence of allowed parameter choices with
$\vartheta+\gamma+\zeta+S/r\to0$, the probability that
$\kappa'$ and $\kappa$ disagree on $Q_n$ tends to zero for
all $n$.  Every subsequence has a further subsequence along which
these restrictions eventually agree almost everywhere, simultaneously
for every $n$.  On that subsequence the configurations converge
vaguely almost everywhere, hence in measure.  Thus $d_0\to0$ and
$d_{\mathrm{root},2}(\kappa',\kappa)\to0$.
If the asserted $\delta_*$ did not exist, a sequence of allowed
choices with the left side of \eqref{eq:low-variance-rooted-threshold}
tending to zero and rooted distance at least $\varepsilon_*$
would contradict this convergence.
\end{proof}

\begin{lemma}[Variance on modified cores]\label{lem:lowvar-core}
Assume that $\#\mathcal G_Q\geq\tfrac54\rho V_{\mathrm{out}}$
for every $S$-cube $Q=u+SI$. Under the hypotheses of
\cref{lem:lowvar-shift}, let $z$ satisfy its
three inequalities, let $\kappa':=\kappa_z$, and fix an integer
$R\geq1$ with $S>4R$. Define
\[
E_R:=\bigl\{t.c:\ c\in C_r,\ Q\in\mathcal Q_z(c),\
 Q\text{ count-compatible},\ t\in Q^{\mathrm{in}},\
 \operatorname{dist}(t,\partial Q^{\mathrm{in}})>R\bigr\}.
\]
Then
\[
\int_{E_R}|N_{B_R}(\kappa')-\rho\Vol_d(B_R)|^2\,d\mu
\leq K_LR^{d-1},
\]
and
\[
\mu(E_R^c)\leq4b_{S,\zeta}+C_d\zeta+C_dR/S+C_dS/r.
\]
There is an integrable function
$C_R(1+N_{B_R}(\kappa)^2)$, independent of the later choices
of $\zeta,S,r,z$, that bounds
$|N_{B_R}(\kappa')-\rho\Vol_d(B_R)|^2$ on $E_R^c$.
\end{lemma}

\begin{proof}
Put $\Lambda'_c:=\Lambda'_{c,z}$.  For each count-compatible
\(Q\in\mathcal Q_z(c)\), define the
periodic point set
\[
\Lambda_Q^{\rm per}:=
 \bigcup_{v\in S'\mathbb Z^d}\bigl((\Lambda'_c\cap Q^{\mathrm{in}})+v\bigr).
\]
Each \(L\)-cell in \(Q^{\mathrm{in}}\) satisfies the exact polynomial
moment identities, so \cref{lem:uniform-periodic-variance} applies to
\(\Lambda_Q^{\rm per}\).
For \(x=t.c\in E_R\), the ball \(B_R+t\) lies in the core.  Equivariance
gives
\[
N_{B_R}(\kappa'(x))
=\#\bigl(\operatorname{supp}\kappa'(c)\cap(B_R+t)\bigr),
\]
and \(\Lambda_Q^{\rm per}\) agrees with
\(\operatorname{supp}\kappa'(c)\) on \(B_R+t\).
For fixed \(c,Q\), enlarge the domain of integration from
\(\{t\in Q^{\mathrm{in}}:\operatorname{dist}(t,\partial Q^{\mathrm{in}})>R\}\)
to \(Q^{\mathrm{in}}\).  The integrand is nonnegative, and the latter
cube is a fundamental domain for \(\Lambda_Q^{\rm per}\).
By \cref{lem:uniform-periodic-variance}, the integral is at most
\(K_LR^{d-1}S'^d\).
The cores are disjoint.  The cross-section formula gives
\[
\int_{C_r}\sum_{Q}(S')^d\,d\nu_r(c)
\leq
\int_{C_r}\Vol_d(W_c)\,d\nu_r(c)=1,
\]
where the sum is over modified $S$-cubes $Q\subset W_c$.  Hence
\[
\int_{E_R}
\left|N_{B_R}(\kappa'(x))-\rho\Vol_d(B_R)\right|^2\,d\mu(x)
\leq K_LR^{d-1}.
\]

If \(x=t.c\notin E_R\), then either $t$ lies in an $S$-cube
$Q\in\mathcal Q_z(c)$ for which \eqref{eq:count-compatible-S-cube}
fails, or $t$ lies in an outer shell, within distance $R$ of the
boundary of a core, or in the part of $W_c$ not covered by
$\mathcal Q_z(c)$.
By \eqref{eq:voronoi-grid-cover} and
$e_1(z)\leq4b_{S,\zeta}$ from \cref{lem:lowvar-shift},
\[
\mu(E_R^c)
\leq 4b_{S,\zeta}+C_d\zeta+C_d\frac RS+C_d\frac Sr.
\]
Since $B_R\subset Q_{2R+2}$, \eqref{eq:low-variance-inserted-count}
gives $N_{B_R}(\kappa_z^+)\leq D_{2R+2}$, independently of
$q,\vartheta,\gamma,\zeta,S,r$.
Therefore
\[
N_{B_R}(\kappa')\leq N_{B_R}(\kappa)+D_{2R+2},
\]
and, on $E_R^c$,
\[
\left|N_{B_R}(\kappa')-\rho\Vol_d(B_R)\right|^2
\leq C_R\bigl(1+N_{B_R}(\kappa)^2\bigr),
\]
where $C_R$ depends only on $R,d,\rho$ and the separation constants.
The function \(1+N_{B_R}(\kappa)^2\) is integrable.
\end{proof}

\begin{proof}[Proof of \cref{prop:low-variance-density}]
Fix an open neighbourhood $\mathcal U$ of $\kappa$ and
$\varepsilon_{\mathcal U}>0$ with
$B_{d_{\mathrm{root},2}}(\kappa,\varepsilon_{\mathcal U})\subset\mathcal U$.
Let $\delta_{\mathcal U}:=\delta_*(\kappa,\varepsilon_{\mathcal U})$
be given by \cref{lem:lowvar-rooted}.
Choose \(D\) and
\(\gamma>0\) so that
\[
\vartheta+\gamma<\frac{\delta_{\mathcal U}}4.
\]
With this \(D\), the constants
\(\varepsilon_0,q_0,\sigma_d\) are fixed.  Choose \(q\geq q_0\) so that
\[
\mu(\mathcal E_L)<\gamma,
\qquad
L^{-d}\int_{\mathcal E_L}N_{LI}(\kappa(x))\,d\mu(x)<\gamma,
\qquad
K_{\mathrm{lat}}q^{-1/d}\leq\varepsilon.
\]
This fixes \(L\) and \(K_L\).  Since \(\beta>(d-1)/d\), choose
an integer \(R>N\) such that
\[
2K_LR^{d-1}<j^{-1}\Vol_d(B_R)^\beta.
\]
For this $R$, the function
$C_R(1+N_{B_R}(\kappa)^2)$ in \cref{lem:lowvar-core}
is integrable. Choose $\delta_R>0$ such that
\[
\mu(E)<\delta_R\quad\Longrightarrow\quad
C_R\int_E(1+N_{B_R}(\kappa)^2)\,d\mu
<K_LR^{d-1}.
\]
Choose \(0<\zeta<1/4\) with
\[
\zeta<\frac{\delta_{\mathcal U}}4,
\qquad
C_d\zeta<\frac{\delta_R}{3}.
\]
For each sufficiently large integer \(\ell\), choose
\(\ell'<\ell\) of the same parity as \(\ell\) with
\[
\bigl|\ell'-\ell(1-3\zeta/2)^{1/d}\bigr|\leq2.
\]
Then \(\ell'/\ell\to(1-3\zeta/2)^{1/d}\), and
\[
\zeta\leq1-(\ell'/\ell)^d\leq2\zeta
\]
for all sufficiently large \(\ell\).  Since \(b_{S,\zeta}\to0\)
for fixed \(\zeta\) as \(S=\ell L\to\infty\), we may also require
\[
S>4R,
\qquad
4b_{S,\zeta}+C_d\frac RS<\frac{\delta_R}{3},
\qquad\#\mathcal G_Q\geq\tfrac54\rho V_{\mathrm{out}}
\quad(Q=u+SI).
\]
Finally choose \(r\geq S^2\) so that
\[
C_dS/r<\delta_R/3,\qquad
S/r<\delta_{\mathcal U}/4.
\]
Choose $z$ by \cref{lem:lowvar-shift} and set $\kappa':=\kappa_z$.  The bounds give
\(\mu(E_R^c)<\delta_R\), while
\cref{lem:lowvar-rooted} gives
\(d_{\mathrm{root},2}(\kappa',\kappa)<\varepsilon_{\mathcal U}\).
By the defining property of \(\delta_R\),
\[
\int_{E_R^c}
\left|N_{B_R}(\kappa')-\rho\Vol_d(B_R)\right|^2\,d\mu
\leq K_LR^{d-1}.
\]
Consequently \(\kappa'\in\mathcal U\), and
\[
\mathcal V_{\kappa'}(R)
\leq2K_LR^{d-1}
<j^{-1}\Vol_d(B_R)^\beta.
\]
\end{proof}

\subsection{High-variance approximation}\label{subsec:high-variance}

In this subsection we prove density of realizations with arbitrarily
large variance at prescribed scales. We preserve the number of points
in each large cube but move them into smaller subcubes. A random
translation of this arrangement produces a count fluctuation of order
the squared expected number of points in an observation ball.

\begin{proposition}[High-variance density]\label{prop:high-variance-density}
For every \(a\in\Ad_d(X,\mu)\), \(\alpha<2\), and \(j,N\geq1\), the set
\(\smash{\mathcal H^a_{\alpha,j,N}}\) is dense in
\(\smash{\RR^a_{\rho,2}(X,\mu)}\).
\end{proposition}

\begin{proof}
Fix \(\kappa\in\smash{\RR^a_{\rho,2}(X,\mu)}\) and an open neighbourhood
\(\mathcal U\) of \(\kappa\).  It suffices to find
\(\kappa'\in\mathcal U\cap\mathcal H^a_{\alpha,j,N}\).

\smallskip\noindent\emph{Definition of \(\kappa_{\delta,S,r}\).}
By \cref{lem:concrete-version} and essential freeness, use a pointwise equivariant Borel
representative of $\kappa$ on an invariant conull Borel set $X_0$
on which $a$ is free.  Write
\(\Lambda_x:=\operatorname{supp}\kappa(x)\).  For
\(0<\delta<6^{-d}\), \(S>0\), and \(r\geq S^2\), choose an
$r$-separated cocompact Borel cross-section $C_r\subset X_0$,
and equip it with the Voronoi cells $W_c$ of
\eqref{eq:labelled-voronoi}. Normalize the transverse measure $\nu_r$ by
\[
\int_X H\,d\mu
=
\int_{C_r}\int_{W_c}H(t.c)\,dt\,d\nu_r(c)
\]
for nonnegative Borel $H$.

With \(h:=\delta^{1/d}S\), for each \(c\in C_r\) put
\[
\mathcal Q_S(c):=\{v+[0,S)^d:v\in S\mathbb Z^d,\ v+[0,S)^d\subset W_c\}.
\]
By \eqref{eq:voronoi-grid-cover},
\[
\Vol_d\!\left(
W_c\setminus\bigcup_{Q\in\mathcal Q_S(c)}Q
\right)
\leq
C_d\frac Sr\,\Vol_d(W_c).
\]
The cross-section formula gives
\[
S^{-d}\left(1-C_d\frac Sr\right)
\leq\int_{C_r}\#\mathcal Q_S(c)\,d\nu_r(c)
\leq S^{-d}.
\]
For $Q=a_Q+[0,S)^d$, set $e:=(1,\ldots,1)$.  Since $h<S/6$,
the cubes
\[
\begin{aligned}
U_Q&:=a_Q+\frac S4e+[0,3h)^d,\\
H_Q&:=a_Q+\frac S4e+[0,h)^d,\\
D_Q&:=a_Q+\frac S4e+2he+[0,h)^d
\end{aligned}
\]
satisfy $H_Q\subset U_Q$ and $\overline{U_Q}\subset\operatorname{int}Q$, and
$D_Q\subset U_Q\setminus H_Q$.  Define
\[
T_Q(y):=a_Q+\frac S4e+2he+\frac13
\left(y-a_Q-\frac S4e\right),
\qquad y\in U_Q.
\]
For each $c\in C_r$, set
\[
U_c:=\bigcup_{Q\in\mathcal Q_S(c)}U_Q,\qquad
\widetilde\Lambda_c:=
\bigl(\Lambda_c\cap(W_c\setminus U_c)\bigr)
\cup\bigcup_{Q\in\mathcal Q_S(c)}T_Q(\Lambda_c\cap U_Q).
\]
Define
\[
\kappa_{\delta,S,r}(x):=
\sum_{u\in(C_r)_x}\ \sum_{v\in\widetilde\Lambda_{u.x}}\delta_{u+v},
\qquad x\in X_0.
\]
The map $c\mapsto\widetilde\Lambda_c$ is Borel. The condition
$Q\subset W_c$ is Borel by the half-space description of $W_c$
in \eqref{eq:labelled-voronoi} and the cube supremum formula
\eqref{eq:cube-linear-supremum}. Point restrictions and affine
images of locally finite configurations are also Borel.
The cells $u+W_{u.x}$ partition $\mathbb R^d$, and only finitely many
intersect a compact set.  Thus the sum defines a locally finite simple
counting measure.  The relation
$\{(x,u)\in X_0\times\mathbb R^d:u.x\in C_r\}$ is
Borel with countable sections, so Lusin--Novikov
\cite[Theorem~18.10]{Kec95} gives Borel partial maps
$u_j:D_j\to\mathbb R^d$ whose values enumerate $(C_r)_x$.
For each bounded Borel set $D$,
\[
N_D(\kappa_{\delta,S,r}(x))
=\sum_{j:x\in D_j}
  \#(\widetilde\Lambda_{u_j(x).x}\cap(D-u_j(x))),
\]
which is Borel because $c\mapsto\widetilde\Lambda_c$
is a Borel locally finite point-set map.  The return-time
identity $(C_r)_{v.x}=(C_r)_x-v$ gives equivariance.
Since $T_Q:U_Q\to D_Q\subset U_Q$
is injective and
\[
\#\widetilde\Lambda_c=N_{W_c}(\kappa(c)),
\]
the cross-section mass-transport formula gives intensity $\rho$.
Moreover,
\(\|T_Q(y)-y\|\leq C_dS\).  For every \(m\geq1\), choose
\(m'=m'(m,S)\) such that the \(C_dS\)-neighbourhood of \(Q_m\)
is contained in \(Q_{m'}\).  Injectivity of the modifications gives
\[
N_{Q_m}(\kappa_{\delta,S,r}(x))\leq N_{Q_{m'}}(\kappa(x))
\qquad(x\in X_0).
\]
Hence \(N_{Q_m}\circ\kappa_{\delta,S,r}\in L^2(\mu)\).

\smallskip\noindent\emph{Closeness in the rooted topology.}
Fix \(m\geq1\) and restrict to $S>\max\{1,6m,2\sqrt d\}$.  For
$Q=a_Q+[0,S)^d\in\mathcal Q_S(c)$ put
\[
a_U:=a_Q+\frac S4e,\qquad a_D:=a_U+2he,
\]
and define the closed cubes
\[
A_{Q,m}:=a_D+[-m/2,h+m/2]^d,\qquad
B_{Q,m}:=a_U+[-3m/2,3h+3m/2]^d.
\]
Since $h<S/6$ and $S>6m$, we have
$A_{Q,m},B_{Q,m}\subset\operatorname{int}Q\subset W_c$.  Moreover,
\[
\begin{aligned}
Q_m\cap(U_Q-t)\ne\varnothing
&\ \Longrightarrow\ t\in B_{Q,m},\\
Q_m\cap(D_Q-t)\ne\varnothing
&\ \Longrightarrow\ t\in A_{Q,m}.
\end{aligned}
\]
The map $(c,t)\mapsto t.c$ is a Borel bijection from
$\{(c,t):c\in C_r,\ t\in W_c\}$ onto $X_0$.  Define
\[
F_m:=\{t.c:c\in C_r,\ Q\in\mathcal Q_S(c),\
                         t\in B_{Q,m}\}.
\]
If $x\notin F_m$, the measures $\kappa_{\delta,S,r}(x)$ and
$\kappa(x)$ agree on $Q_m$.  The cross-section formula gives
\[
\mu(F_m)
=(3h+3m)^d\int_{C_r}\#\mathcal Q_S(c)\,d\nu_r(c)
\leq \left(3\delta^{1/d}+\frac{3m}{S}\right)^d
\leq C_d\delta+\frac{C_{d,m}}S.
\]

For $x=t.c$, with $t\in W_c$, define
\[
\begin{aligned}
I_m(x)&:=\sum_{Q\in\mathcal Q_S(c)}
\#\bigl(T_Q(\Lambda_c\cap U_Q)\cap(t+Q_m)\bigr),\\
O_m(x)&:=\sum_{Q\in\mathcal Q_S(c)}
\#\bigl(\Lambda_c\cap U_Q\cap(t+Q_m)\bigr).
\end{aligned}
\]
An inserted or removed point in $u+U_Q$, with
$u\in(C_r)_c\setminus\{0\}$ and $Q\in\mathcal Q_S(u.c)$, cannot
lie in $t+Q_m$: otherwise $t-u\in B_{Q,m}\subset W_{u.c}$,
contradicting the unique representation $x=t.c$.
Unchanged points cancel in the count difference. Thus
\[
N_{Q_m}(\kappa_{\delta,S,r}(x))-N_{Q_m}(\kappa(x))
=I_m(x)-O_m(x),
\qquad
O_m\leq (N_{Q_m}\circ\kappa)\mathbf1_{F_m}.
\]
If $x=t.c$ and moved points from $U_Q$ enter $t+Q_m$, then
$t\in A_{Q,m}$.  Their preimages lie in $s+Q_{3m}$, where
\[
s=3t+a_U-3a_D.
\]
This affine map sends $A_{Q,m}$ bijectively onto $B_{Q,m}$, with
$ds=3^d\,dt$.  For each $c$, the sets $A_{Q,m}$ lie in pairwise disjoint cubes
$Q\subset W_c$. Since $(c,t)\mapsto t.c$ is injective on
$\{(c,t):c\in C_r,\ t\in W_c\}$, the sets
$\{t.c:t\in A_{Q,m}\}$ are disjoint as $(c,Q)$ varies.
The same argument applies to $B_{Q,m}$, so
\[
\begin{split}
\int_X I_m^2\,d\mu
&\leq
\int_{C_r}\sum_{Q\in\mathcal Q_S(c)}
\int_{A_{Q,m}}
\bigl(N_{Q_{3m}}(\kappa((3t+a_U-3a_D).c))\bigr)^2
\,dt\,d\nu_r(c)\\
&=3^{-d}\int_{F_m}(N_{Q_{3m}}\circ\kappa)^2\,d\mu.
\end{split}
\]
It follows that
\[
\begin{split}
&\|N_{Q_m}\circ\kappa_{\delta,S,r}
                  -N_{Q_m}\circ\kappa\|_2^2\\
&\qquad\leq
2(1+3^{-d})\int_{F_m}(N_{Q_{3m}}\circ\kappa)^2\,d\mu.
\end{split}
\]
The function \((N_{Q_{3m}}\circ\kappa)^2\) is integrable.
By the bound on $\mu(F_m)$ and absolute continuity of the integral,
\[
\lim_{\delta\downarrow0}\,
\limsup_{S\to\infty}\,
\sup_{r\geq S^2}
\left\|
N_{Q_m}\circ\kappa_{\delta,S,r}
-
N_{Q_m}\circ\kappa
\right\|_{L^2(\mu)}
=0.
\]
If \(\varphi\in C_c(\mathbb R^d)\) and
\(\operatorname{supp}\varphi\subset Q_m\), then
\(S_\varphi(\kappa_{\delta,S,r})=S_\varphi(\kappa)\) on
\(X\setminus F_m\).  The bound $\mu(F_m)\leq C_d\delta+C_{d,m}/S$ and the subsequence
argument in the proof of \cref{lem:lowvar-rooted} give
\(\kappa_{\delta,S,r}\to\kappa\) in measure for the vague topology.
Together with the local-count \(L^2\)-convergence, the definition of
\(d_{\mathrm{root},2}\) gives
\[
\lim_{\delta\downarrow0}\,
\limsup_{S\to\infty}\,
\sup_{r\geq S^2}
 d_{\mathrm{root},2}(\kappa_{\delta,S,r},\kappa)=0.
\]
Let \(\varepsilon_{\mathcal U}>0\) satisfy
\(B_{d_{\mathrm{root},2}}(\kappa,\varepsilon_{\mathcal U})\subset\mathcal U\), and take \(0<\delta<6^{-d}\) so small that, for all sufficiently large $S$ and
every $r\geq S^2$,
\[
d_{\mathrm{root},2}(\kappa_{\delta,S,r},\kappa)
<\varepsilon_{\mathcal U}.
\]
\smallskip\noindent\emph{Variance from the empty region.}
Since \(\alpha<2\), choose an integer \(R>N\) so large that
\[
\frac{\delta}{4}\,
\rho^2\Vol_d(B_R)^2
>
j\,\Vol_d(B_R)^\alpha.
\]
Choose $S$ sufficiently large that
$d_{\mathrm{root},2}(\kappa_{\delta,S,r},\kappa)<\varepsilon_{\mathcal U}$
for every $r\ge S^2$, and
\[
h=\delta^{1/d}S>2R,
\qquad
\left(1-\frac{2R}{h}\right)^d\geq\frac12.
\]
Then choose \(r\geq S^2\) so that
\[
1-C_dS/r\geq\frac12.
\]
With
\(\kappa':=\kappa_{\delta,S,r}\), we have \(\kappa'\in\mathcal U\).
Let \(E_R^0\subset X\) consist of the points \(x=t.c\) for which
\(t\) lies in some \(H_Q\) and
\(\operatorname{dist}(t,\partial H_Q)>R\).  Since \(T_Q(U_Q)=D_Q\subset U_Q\setminus H_Q\),
\[
N_{B_R}(\kappa_{\delta,S,r}(x))=0
\qquad(x\in E_R^0).
\]
The lower bound for $\int_{C_r}\#\mathcal Q_S(c)\,d\nu_r(c)$ gives
\[
\begin{split}
\mu(E_R^0)
&\geq
\frac{(h-2R)^d}{S^d}
\left(1-C_d\frac Sr\right)\\
&=
\delta
\left(1-\frac{2R}{h}\right)^d
\left(1-C_d\frac Sr\right).
\end{split}
\]
By the choices of \(S\) and \(r\), the last expression is at least
\(\delta/4\).

Since \(\kappa'\) has intensity \(\rho\),
\[
\begin{split}
\mathcal V_{\kappa'}(R)
&=
\int_X
\left|
N_{B_R}(\kappa'(x))
-\rho\Vol_d(B_R)
\right|^2\,d\mu(x)\\
&\geq
\mu(E_R^0)\,
\rho^2\Vol_d(B_R)^2\\
&\geq
\frac{\delta}{4}\,
\rho^2\Vol_d(B_R)^2
>
j\,\Vol_d(B_R)^\alpha.
\end{split}
\]
\end{proof}

\subsection{Density of generating realizations}\label{subsec:sparse-coding}

We prove that every realization can be approximated by one from
which the original state can be recovered. For a separated cocompact
cross-section, we retain all but three points in certain Voronoi
cells and place three points near the corresponding return vector.
Two of these points have a fixed nonzero difference, and the third
records the cross-section element through a Borel injection into a
Cantor set. \Cref{lem:sparse-markers,lem:sparse-recovery} give,
respectively, the construction and the recovery map.

\begin{proposition}[Density of generating realizations]\label{prop:dense-generation}
Let \(a\in\Ad_d(X,\mu)\), let
\(\kappa\in\smash{\RR^a_{\rho,2}(X,\mu)}\), and let
\(\mathcal U\) be an open neighbourhood of \(\kappa\).  Then
\[
\mathcal U\cap\RR^{a,\mathrm{iso}}_{\rho,2}(X,\mu)\neq\varnothing.
\]
Consequently,
\(\smash{\RR^{a,\mathrm{iso}}_{\rho,2}(X,\mu)}\) is dense in
\(\smash{\RR^a_{\rho,2}(X,\mu)}\).
\end{proposition}

\subsubsection{Proof of Proposition~\ref{prop:dense-generation}}

\begin{lemma}[Sparse marker replacement]\label{lem:sparse-markers}
Let $a\in\Ad_d(X,\mu)$, let
$\kappa\in\RR^a_{\rho,2}(X,\mu)$, and fix $\varepsilon>0$.
For every sufficiently large $r$ there exist an $r$-separated
cocompact cross-section $C_r$ with transverse measure $\nu_r$,
a Borel equivariant realization
$\kappa_r\in\RR^a_{\rho,2}(X,\mu)$,
vectors $0\ne v_*\in B_\varepsilon$ and $a_0\in B_\varepsilon$,
a Borel set $E'_r\subset C_r$ with $\nu_r(E'_r)>0$,
a compact Lebesgue-null Cantor set $K\subset B_\varepsilon$,
and a Borel injection $w:E'_r\to K$ satisfying the following.
\begin{enumerate}[label=\textup{(\roman*)}]
\item For $\mu$-almost every $x$, the ordered pairs of distinct
points of $\operatorname{supp}\kappa_r(x)$ having difference $v_*$
are precisely
\[
(u+a_0,u+a_0+v_*),
\qquad u\in(C_r)_x,\quad u.x\in E'_r.
\]
For each such $u$,
\[
\operatorname{supp}\kappa_r(x)\cap(u+a_0+K)
=\{u+a_0+w(u.x)\}.
\]
\item For the removed and inserted counting measures
\[
\kappa_r^-(x):=\delta_{\operatorname{supp}\kappa(x)
  \setminus\operatorname{supp}\kappa_r(x)},\qquad
\kappa_r^+(x):=\delta_{\operatorname{supp}\kappa_r(x)
  \setminus\operatorname{supp}\kappa(x)}.
\]
For every fixed $m$ and all sufficiently large $r$,
\[
\int_XN_{Q_m}(\kappa_r^\pm(x))\,d\mu(x)
\leq3m^d\nu_r(C_r),\qquad
N_{Q_m}(\kappa_r^+(x))\leq3,
\qquad
\nu_r(C_r)\leq\Vol_d(B_{r/2})^{-1}.
\]
\end{enumerate}
\end{lemma}

\begin{lemma}[Recovery from sparse markers]\label{lem:sparse-recovery}
The realization $\kappa_r$ in \cref{lem:sparse-markers} is generating.
\end{lemma}

\begin{proof}
Fix $\varepsilon>0$ and take $\kappa_r$ from
\cref{lem:sparse-markers,lem:sparse-recovery}.  For fixed $m$, set
$N_{r,m}^\pm:=N_{Q_m}\circ\kappa_r^\pm$.
Since $\nu_r(C_r)\leq\Vol_d(B_{r/2})^{-1}\to0$,
\[
\int_XN_{r,m}^\pm\,d\mu\longrightarrow0
\qquad(r\to\infty).
\]
The inequalities
\[
(N_{r,m}^-)^2\leq(N_{Q_m}\circ\kappa)^2
 \mathbf1_{\{N_{r,m}^->0\}},\qquad
(N_{r,m}^+)^2\leq3N_{r,m}^+
\]
and $\mu\{N_{r,m}^->0\}\leq\int_XN_{r,m}^-\,d\mu$
imply $\|N_{r,m}^-\|_2+\|N_{r,m}^+\|_2\to0$.
For $\varphi\in C_c(\mathbb R^d)$ with support in $Q_m$,
\[
|S_\varphi(\kappa_r)-S_\varphi(\kappa)|
\leq\|\varphi\|_\infty(N_{r,m}^-+N_{r,m}^+).
\]
Moreover,
\[
\mu\{\kappa_r|_{Q_m}\ne\kappa|_{Q_m}\}
\leq\int_X(N_{r,m}^-+N_{r,m}^+)\,d\mu\longrightarrow0.
\]
Hence $\kappa_r\to\kappa$ in measure, and
$d_{\mathrm{root},2}(\kappa_r,\kappa)\to0$. Choose $r$
large enough that $\kappa_r\in\mathcal U$.
\end{proof}

\subsubsection{Proof of Lemma~\ref{lem:sparse-markers}}

We first choose a vector $v_*$ that is not the difference of any
two original configuration points, almost surely. The inserted
pairs with difference $v_*$ can then be distinguished from the
retained points.

By \cref{lem:concrete-version}, we may restrict $a$ and $\kappa$
to an invariant conull Borel set $X_0$ on which $a$ is free and
$\kappa$ is Borel and pointwise equivariant.  Write
$\Lambda_x:=\operatorname{supp}\kappa(x)$ for $x\in X_0$.
For a bounded Borel set $D\subset\mathbb R^d$, define
\[
\alpha(D):=
\int_X\sum_{p\in\Lambda_x\cap Q_1}
 \sum_{\substack{q\in\Lambda_x\\q\ne p}}
 \mathbf1_D(q-p)\,d\mu(x).
\]
Its integrand is bounded by
$N_{Q_1}(\kappa(x))N_{Q_1+D}(\kappa(x))$, so $\alpha$ is
locally finite by the local second-moment assumption.  Choose
$v_*\in B_\varepsilon\setminus\{0\}$ with $\alpha(\{v_*\})=0$.
Almost surely, $\Lambda_x$ has no ordered pair $(p,q)$ with
$q-p=v_*$ and $p\in Q_1$.  Translation invariance and the
integer translates of $Q_1$ show that, after restricting $X_0$
further, no two points of $\Lambda_x$ differ by $v_*$ for any
$x\in X_0$.

\smallskip\noindent\emph{Cross-section and three-point sets.}
Fix $r>8\varepsilon$. Choose an $r$-separated cocompact
Borel cross-section $C_r\subset X_0$ and its Voronoi cells $W_c$
by \eqref{eq:labelled-voronoi}. The return relation
\[
\mathscr R_r:=\{(x,u)\in X_0\times\mathbb R^d:u.x\in C_r\}
\]
is Borel because the action and $C_r$ are Borel.  Each
section $(C_r)_x$ is $r$-separated and therefore countable.
Lusin--Novikov \cite[Theorem~18.10]{Kec95} gives Borel sets
$D_j\subset X_0$ and Borel maps
$u_j:D_j\to\mathbb R^d$ whose disjoint graphs cover
$\mathscr R_r$.  In particular,
$(C_r)_x=\{u_j(x):x\in D_j\}$.
Let $\nu_r$ be the transverse measure on $C_r$.  The
cross-section and mass-transport formulas
\cite[Theorem~1.14]{BHK25} give
\[
1=\int_{C_r}\Vol_d(W_c)\,d\nu_r(c),\qquad
\rho=\int_{C_r}\#(\Lambda_c\cap W_c)\,d\nu_r(c).
\]
Consequently,
$\nu_r(C_r)\leq\Vol_d(B_{r/2})^{-1}$.
The relation
\[
\mathscr P:=\{(x,p)\in X_0\times\mathbb R^d:
 p\in\Lambda_x\}
\]
is Borel: since $\Lambda_x$ is closed,
$p\in\Lambda_x$ if and only if
$N_{B_{1/n}}(\theta_p\kappa(x))>0$ for every $n\geq1$.
The latter condition is Borel in $(x,p)$ by the Borel map
$(x,p)\mapsto\theta_p\kappa(x)$.
Its sections are countable by local finiteness.  Applying
Lusin--Novikov to $\mathscr P$, choose Borel sets
$H_j\subset X_0$ and Borel maps $p_j:H_j\to\mathbb R^d$
whose disjoint graphs cover $\mathscr P$.
Thus $\Lambda_x=\{p_j(x):x\in H_j\}$.

We need a positive transverse-measure set of cells containing at
least three original points. Take $r$ sufficiently large that
$2\nu_r(C_r)<\rho$, and set
\[
E_r:=\{c\in C_r:\#(\Lambda_c\cap W_c)\geq3\}.
\]
The displayed description of $\Lambda_c$ and the Borel
relation defining $W_c$ show that $E_r$ is Borel.
If $\nu_r(E_r)=0$, the second integral above is at most
$2\nu_r(C_r)<\rho$, a contradiction.  Thus $\nu_r(E_r)>0$.
For each $c\in E_r$, let $F(c)$ consist of the three
points $p_j(c)\in W_c$ with the smallest indices $j$.
The maps $p_j$ make $c\mapsto F(c)$ Borel.
Choose a compact Lebesgue-null Cantor set
\[
K\subset B_\varepsilon\setminus\{-v_*,0,v_*,2v_*\}
\]
and a Borel injection $w:E_r\to K$, using
\cite[Theorem~15.6]{Kec95} when $E_r$ is uncountable.

\smallskip\noindent\emph{Choice of the inserted points.}
For $c\in E_r$ and $b\in B_\varepsilon$, put
$A_b(c):=\{b,b+v_*,b+w(c)\}$. We choose one vector $b$
for all cells, excluding the values for which an inserted point
would coincide with an original point, create an unintended pair
of difference $v_*$, or place an original point in $b+K$.  Let $\mathcal B_c$
be the set of $b\in B_\varepsilon$ satisfying at least one of
\[
\begin{aligned}
&A_b(c)\cap\Lambda_c\ne\varnothing,\\
&\exists\,p\in A_b(c),\ q\in\Lambda_c\setminus F(c):
  \quad q-p\in\{-v_*,v_*\},\\
&(b+K)\cap\Lambda_c\ne\varnothing.
\end{aligned}
\]
The first two conditions hold for at most countably many $b$.
The set of $b$ satisfying the third condition is contained in
$\bigcup_{q\in\Lambda_c}(q-K)$ and is Lebesgue null.
The set $\{(c,b):c\in E_r,\ b\in\mathcal B_c\}$
is Borel.  Indeed, in each condition defining $\mathcal B_c$
we can replace a quantifier over $\Lambda_c$ by the
countable disjunction over $p_j(c)$ for which $c\in H_j$;
the condition $p_j(c)\in F(c)$ is Borel by its construction, and the
third condition becomes $p_j(c)-b\in K$ for some such $j$.
Fubini's theorem therefore gives
\[
\int_{B_\varepsilon}
 \nu_r\{c\in E_r:b\in\mathcal B_c\}\,db=0.
\]
Fix $a_0\in B_\varepsilon$ such that
$\nu_r\{c\in E_r:a_0\in\mathcal B_c\}=0$,
and define
$E'_r:=\{c\in E_r:a_0\notin\mathcal B_c\}$.
Then $\nu_r(E'_r)=\nu_r(E_r)>0$.

\smallskip\noindent\emph{Definition and measurability of the modified realization.}
We now replace $F(c)$ by $A_{a_0}(c)$ in each cell with
$c\in E'_r$, leaving all other cells unchanged. For $c\in C_r$,
define the finite set
\[
\Lambda_c^W:=
\begin{cases}
((\Lambda_c\cap W_c)\setminus F(c))\cup A_{a_0}(c),
 &c\in E'_r,\\
\Lambda_c\cap W_c,&c\notin E'_r.
\end{cases}
\]
For $x\in X_0$, set
\[
\operatorname{supp}\kappa_r(x):=
 \bigcup_{u\in(C_r)_x}(u+\Lambda_{u.x}^W).
\]
Since $a_0,v_*,w(c)\in B_\varepsilon$, we have
$A_{a_0}(c)\subset B_{2\varepsilon}\subset W_c$.
The three points of $A_{a_0}(c)$ are distinct because
$w(c)\notin\{0,v_*\}$.  For $c\in E'_r$, the definition of $\mathcal B_c$ gives
$A_{a_0}(c)\cap\Lambda_c=\varnothing$.
Thus $\Lambda_c^W$ contains exactly the points of
$\Lambda_c\cap W_c$ outside $F(c)$, together with the three
points of $A_{a_0}(c)$.  The partition by the sets
$u+W_{u.x}$ makes the union defining $\kappa_r(x)$
disjoint.  It is locally finite: the retained points form a
subset of $\Lambda_x$, while any point of
$u+A_{a_0}(u.x)$ is within $2\varepsilon$ of an
$r$-separated return vector $u$.  The maps $p_j$, the
Borel relation defining $W_c$, and the three points of
$A_{a_0}(c)$ show that $c\mapsto\Lambda_c^W$ is a
Borel finite-set map.  For every bounded Borel set $D$,
\[
N_D(\kappa_r(x))
=\sum_{j:x\in D_j}
  \#\bigl((u_j(x)+\Lambda^W_{u_j(x).x})\cap D\bigr).
\]
Each summand is Borel, so $\kappa_r$ is Borel.
The identities
$(C_r)_{t.x}=(C_r)_x-t$ and $(u-t).(t.x)=u.x$
prove equivariance.

\smallskip\noindent\emph{Point counts and second moments.}
Each replacement removes and inserts exactly three points.
Consequently, for every $c\in C_r$ we have
$\#\Lambda_c^W=\#(\Lambda_c\cap W_c)$.  The cross-section
mass-transport formula gives
\[
\int_XN_{Q_m}(\kappa_r(x))\,d\mu(x)
=m^d\int_{C_r}\#\Lambda_c^W\,d\nu_r(c)
=\rho m^d.
\]
Translation invariance of the intensity measure then gives
intensity $\rho$ on every bounded Borel set.
The points of $\Lambda_x$ removed by the construction are
precisely those in
\[
\bigcup_{\substack{u\in(C_r)_x\\u.x\in E'_r}}
 (u+F(u.x)),
\]
and the points of $\operatorname{supp}\kappa_r(x)$ not in
$\Lambda_x$ are precisely those in
\[
\bigcup_{\substack{u\in(C_r)_x\\u.x\in E'_r}}
 (u+A_{a_0}(u.x)).
\]
Both unions are disjoint.  Applying mass transport to their
counting measures yields, for every $m\geq1$,
\[
\int_XN_{Q_m}(\kappa_r^-(x))\,d\mu(x)
=\int_XN_{Q_m}(\kappa_r^+(x))\,d\mu(x)
=3m^d\nu_r(E'_r).
\]
If $p\in(u+A_{a_0}(u.x))\cap Q_m$ and
$q\in(v+A_{a_0}(v.x))\cap Q_m$, then
$\|u-p\|,\|v-q\|<2\varepsilon$ and
$\|u-v\|\leq\sqrt d\,m+4\varepsilon$.
For $r>\sqrt d\,m+4\varepsilon$, the $r$-separation
forces $u=v$.  Thus $N_{Q_m}(\kappa_r^+(x))\leq3$.
For each fixed $r,m$, the number of $r$-separated vectors in
$Q_m+B_{2\varepsilon}$ is bounded by a finite constant
$C_{r,m}$, independently of $x$.  Therefore
\[
N_{Q_m}(\kappa_r(x))
\leq N_{Q_m}(\kappa(x))+3C_{r,m},
\]
which proves that $\kappa_r$ has local second moments.

\smallskip\noindent\emph{Identification of the pairs with difference $v_*$.}
We verify that the pairs with difference $v_*$ are exactly
those asserted in part (i). Fix $x\in X_0$.
No two points of $\Lambda_x$ differ by $v_*$, so the
points retained from $\Lambda_x$ contain no such pair.
For each $u\in(C_r)_x$ with $c:=u.x\in E'_r$,
the only ordered pair of points of $u+A_{a_0}(c)$
with difference $v_*$ is
$(u+a_0,u+a_0+v_*)$:
all other possibilities would require
$w(c)\in\{-v_*,0,v_*,2v_*\}$.
The second condition defining $\mathcal B_c$ excludes a pair
of difference $v_*$ between $u+A_{a_0}(c)$ and
$u+((\Lambda_c\cap W_c)\setminus F(c))$.
If $q\in\Lambda_x\setminus(u+W_c)$, then
$\|q-u\|\geq r/2$ because $B_{r/2}^{\circ}\subset W_c$.
For $p\in u+A_{a_0}(c)$ we have $\|p-u\|<2\varepsilon$,
so $\|p-q\|\geq r/2-2\varepsilon>\|v_*\|$.
If $u,v\in(C_r)_x$ are distinct and
$p\in u+A_{a_0}(u.x)$, $q\in v+A_{a_0}(v.x)$, then
$\|p-q\|\geq r-4\varepsilon>\|v_*\|$.
Thus the ordered pairs of difference $v_*$ in
$\operatorname{supp}\kappa_r(x)$ are exactly those in (i).

For $u\in(C_r)_x$ with $c:=u.x\in E'_r$, we have
$(a_0+K)\cap\Lambda_c=\varnothing$ by the definition of
$E'_r$.
The set $a_0+K$ lies in $B_{2\varepsilon}\subset W_c$.
Since $0,v_*\notin K$, its intersection with
$\Lambda_c^W$ is $\{a_0+w(c)\}$.  The partition
by the sets $u+W_{u.x}$ therefore gives
\[
\operatorname{supp}\kappa_r(x)\cap(u+a_0+K)
=\{u+a_0+w(c)\}.
\]
This proves (i) and (ii). \hfill$\square$

\subsubsection{Proof of Lemma~\ref{lem:sparse-recovery}}

First we show that almost every configuration contains an inserted
triple. We then construct a Borel map recovering $x$ from any such
triple in $\kappa_r(x)$.

Let
\[
S_r:=\{x\in X_0:\exists\,u\in(C_r)_x
  \text{ with }u.x\in E'_r\}.
\]
Retain the partial Borel maps $u_j:D_j\to\mathbb R^d$
from the proof of \cref{lem:sparse-markers}.  They give
\[
S_r=\bigcup_{j\geq1}\{x\in D_j:u_j(x).x\in E'_r\}.
\]
Thus $S_r$ is Borel.  The return-time identity makes it invariant.  Every point
$t.c$ with $c\in E'_r$ and $t\in W_c$ belongs to $S_r$.
The cross-section formula therefore gives
\[
\mu(S_r)\geq
\int_{E'_r}\Vol_d(W_c)\,d\nu_r(c)
\geq\Vol_d(B_{r/2})\,\nu_r(E'_r)>0.
\]
Since $a$ is ergodic, $\mu(S_r)=1$.

\smallskip\noindent\emph{Borel reconstruction.}
Fix $x\in S_r$ and write
$\Lambda_x^r:=\operatorname{supp}\kappa_r(x)$.
By \cref{lem:sparse-markers}(i), there is an ordered pair
$(z,z+v_*)\in (\Lambda_x^r)^2$ and exactly one point
$p\in\Lambda_x^r\cap(z+K)$.  For any such pair, the same lemma
gives $c\in E'_r$ and $u\in(C_r)_x$ with
\[
z=u+a_0,\qquad p=z+w(c),\qquad c=u.x.
\]
The injection $w$ determines $c$ uniquely from $p-z$.
It follows that
\[
x=(a_0-z).\,w^{-1}(p-z).
\]
Lusin--Souslin \cite[Theorem~15.1]{Kec95} implies that
$w(E'_r)$ is Borel and $w^{-1}:w(E'_r)\to E'_r$ is Borel.
Let $\mathscr H$ be the set of triples
$(\omega,z,p)\in\operatorname{Conf}(\mathbb R^d)
\times\mathbb R^d\times\mathbb R^d$ satisfying
\[
z,z+v_*,p\in\operatorname{supp}\omega,\qquad
p-z\in w(E'_r),\qquad N_{z+K}(\omega)=1.
\]
The relation $\mathscr H$ is Borel: point incidence in a
locally finite configuration is Borel, and
$N_{z+K}(\omega)=N_K(\theta_z\omega)$ is Borel by joint
continuity of translation and Borel measurability of $N_K$.
For each $\omega$, the section $\mathscr H_\omega$ is
countable because both $z$ and $p$ belong to the countable
set $\operatorname{supp}\omega$.  Lusin--Novikov
\cite[Theorem~18.10]{Kec95} therefore shows that the
projection $\Omega$ of $\mathscr H$ onto the configuration
space is Borel and provides a Borel map
$\omega\mapsto(z(\omega),p(\omega))$ on $\Omega$ with
$(\omega,z(\omega),p(\omega))\in\mathscr H$.
Consequently,
\[
\omega\longmapsto
(a_0-z(\omega)).\,w^{-1}(p(\omega)-z(\omega))
\]
is a Borel map from $\Omega$ to $X_0$.  For every
$x\in S_r$, \cref{lem:sparse-markers}(i) implies that
$\kappa_r(x)\in\Omega$ and every pair
$(z,p)\in\mathscr H_{\kappa_r(x)}$ gives the value $x$.
Thus this map sends $\kappa_r(x)$ to $x$ on the conull
set $S_r$, and $\kappa_r$ is generating. \hfill$\square$

The density statements in
\cref{prop:low-variance-density,prop:high-variance-density,prop:dense-generation},
together with \cref{prop:ihyp-baire-reduction}, give
\[
\{\kappa\in\RR^{a,\mathrm{iso}}_{\rho,2}(X,\mu):
  \kappa_*\mu\in\IHyp_\rho\}
\]
is residual in \(\smash{\RR^{a,\mathrm{iso}}_{\rho,2}(X,\mu)}\) for every
$a\in\Ad_d(X,\mu)$, and \(\IHyp_\rho\) is residual in
\(\smash{\EPP_{\rho,2}}\).  Taking $\alpha=1$ in the definition of
intermittent hyperuniformity gives
\[
\limsup_{R\to\infty}
\frac{V_\eta(R)}{\Vol_d(B_R)}=+\infty
\qquad(\eta\in\IHyp_\rho),
\]
so \(\Hyp_\rho\cap\smash{\EPP_{\rho,2}}\) is meagre in
\(\smash{\EPP_{\rho,2}}\).

\section{Maximal and linear rigidity}\label{sec:rigidity}

We prove that maximal rigidity is a $G_\delta$ condition and obtain
density by sparse marker coding.  For linear maximal rigidity,
polynomial moment identities allow local linear statistics to be
approximated in $L^2$ by exterior statistics.

\subsection{Maximal rigidity}\label{subsec:maximal-rigidity}

\begin{proposition}\label{prop:mrig-gdelta}
The set \(\MRig_\rho\) is \(G_\delta\) in \(\PP_{\rho,2}\).  Therefore, for every
\(a\in\Ad_d(X,\mu)\),
\[
\{\kappa\in\RR^a_{\rho,2}(X,\mu):\kappa_*\mu\in\MRig_\rho\}
\]
is \(G_\delta\) in \(\smash{\RR^a_{\rho,2}(X,\mu)}\).
\end{proposition}

\begin{proof}
Fix an integer $m\geq1$.  Choose a countable family of
nonnegative functions $f_{m,\ell}\in C_c(\operatorname{int}B_m)$ such that
every nonnegative $f\in C_c(\operatorname{int}B_m)$ admits uniform
approximants from the family whose supports lie in a common compact subset
of $\operatorname{int}B_m$, and put
\[
F_{m,\ell}(\omega):=e^{-S_{f_{m,\ell}}(\omega)}.
\]
For $f\geq0$ continuous and compactly supported, $e^{-S_f}$ is bounded
and continuous in the vague topology, and
$e^{-S_f}e^{-S_g}=e^{-S_{f+g}}$.  Since Laplace functionals separate
locally finite counting measures, the functions $F_{m,\ell}$ generate
the $\sigma$-algebra of $\omega|_{\operatorname{int}B_m}$.
Since $N_{\partial B_m}=0$ almost surely by \cref{lem:null-window}, it also
generates the $\sigma$-algebra of $\omega|_{B_m}$ modulo null sets.

Choose an analogous countable family of nonnegative functions in
$C_c(\mathbb R^d\setminus B_m)$, and let $\mathcal A_m^{\rm out}$ be the
unital algebra of finite rational linear combinations of finite products
of their Laplace functionals.  This is a countable family of bounded
continuous functions.  The algebra generates $\mathcal F_{B_m^c}$: compactly supported
continuous functions on $B_m^c$ determine the exterior configuration.
Since this algebra contains the constants and generates
$\mathcal F_{B_m^c}$, the functional monotone-class theorem gives its
density in $L^2(\mathcal F_{B_m^c},\eta)$.
Define
\[
D_{m,\ell}(\eta)
:=
\inf_{G\in\mathcal A_m^{\rm out}}
\int |F_{m,\ell}-G|^2\,d\eta.
\]
For fixed $G\in\mathcal A_m^{\rm out}$, the function
$|F_{m,\ell}-G|^2$ is bounded and continuous on configuration space.
Consequently the map
$\eta\mapsto\int|F_{m,\ell}-G|^2\,d\eta$ is continuous,
even in the weak topology.  Since
\[
\{\eta:D_{m,\ell}(\eta)<\varepsilon\}
=\bigcup_{G\in\mathcal A_m^{\rm out}}
\left\{\eta:
\int|F_{m,\ell}-G|^2\,d\eta<\varepsilon\right\},
\]
this set is open.

Because the $L^2(\eta)$-closure of $\mathcal A_m^{\rm out}$ is
$L^2(\mathcal F_{B_m^c},\eta)$, the equality
$D_{m,\ell}(\eta)=0$ is equivalent to $F_{m,\ell}$ being
$\mathcal F_{B_m^c}$-measurable modulo $\eta$-null sets.  The inside
functionals generate the restriction $\omega|_{B_m}$, so
\[
\MRig_\rho
=
\bigcap_{m,\ell,k\geq1}
\{\eta\in\PP_{\rho,2}:D_{m,\ell}(\eta)<k^{-1}\}.
\]
If $D_{m,\ell}(\eta)=0$ for all $\ell$,
$\omega|_{B_m}$ is a Borel function of $\omega|_{B_m^c}$ almost surely, by
the Doob--Dynkin lemma.  For bounded Borel $A\subset B_m$ we have
$\mathcal F_{B_m^c}\subset\mathcal F_{A^c}$, so $\omega|_A$ is then also
determined by $\omega|_{A^c}$.  Hence $\MRig_\rho$ is $G_\delta$, and
continuity of $\kappa\mapsto\kappa_*\mu$ from
\cref{prop:global-polish} gives the rooted assertion.
\end{proof}

\begin{proposition}\label{prop:mrig-dense-rooted}
Let \(a\in\Ad_d(X,\mu)\), let
\(\kappa\in\smash{\RR^a_{\rho,2}(X,\mu)}\), and let
\(\mathcal U\) be a neighbourhood of \(\kappa\).  Then
\[
\{\lambda\in\mathcal U\cap\RR^{a,\mathrm{iso}}_{\rho,2}(X,\mu):
  \lambda_*\mu\in\MRig_\rho\}
\neq\varnothing.
\]
\end{proposition}

\begin{proof}
Let $\kappa_r$ be the realizations from the proof of
\cref{prop:dense-generation}.  Choose $r$ sufficiently large that
$\kappa_r\in\mathcal U$, using $\kappa_r\to\kappa$ in the rooted local
$L^2$ topology.  Let $E_r'\subset C_r$ be the Borel set on which the three-point
replacement of \cref{lem:sparse-markers} is performed. It has positive
transverse measure. Define
\[
T_x:=\{u\in(C_r)_x:u.x\in E_r'\}.
\]
For \(u\in T_x\), with \(c:=u.x\), the three inserted points are
\[
u+a_0,\qquad
u+a_0+v_*,
\qquad
u+a_0+w(c).
\]
By \cref{lem:sparse-markers}(i), for almost every $x$, every
ordered pair $(z,z+v_*)$ of points in
$\operatorname{supp}\kappa_r(x)$ has the form
$(u+a_0,u+a_0+v_*)$ for some $u\in T_x$.

The translation-invariant process $x\mapsto\delta_{T_x}$ has intensity
$\nu_r(E_r')>0$.  The cardinality $\#T_x$ is invariant under $x\mapsto t.x$ and hence,
by ergodicity, equals a constant $k\in\mathbb N_0\cup\{\infty\}$
for almost every $x$.
If $k<\infty$, then
\[
\nu_r(E_r')\Vol_d(Q_n)
=\int_X\#(T_x\cap Q_n)\,d\mu(x)\leq k
\qquad(n\geq1),
\]
a contradiction as $n\to\infty$.  Thus $T_x$ is infinite almost surely.
Since $T_x$ is $r$-separated, it is unbounded.

Fix $m\geq1$ and put $A:=B_m$.  Set
\[
M:=\{a_0,a_0+v_*\}\cup(a_0+K).
\]
Since $M$ is bounded and $T_x$ is unbounded, almost surely some
$u\in T_x$ satisfies $u+M\subset A^c$.  The first two points of the
corresponding triple are $(z,z+v_*)$, where $z=u+a_0$;
its third point is in $z+K$.
Let $\omega^{\rm out}:=\kappa_r(x)|_{B_m^c}$, and define
\[
\begin{split}
Z_m(\omega^{\rm out}):=\bigl\{z\in\operatorname{supp}(\omega^{\rm out}):{}&
 z+v_*\in\operatorname{supp}(\omega^{\rm out}),\\[-2pt]
 & (z+K)\cap B_m=\varnothing,\quad
 N_{z+K}(\omega^{\rm out})=1\bigr\}.
\end{split}
\]
For $x$ in the conull set on which \cref{lem:sparse-markers}(i)
holds, every $z\in Z_m(\omega^{\rm out})$ determines an ordered pair
$(z,z+v_*)$ in $\operatorname{supp}\kappa_r(x)$.  By that lemma
there is $u\in(C_r)_x$ with $c:=u.x\in E'_r$ and $z=u+a_0$.
The same lemma gives
$\operatorname{supp}\kappa_r(x)\cap(z+K)=\{z+w(c)\}$.
Since $z+K\subset B_m^c$, this is also the unique point of
$\omega^{\rm out}$ in $z+K$.  Conversely, if $u\in T_x$ satisfies $u+M\subset B_m^c$, then
$u+a_0\in Z_m(\omega^{\rm out})$.  Thus
$Z_m(\kappa_r(x)|_{B_m^c})\ne\varnothing$ for almost every $x$.

Consider the subset of the product of configuration space and
$\mathbb R^d$
\[
\mathscr Z_m:=\{(\omega,z):z\in Z_m(\omega)\}.
\]
This subset is Borel.  Indeed, the conditions
$z,z+v_*\in\operatorname{supp}\omega$ are Borel in
$(\omega,z)$ by joint Borel measurability of point incidence,
$N_{z+K}(\omega)=N_K(\theta_z\omega)$ is Borel, and
$\{z:(z+K)\cap B_m=\varnothing\}$ is open because $K$ and
$B_m$ are compact.  Each vertical section of $\mathscr Z_m$ is
countable, since it is contained in $\operatorname{supp}\omega$.
Lusin--Novikov implies that
$\Omega_m:=\{\omega:Z_m(\omega)\ne\varnothing\}$ is Borel
and supplies a Borel map $z: \Omega_m\to\mathbb R^d$
satisfying $z(\omega)\in Z_m(\omega)$.  The subset
\[
\{(\omega,p):\omega\in\Omega_m,\quad
p\in\operatorname{supp}\omega\cap(z(\omega)+K)\}
\]
of $\Omega_m\times\mathbb R^d$ is Borel and has singleton
vertical sections.  Another application of Lusin--Novikov gives its
unique point as a Borel map $p:\Omega_m\to\mathbb R^d$.
For $\omega=\kappa_r(x)|_{B_m^c}$ and almost every $x$, the preceding
identities give $p(\omega)-z(\omega)=w(c)$ and
\[
x=(a_0-z(\omega)).\,w^{-1}(p(\omega)-z(\omega)).
\]
By Lusin--Souslin \cite[Theorem~15.1]{Kec95}, $w(E'_r)$ is Borel
and $w^{-1}:w(E'_r)\to E'_r$ is Borel. The set
$\Omega_m':=\{\omega\in\Omega_m:p(\omega)-z(\omega)\in w(E'_r)\}$
is Borel and contains $\kappa_r(x)|_{B_m^c}$ for almost every $x$.
The displayed formula defines a Borel reconstruction on $\Omega_m'$.
Thus $\kappa_r(x)$ is determined by its restriction to $B_m^c$.  For any bounded
Borel $A$, choose $m$ with $A\subset B_m$; then
$\mathcal F_{B_m^c}\subset\mathcal F_{A^c}$.  Hence the restriction
to $A$ is determined by its exterior, and $(\kappa_r)_*\mu$ is
maximally rigid.
\end{proof}

\begin{corollary}\label{cor:mrig-residual}
For every \(a\in\Ad_d(X,\mu)\),
\[
\{\kappa\in\RR^{a,\mathrm{iso}}_{\rho,2}(X,\mu):
  \kappa_*\mu\in\MRig_\rho\}
\]
is residual in
\(\smash{\RR^{a,\mathrm{iso}}_{\rho,2}(X,\mu)}\).  Moreover,
\(\MRig_\rho\cap\smash{\EPP_{\rho,2}}\) is residual in
\(\smash{\EPP_{\rho,2}}\).
\end{corollary}

\begin{proof}
The rooted statement follows from
\cref{prop:mrig-gdelta,prop:mrig-dense-rooted}.  The global statement
follows from \cref{prop:local-global}.
\end{proof}

\subsection{Linear maximal rigidity}\label{subsec:linear-maximal-rigidity}

\begin{proposition}\label{prop:linear-maximal-rigidity}
The set \(\LMRig_\rho\) is \(G_\delta\) in \(\PP_{\rho,2}\).  For every
\(a\in\Ad_d(X,\mu)\),
\[
\{\kappa\in\RR^{a,\mathrm{iso}}_{\rho,2}(X,\mu):
  \kappa_*\mu\in\LMRig_\rho\}
\]
is residual in \(\smash{\RR^{a,\mathrm{iso}}_{\rho,2}(X,\mu)}\), while
\(\LMRig_\rho\cap\smash{\EPP_{\rho,2}}\) is residual in
\(\smash{\EPP_{\rho,2}}\).
\end{proposition}

\medskip

For $m\geq1$ and $\mathbf j\in\mathbb N_0^d$, set
\[
\gamma_{m,\mathbf j}(u):=u^{\mathbf j}\mathbf1_{B_m}(u).
\]
Choose a countable family
$\mathcal G_m\subset C_c(\mathbb R^d\setminus B_m)$ whose rational span
uniformly approximates every $f\in C_c(\mathbb R^d\setminus B_m)$
with approximating supports in a compact set depending on $f$. Define
\[
D_{m,\mathbf j}(\eta)
:=
\inf_{g\in\operatorname{span}_{\mathbb Q}\mathcal G_m}
\|Z_{\gamma_{m,\mathbf j}}^\eta-Z_g^\eta\|_{L^2(\eta)}.
\]
\begin{lemma}\label{lem:linear-gdelta}
The set $\LMRig_\rho$ is $G_\delta$ in $\PP_{\rho,2}$.
\end{lemma}

\begin{proof}
For fixed $g$, put $h:=\gamma_{m,\mathbf j}-g$.  Then
\[
\|Z_h^\eta\|_2^2
=
\int S_h^2\,d\eta-\left(\rho\int h\right)^2.
\]
The function $h$ is bounded and compactly supported and is continuous
outside $\partial B_m$.  By \cref{lem:null-window}, the configuration
has no point on this sphere almost surely under any measure in
$\PP_{\rho,2}$.  The local $L^2$ characterization in
\cref{prop:local-L2-characterization} therefore gives continuity of
$\smash{\int}S_h^2\,d\eta$.  Consequently
$\{\eta:D_{m,\mathbf j}(\eta)<\varepsilon\}$ is open: it is a union,
over $g$ in the countable rational span, of strict sublevel sets of
continuous functions.

For bounded Borel $W$, define the control measure on Borel sets
$C\subset\mathbb R^d$ by
\[
\nu_{W,\eta}(C)
:=
\int N_W(\omega)\,\omega(C\cap W)\,d\eta(\omega).
\]
This is a finite positive measure, because
$\nu_{W,\eta}(W)=\int N_W^2\,d\eta<\infty$.
For bounded Borel $f,h$ vanishing outside $W$, Cauchy--Schwarz in each configuration
gives
\begin{equation}\label{eq:control-measure}
\|Z_f^\eta-Z_h^\eta\|_2^2
\leq
\int_W|f-h|^2\,d\nu_{W,\eta}.
\end{equation}
Indeed, the left-hand side is at most the second moment of $S_{f-h}$, while
$|S_{f-h}|^2\leq N_WS_{|f-h|^2}$.

Let $g$ be bounded, compactly supported, and zero on $B_m$,
and choose a bounded Borel set $W$ containing its support and a neighbourhood
of that support.  By \cref{lem:null-window},
$\nu_{W,\eta}(\partial B_m)=0$.  For
$U_\varepsilon:=\{u\in W:\operatorname{dist}(u,\partial B_m)<\varepsilon\}$,
\[
\|Z_{g\mathbf1_{U_\varepsilon}}^\eta\|_2^2
\leq\|g\|_\infty^2\nu_{W,\eta}(U_\varepsilon)
\longrightarrow0\qquad(\varepsilon\downarrow0),
\]
by \eqref{eq:control-measure}.  On the remaining compact
region, bounded Borel functions can be approximated in
$L^2(\nu_{W,\eta})$ by continuous functions supported away from $B_m$.
The defining approximation property of $\mathcal G_m$ then approximates
these functions uniformly, and hence in the same $L^2$ space, by its
rational span.  Another application of \eqref{eq:control-measure} yields
\begin{equation}\label{eq:exterior-linear-span}
\overline{\operatorname{span}}
\{Z_g^\eta:g\in\mathcal G_m\}
=\mathcal H_{B_m^c}^\eta.
\end{equation}

By \eqref{eq:exterior-linear-span}, $D_{m,\mathbf j}(\eta)=0$ if and only if
$Z_{\gamma_{m,\mathbf j}}^\eta\in\mathcal H_{B_m^c}^\eta$.  Hence
\[
\LMRig_\rho
=
\bigcap_{m\geq1}
\bigcap_{\mathbf j\in\mathbb N_0^d}
\bigcap_{k\geq1}
\{\eta\in\PP_{\rho,2}:D_{m,\mathbf j}(\eta)<k^{-1}\}.
\]
Let $\eta$ belong to the right-hand side.  Then
$\smash{Z_{\gamma_{m,\mathbf j}}^\eta}\in\smash{\mathcal H_{B_m^c}^\eta}$ for every
$m,\mathbf j$, hence the same holds for $Z_{p\mathbf1_{B_m}}^\eta$ for every
polynomial $p$.  The control measure $\nu_{B_m,\eta}$ is finite.  Polynomials are
uniformly dense in $C(B_m)$ by Stone--Weierstrass, and continuous
functions are dense in $L^2(\nu_{B_m,\eta})$.  Thus every bounded
Borel function $f$ vanishing outside $B_m$ admits polynomials $p_n$ with
\[
\int_{B_m}|f-p_n|^2\,d\nu_{B_m,\eta}\longrightarrow0
\qquad(n\to\infty).
\]
By \eqref{eq:control-measure},
$Z_{p_n\mathbf1_{B_m}}^\eta\to \smash{Z_f^\eta}$ in $L^2(\eta)$, so
$\smash{Z_f^\eta}\in\smash{\mathcal H_{B_m^c}^\eta}$.  If $A$ is bounded and Borel, choose $m$ with
$A\subset B_m$; then $\smash{\mathcal H_{B_m^c}^\eta}\subset\smash{\mathcal H_{A^c}^\eta}$,
and $\eta$ is linearly maximally rigid.

Conversely, if \(\eta\in\LMRig_\rho\), then
\(\smash{Z_{\gamma_{m,\mathbf j}}^\eta}\in\smash{\mathcal H_{B_m^c}^\eta}\) for every
\(m,\mathbf j\).  The description of \(\smash{\mathcal H_{B_m^c}^\eta}\) just proved gives
\(D_{m,\mathbf j}(\eta)=0\).  This proves the reverse inclusion and the \(G_\delta\) assertion.
\end{proof}

\begin{lemma}\label{lem:linear-density}
Let $a\in\Ad_d(X,\mu)$, $m\geq1$, $\mathbf j\in\mathbb N_0^d$,
and $\varepsilon>0$. Then
\[
\{\kappa\in\RR^a_{\rho,2}(X,\mu):
  D_{m,\mathbf j}(\kappa_*\mu)<\varepsilon\}
\]
is dense in $\RR^a_{\rho,2}(X,\mu)$.
\end{lemma}

\begin{proof}
Fix $a\in\Ad_d(X,\mu)$, $m\geq1$, $\mathbf j\in\mathbb N_0^d$, and
$\varepsilon>0$.  Let $\kappa\in\smash{\RR^a_{\rho,2}(X,\mu)}$ and let
$\mathcal U$ be a neighbourhood of $\kappa$.  Choose
$\varepsilon_{\mathcal U}>0$ with
$B_{d_{\mathrm{root},2}}(\kappa,\varepsilon_{\mathcal U})\subset\mathcal U$.
Let $s:=|\mathbf j|$ and fix an integer $p>d+s$.
Since $2p>d+1$, the cell-replacement construction of
\cref{subsec:low-variance-surgery} applies with this value of $p$.
In every modified $L$-cell, the identity
\eqref{eq:modified-core-cell-moments} holds for all polynomials of degree
less than $p$.  The construction preserves intensity and satisfies
$N_{Q_m}(\kappa')\leq N_{Q_m}(\kappa)+D_m$ for each fixed $m$;
see \cref{lem:lowvar-inserted}.

Let \(\varphi\in C_c^\infty(\mathbb R^d)\) satisfy
\(\varphi=1\) on \(B_1\) and \(\operatorname{supp}\varphi\subset B_2\), and set
\[
f_R(u):=u^{\mathbf j}\varphi(u/R).
\]
Fix $\delta_{\mathcal U}:=\delta_*(\kappa,\varepsilon_{\mathcal U})$
from \cref{lem:lowvar-rooted}.
Choose the repair cube $D$ and the tolerance $\gamma>0$ so that
$\vartheta:=\Vol_d(D)$ satisfies
$\vartheta+\gamma<\delta_{\mathcal U}/4$.
The choice of $D$ fixes the moment-completion threshold $\varepsilon_0$ in
\cref{prop:relative-moment-repair}.  Choose $q$ large enough to satisfy
that threshold and the estimates of
\cref{lem:weighted-cell-ergodic} with tolerance $\gamma$,
and put $L=(q/\rho)^{1/d}$.  Each modified $L$-cell contains exactly $q=\rho L^d$ points by
\eqref{eq:modified-core-cell-moments}.  Its point measure therefore differs
from $\rho\Vol_d|_Q$ by a signed measure of total variation at most
$2\rho L^d$.

With $L$ fixed, choose $R>\max\{m,L\}$ so large that
\[
C L^p R^{d+s-p}<\frac{\varepsilon}{4},
\]
where $C$ depends on $d,\rho,p,\mathbf j,\varphi$ but is independent of
$R$ and the later choices of $\zeta,S,r,z$. Such an $R$
exists because $p>d+s$.  The local bound on inserted points in the
construction gives an integrable majorant
\[
Y_R:=C_R\bigl(1+N_{B_{2R+1}}(\kappa)^2\bigr).
\]
The constant $C_R$ does not depend on $\zeta,S,r,z$.  Since $\kappa$ has local second
moments, $Y_R\in L^1(\mu)$.  Absolute continuity of its integral
provides $\delta_R>0$ such that
\[
\mu(E)<\delta_R
\quad\Longrightarrow\quad
\int_EY_R\,d\mu<\frac{\varepsilon^2}{4}.
\]
Choose $0<\zeta<1/4$ so that
\[
\zeta<\delta_{\mathcal U}/4,
\qquad C_d\zeta<\delta_R/3.
\]
By \cref{lem:weighted-cell-ergodic}, $b_{S,\zeta}\to0$ as
$S\to\infty$ for this fixed $\zeta$, independently of $r$.
Choose integers $\ell>\ell'$ with $\ell-\ell'$ even, and put
$S:=\ell L$ and $S':=\ell'L$.  Choose $\ell,\ell'$ sufficiently large
that the shell lattice satisfies the capacity bound
$\#\mathcal G_Q\geq5\rho V_{\mathrm{out}}/4$ for every $S$-cube,
and that
\[
\zeta\leq1-(S'/S)^d\leq2\zeta,
\qquad S>4R+2\sqrt d\,L,
\]
and
\[
4b_{S,\zeta}+C_d\frac{2R+\sqrt d\,L}{S}<\delta_R/3.
\]
Finally choose $r\geq S^2$ so that
\[
S/r<\delta_{\mathcal U}/4,
\qquad C_d S/r<\delta_R/3.
\]
Apply the replacement rule in
\cref{subsec:low-variance-surgery} with the chosen parameters,
obtaining a realization $\kappa'$ and a shift $z\in[0,S)^d$.
Set
\[
E_{\rm exc}:=X\setminus
\left\{t.c:
\begin{array}{l}
c\in C_r,\ Q\in\mathcal Q_z(c),\ Q\text{ count-compatible},\\
t\in Q^{\mathrm{in}},\quad
\operatorname{dist}(t,\partial Q^{\mathrm{in}})>2R+\sqrt d\,L
\end{array}
\right\}.
\]
The Voronoi covering estimate \eqref{eq:voronoi-grid-cover} and
the bound $e_1(z)\leq4b_{S,\zeta}$ in \cref{lem:lowvar-shift} give
\[
\mu(E_{\rm exc})\leq
4b_{S,\zeta}+C_d\zeta
+C_d\frac{2R+\sqrt d\,L}{S}+C_d\frac Sr<\delta_R.
\]
Also $\vartheta+\gamma+\zeta+S/r<\delta_{\mathcal U}$, so
$d_{\mathrm{root},2}(\kappa',\kappa)<\varepsilon_{\mathcal U}$.
Thus $\kappa'\in\mathcal U$, and the bound by $Y_R$ gives
\[
\int_{E_{\rm exc}}
\bigl|Z_{f_R}^{(\kappa')_*\mu}(\kappa'(x))\bigr|^2\,d\mu(x)
<\frac{\varepsilon^2}{4}.
\]

For $x\notin E_{\rm exc}$, the support of $f_R$ lies in the union
of modified $L$-cells in a single core: the extra margin
$\sqrt d\,L$ ensures that every cell intersecting $B_{2R}$ belongs to that
core.  For each of these cells $Q$, define
\[
\Delta_Q:=\kappa'(x)|_Q-\rho\Vol_d|_Q.
\]
By \eqref{eq:modified-core-cell-moments},
\[
\int P\,d\Delta_Q=0\qquad(\deg P<p).
\]
Also $\|\Delta_Q\|_{\mathrm{TV}}\leq2\rho L^d$, since the cell
contains exactly $\rho L^d$ points.
Let $c_Q$ be the centre of $Q$.  The Taylor polynomial
$T^{p-1}_{c_Q}f_R$ has degree less than $p$ and integrates to zero
against $\Delta_Q$.  Since every derivative of $f_R$ of order $p$ is
bounded by $C R^{s-p}$, Taylor's theorem and the preceding total
variation bound give
\[
\left|\int f_R\,d\Delta_Q\right|
\leq
\|\Delta_Q\|_{\mathrm{TV}}
\sup_Q|f_R-T^{p-1}_{c_Q}f_R|
\leq
C L^{d+p}R^{s-p}.
\]
Since $R>L$, at most $C_d(R/L)^d$ cells intersect $B_{2R}$.
Each contributes at most $CL^{d+p}R^{s-p}$; the other cells contribute
zero because $f_R$ vanishes there.  Therefore
\[
\left|
\int f_R\,d(\kappa'(x)-\rho\Vol_d)
\right|
\leq
C L^p R^{d+s-p}
<
\frac{\varepsilon}{4}.
\]
Thus
\[
\int_{X\setminus E_{\rm exc}}
\bigl|Z_{f_R}^{(\kappa')_*\mu}(\kappa'(x))\bigr|^2\,d\mu(x)
\leq\frac{\varepsilon^2}{16}.
\]
Combining the two parts gives
\[
\|Z_{f_R}^{(\kappa')_*\mu}\|_{L^2((\kappa')_*\mu)}
<\varepsilon.
\]

Since $R>m$ and $\varphi=1$ on $B_1$, we have
$f_R(u)=u^{\mathbf j}$ on $B_m$.  The function
\[
g_R:=f_R-\gamma_{m,\mathbf j}
\]
is bounded, compactly supported, and vanishes on \(B_m\).  Hence $Z_{-g_R}^{(\kappa')_*\mu}\in
\mathcal H_{B_m^c}^{(\kappa')_*\mu}$.
Since $D_{m,\mathbf j}$ is the distance to this closed subspace,
\[
D_{m,\mathbf j}((\kappa')_*\mu)
\leq
\|Z_{\gamma_{m,\mathbf j}}^{(\kappa')_*\mu}
      -Z_{-g_R}^{(\kappa')_*\mu}\|_2
=
\|Z_{f_R}^{(\kappa')_*\mu}\|_2
<
\varepsilon.
\]
Since $\mathcal U$ and $\varepsilon$ were arbitrary, for every
$k\geq1$ the preimage of
$\{\eta:D_{m,\mathbf j}(\eta)<k^{-1}\}$ in
$\RR^a_{\rho,2}(X,\mu)$ is dense.  It is open by
the $G_\delta$ argument above and continuity of the push-forward map.
\end{proof}

\begin{proof}[Proof of \cref{prop:linear-maximal-rigidity}]
By \cref{lem:linear-gdelta}, the preimage of $\LMRig_\rho$
is a $G_\delta$ in each rooted realization space. By
\cref{lem:linear-density}, its defining open sets are dense.
Intersecting with the dense $G_\delta$ of generating realizations
from \cref{prop:dense-generation} proves the fixed-action assertion.
The global conclusion follows from \cref{prop:local-global}.
\end{proof}

\emph{Proof of \cref{thm:main-rigidity,cor:main-intersection,prop:main-fixed-action}.}
For each $a\in\Ad_d(X,\mu)$, the intermittent-hyperuniform and linearly
maximally rigid realizations are residual in
$\smash{\RR^{a,\mathrm{iso}}_{\rho,2}(X,\mu)}$ by the end of
\cref{sec:intermittent-hyperuniformity} and
\cref{prop:linear-maximal-rigidity}; their intersection gives
\cref{prop:main-fixed-action}.  The classes \(\IHyp_\rho\) and
\(\LMRig_\rho\), intersected with \(\EPP_{\rho,2}\), are residual in
$\EPP_{\rho,2}$, which proves \cref{thm:main-rigidity} and the first
assertion of \cref{cor:main-intersection}.  Finally,
$\Hyp_\rho\cap\EPP_{\rho,2}$ is meagre by the end of
\cref{sec:intermittent-hyperuniformity}; intersecting its residual
complement with $\MRig_\rho\cap\EPP_{\rho,2}$ gives the last assertion.

\section{Periodic approximation and pairwise negative correlation}
\label{sec:pnc}

The approximation of measure-preserving transformations by periodic
ones, and its use in spaces of actions, are developed in the
Rokhlin--Halmos tradition; see \cite[Chaps.~1--2]{Kec10}.
We adapt the block-approximation method to point-process measures,
controlling both intensity and local second moments.  Periodic
approximation proves density of hyperuniform processes and, together
with the covariance representation below, excludes pairwise negative
correlation on a dense set.

\subsection{Periodic approximation}\label{subsec:periodic-approximation}

For a full-rank lattice $\mathbf L<\mathbb R^d$, write
\[
\mathbf L^*:=\{\xi\in\mathbb R^d:
\langle\xi,t\rangle\in\mathbb Z\text{ for every }t\in\mathbf L\},
\qquad
\lambda_1(\mathbf L):=\min_{t\in\mathbf L\setminus\{0\}}\|t\|.
\]
A nonempty locally finite set \(P\subset\mathbb R^d\) is \emph{periodic} if its
stabilizer
\[
\operatorname{Stab}(P):=\{t\in\mathbb R^d:P+t=P\}
\]
contains a full-rank lattice.  Its stabilizer is then itself a full-rank
lattice, called the \emph{period lattice} of \(P\).  If
\(\mathbf L\subseteq\operatorname{Stab}(P)\) is a full-rank lattice and
\(D_\mathbf L\) a half-open fundamental domain, define the \emph{density}
of \(P\) by
\[
\operatorname{dens}(P):=
\frac{\#(P\cap D_\mathbf L)}{\operatorname{covol}(\mathbf L)}.
\]
This is independent of the lattice and fundamental domain.  Associate to $P$ the probability measure
\[
\eta_P
:=
\frac1{\operatorname{covol}(\mathbf L)}
\int_{D_\mathbf L}\delta_{\theta_t\delta_P}\,dt.
\]
The measure \(\eta_P\) is the unique translation-invariant probability
measure supported on the translation orbit of $\delta_P$, has intensity \(\operatorname{dens}(P)\), and is
independent of the chosen lattice
\(\mathbf L\subseteq\operatorname{Stab}(P)\).

A function \(h:\operatorname{Conf}(\mathbb R^d)\to\mathbb R\) is
\emph{local} if there is a bounded Borel set \(W\subset\mathbb R^d\)
such that \(h(\omega)=h(\omega')\) whenever
\(\omega|_W=\omega'|_W\). Such a set \(W\) is a
\emph{dependence window} for \(h\).

\begin{lemma}[Periodic approximation]\label{lem:periodic-block-approx}
Let \(\mathbf L_0<\mathbb R^d\) be a full-rank lattice, let
\(s_n\to\infty\), and set \(\mathbf L_n:=s_n\mathbf L_0\).  Suppose
\(N_n\in\mathbb N\) and
\[
\rho_n:=\frac{N_n}{\operatorname{covol}(\mathbf L_n)}\leq\rho,
\qquad
\rho_n\longrightarrow\rho\qquad(n\to\infty).
\]
Let \(\eta\in\PP_{\rho,2}\), let
\(h_1,\ldots,h_k\) be bounded continuous local functions, let
\(Q_1',\ldots,Q_\ell'\) be bounded half-open cubes, and let
\(\varepsilon>0\).  There is \(n_0\) such that, for every
\(n\geq n_0\), there is a simple \(\mathbf L_n\)-periodic set \(P_n\)
with exactly \(N_n\) points in a fundamental domain and
\[
\left|\int h_i\,d\eta_{P_n}-\int h_i\,d\eta\right|<\varepsilon
\quad(1\leq i\leq k),
\]
\[
\left|\int N_{Q_j'}^2\,d\eta_{P_n}-\int N_{Q_j'}^2\,d\eta\right|<\varepsilon
\quad(1\leq j\leq\ell).
\]
The measure \(\eta_{P_n}\) has intensity \(\rho_n\).
\end{lemma}

\begin{proof}
Write $V_n:=\operatorname{covol}(\mathbf L_n)$ and choose a half-open
fundamental parallelepiped $D$ of $\mathbf L_0$, so $s_nD$ is one for
$\mathbf L_n$.  Choose $r>0$ such that all dependence windows of the
functions $h_i$ and all the cubes $Q_j'$ are contained in $B_r$.
Fix $L>10r$.  For an $L$-cube $Q=z+[0,L)^d$, define
\[
Q^{\mathrm{ret}}:=z+[2r,L-2r)^d,\qquad
Q^{\mathrm{safe}}:=z+[3r,L-3r)^d.
\]
Their relative volumes are
\[
 b_L:=\left(1-\frac{4r}{L}\right)^d,\qquad
 a_L:=\left(1-\frac{6r}{L}\right)^d.
\]
For $t\in Q^{\mathrm{safe}}$ we have $t+B_r\subset Q^{\mathrm{ret}}$.
For large \(n\), enumerate the cubes \(Lk+[0,L)^d\),
\(k\in\mathbb Z^d\), contained in \(s_nD\) as
\(Q_1,\ldots,Q_J\), with \(J=J(n,L)\).
Their union leaves uncovered volume
\(O(L/s_n)V_n\), so
\[
 \frac{JL^d}{V_n}=1+O(L/s_n).
\]
On $\operatorname{Conf}(\mathbb R^d)^J$, put
$\mathfrak m_{n,L}:=\eta^{\otimes J}$ and denote the coordinate
configurations by $\omega_1,\ldots,\omega_J$.
Define $P_n^0$ on the fundamental domain by
\[
P_n^0\cap s_nD:=\bigcup_{a=1}^J
(\operatorname{supp}\omega_a\cap Q_a^{\mathrm{ret}}),
\]
and extend it $\mathbf L_n$-periodically.  Put
\(U_n:=\bigcup_{a=1}^J Q_a^{\mathrm{safe}}\).
Distinct retained cubes, including their periodic translates, are
separated by at least $4r$.  For $t\in Q_a^{\mathrm{safe}}$, each test depends only on the points
in $t+B_r\subset Q_a^{\mathrm{ret}}$.  Thus
\[
 h_i(\theta_t\delta_{P_n^0})=h_i(\theta_t\omega_a),\qquad
 N_{Q_j'}(\theta_t\delta_{P_n^0})^2
       =N_{Q_j'}(\theta_t\omega_a)^2.
\]
For \(1\le a\le J\), define
\[
 X_{a,i}=L^{-d}\int_{Q_a^{\mathrm{safe}}}
                  h_i(\theta_t\omega_a)\,dt,
 \qquad
 Y_{a,j}=L^{-d}\int_{Q_a^{\mathrm{safe}}}
                  N_{Q_j'}(\theta_t\omega_a)^2\,dt,
 \qquad
 Z_a=N_{Q_a^{\mathrm{ret}}}(\omega_a).
\]
The definitions and the local identities above give
\[
\begin{aligned}
\frac1{V_n}\int_{U_n}h_i(\theta_t\delta_{P_n^0})\,dt
 &=\frac{L^d}{V_n}\sum_{a=1}^J X_{a,i},\\
\frac1{V_n}\int_{U_n}N_{Q_j'}(\theta_t\delta_{P_n^0})^2\,dt
 &=\frac{L^d}{V_n}\sum_{a=1}^J Y_{a,j}.
\end{aligned}
\]
For fixed \(L\), the vectors
\((X_{a,1},\ldots,X_{a,k},Y_{a,1},\ldots,Y_{a,\ell},Z_a)\)
are independent and identically distributed over \(a\).
Translation invariance gives
\[
 \int X_{a,i}\,d\mathfrak m_{n,L}=a_L\int h_i\,d\eta,\qquad
 \int Y_{a,j}\,d\mathfrak m_{n,L}=a_L\int N_{Q_j'}^2\,d\eta,\qquad
 \int Z_a\,d\mathfrak m_{n,L}=\rho b_LL^d.
\]
Define \(K_n:=\#(P_n^0\cap s_nD)=\sum_{a=1}^J Z_a\).
The variables \(X_{a,i}\) are bounded, while \(Y_{a,j}\) and \(Z_a\)
are integrable by the local second-moment assumption.  For fixed
\(L\), we have \(J(n,L)\to\infty\), so the weak law of large numbers
applies to their sample averages.  By the weak law, choose $\delta_n\downarrow0$ so that the event
$\mathscr A_{n,L}$ on which the sample means of $X_{a,i}$,
$Y_{a,j}$ and $Z_a$ differ from their respective expectations by
at most $\delta_n$ satisfies
$\mathfrak m_{n,L}(\mathscr A_{n,L})\to1$.
On $\mathscr A_{n,L}$,
\[
 \frac{K_n}{V_n}
   =\rho b_L\frac{JL^d}{V_n}+o(1)
   =\rho b_L+O(L/s_n)+o(1),
\]
where the $o(1)$ term is bounded by a constant multiple of $\delta_n$
for fixed $L$.

The relative volume of \(s_nD\setminus U_n\) is
\(O(r/L)+O(L/s_n)\), so its contribution
to each bounded $h_i$ is bounded deterministically by this fraction
times $\|h_i\|_\infty$.  For the integrals of $N_{Q_j'}^2$, the
separation of the cubes $Q_a^{\mathrm{ret}}$ implies that
$t+Q_j'$ intersects at most one of them.  For each $a$, the set of origins $t$
outside the safe cubes for which $t+Q_j'$ intersects
$Q_a^{\mathrm{ret}}$ has volume $O(rL^{d-1})$.
On that set the square count of $P_n^0$ is bounded by
$N_{Q_j'}(\theta_t\omega_a)^2$.  Translation invariance and Fubini give
\[
 \int\frac1{V_n}
   \int_{s_nD\setminus U_n}
      N_{Q_j'}(\theta_t\delta_{P_n^0})^2\,dt\,d\mathfrak m_{n,L}
 \le C_j\frac rL\int N_{Q_j'}^2\,d\eta.
\]
The integrals over $U_n$ differ from the corresponding
integrals against $\eta$ by
$O(1-a_L)+O(L/s_n)+o(1)$ on $\mathscr A_{n,L}$.
Choose $L$ so that the deterministic losses and the sum of the
Markov bounds for the complement integrals are small. For sufficiently
large $n$, $\mathfrak m_{n,L}(\mathscr A_{n,L})>1/2$ and the
union of the excluded Markov events has measure less than $1/2$. The intersection therefore contains a sample satisfying the test
inequalities. For fixed \(L\), the difference \(\rho-\rho b_L=\rho(1-b_L)\) is
positive.  Since \(\rho_n\to\rho\), the bound on \(K_n/V_n\)
gives \(K_n/V_n<\rho_n=N_n/V_n\) for sufficiently large \(n\)
on $\mathscr A_{n,L}$.
For this sample, put \(M_n:=N_n-K_n\).  Then
\[
 0\le \frac{M_n}{V_n}
 =\rho(1-b_L)+O(L/s_n)+o(1)
 \le 2\rho(1-b_L).
\]

If $M_n>0$, let $k_n$ be the least integer with $k_n^d\ge M_n$.
The grid $k_n^{-1}\mathbf L_n/\mathbf L_n$ has $k_n^d$ points.
Since $k_n\leq2M_n^{1/d}$, its separation is at least
\[
 c_{\mathbf L_0}\left(\frac{V_n}{M_n}\right)^{1/d},
 \qquad
 c_{\mathbf L_0}:=
 \frac{\lambda_1(\mathbf L_0)}
 {2\operatorname{covol}(\mathbf L_0)^{1/d}}.
\]
The estimates above permit $L$ to be fixed sufficiently large, and
then $n$ sufficiently large, that $M_n/V_n\leq2\rho(1-b_L)$ and
$c_{\mathbf L_0}(2\rho(1-b_L))^{-1/d}>2r$.
Choose a translate of the grid disjoint from \(P_n^0\) modulo
\(\mathbf L_n\), and let \(G_n\) consist of \(M_n\) of its points,
extended \(\mathbf L_n\)-periodically.  Set \(P_n:=P_n^0\cup G_n\).
This set is simple and \(\mathbf L_n\)-periodic, and
\(\#(P_n\cap s_nD)=K_n+M_n=N_n\).

For a bounded Borel window \(Q\), put
\(I_Q(t):=\#(G_n\cap(t+Q))\).
For each fixed $Q$, separation gives $\sup_t I_Q(t)\le C_Q$, with
$C_Q$ independent of sufficiently large $L,n$.  Since $I_Q(t)^2\le C_Q I_Q(t)$, integration over the period
cell gives
\begin{equation}\label{eq:periodic-insertion-count}
 \frac1{V_n}\int_{s_nD} I_Q(t)^2\,dt
 \le C_Q\,|Q|\frac{M_n}{V_n}.
\end{equation}
Let \(W_i\) be a bounded dependence window for \(h_i\).  If
\((t+W_i)\) contains no inserted point, then
\(h_i(\theta_t\delta_{P_n})=h_i(\theta_t\delta_{P_n^0})\).
Consequently,
\[
\left|\int h_i\,d\eta_{P_n}-\int h_i\,d\eta_{P_n^0}\right|
\leq 2\|h_i\|_\infty\,\Vol_d(W_i)\frac{M_n}{V_n}.
\]
For the functions $N_{Q_j'}^2$, for every $t$,
\[
N_{Q_j'}(\theta_t\delta_{P_n})
=N_{Q_j'}(\theta_t\delta_{P_n^0})+I_{Q_j'}(t).
\]
Put $C(t):=N_{Q_j'}(\theta_t\delta_{P_n^0})$ and
$I(t):=I_{Q_j'}(t)$.  Cauchy--Schwarz gives
\[
\begin{split}
\frac1{V_n}\int_{s_nD}\bigl((C+I)^2-C^2\bigr)\,dt
&\leq 2\left(\frac1{V_n}\int_{s_nD}C^2\,dt\right)^{1/2}
\left(\frac1{V_n}\int_{s_nD}I^2\,dt\right)^{1/2}\\
&\qquad+\frac1{V_n}\int_{s_nD}I^2\,dt.
\end{split}
\]
The choice of $P_n^0$ and \eqref{eq:periodic-insertion-count} give
\[
\begin{aligned}
\frac1{V_n}\int_{s_nD}C(t)^2\,dt
&\leq\int N_{Q_j'}^2\,d\eta+\varepsilon,\\
\frac1{V_n}\int_{s_nD}I(t)^2\,dt
&=O(M_n/V_n).
\end{aligned}
\]
Since $M_n/V_n=O(r/L)+O(L/s_n)+o(1)$, the difference of
the $N_{Q_j'}^2$-integrals caused by the $M_n$ inserted points is
$O((r/L)^{1/2})+o(1)$ as $n\to\infty$
with $L$ fixed.  Taking $L$ large and then $n$ large gives both
required approximations.  The point count per fundamental
domain gives the asserted intensity $\rho_n$.
\end{proof}

Fix a countable family
\((f_k)_{k\geq1}\subset C_c^+(\mathbb R^d)\) such that, for every
\(m\geq1\), every \(f\in C_c^+(\mathbb R^d)\) with
\(\operatorname{supp}f\subset B_m\) is a uniform limit of functions from
this family supported in \(B_{m+1}\), and put
\[
h_k(\omega):=\exp(-S_{f_k}(\omega)).
\]
Each \(h_k\) is a bounded continuous local function.  If
\(\|f-f_k\|_\infty<\varepsilon\) and both functions are supported in
\(B_{m+1}\), then, for every translation-invariant measure \(\eta'\) of intensity at most
\(\rho\),
\[
\left|\int e^{-S_f}\,d\eta'-\int e^{-S_{f_k}}\,d\eta'\right|
\leq \varepsilon\rho\Vol_d(B_{m+1}).
\]
The Laplace-functional criterion \cite{Kal02} shows that \((h_k)\)
determines weak convergence of point-process laws.  Applying
\cref{lem:periodic-block-approx} successively to
\(h_1,\ldots,h_n\) and \(Q_1,\ldots,Q_n\) gives weak convergence together
with convergence of the required local second moments.

A process \(\eta\in\EPP_{\rho,2}\) is \emph{periodic} if there
is a full-rank lattice \(\mathbf L<\mathbb R^d\) such that
\(\theta_t\omega=\omega\) for every \(t\in\mathbf L\), for
\(\eta\)-almost every \(\omega\).  Let \(\PerPP_\rho\) be the class
of periodic processes in \(\EPP_{\rho,2}\).

\begin{proposition}\label{prop:periodic-density}
The class \(\PerPP_\rho\) is contained in \(\Hyp_\rho\) and is dense
in \(\PP_{\rho,2}\).  Consequently
\(\Hyp_\rho\cap\EPP_{\rho,2}\) is dense in \(\EPP_{\rho,2}\).
\end{proposition}

\begin{proof}
Fix $\eta\in\PP_{\rho,2}$. For each $n$, apply
\cref{lem:periodic-block-approx} to the full sequence
\[
s_j:=(j/\rho)^{1/d},\qquad
\mathbf L_j:=s_j\mathbb Z^d,\qquad N_j:=j\quad(j\geq1),
\]
using the first $n$ functions $h_1,\ldots,h_n$, the first
$n$ cubes $Q_1,\ldots,Q_n$, and tolerance $1/n$.
The lemma gives $j_n$ such that every $j\geq j_n$ admits
the required periodic approximation. Choose increasing
integers $q_n\geq j_n$, and put
$S_n:=(q_n/\rho)^{1/d}$,
$\mathbf L_n:=S_n\mathbb Z^d$.
Let $P_n$ be a $\mathbf L_n$-periodic set with exactly $q_n$
points per fundamental domain.
Thus
$\operatorname{dens}(P_n)=\rho$ and $\eta_{P_n}\in\PerPP_\rho$.  The approximations of the Laplace functionals $h_k$ give weak
convergence.  The second-moment estimates then give
$\eta_{P_n}\to\eta$ in $\PP_{\rho,2}$ by
\cref{prop:local-L2-characterization}.  Thus $\PerPP_\rho$ is dense.

Fix $\eta\in\PerPP_\rho$ and choose a full-rank lattice
$\mathbf L$ for which $\theta_t\omega=\omega$ for every
$t\in\mathbf L$, for $\eta$-almost every $\omega$.
Let $m_{\mathbf L}$ be Haar probability on
$\mathbb R^d/\mathbf L$.  For every bounded continuous function $h$
on configuration space, the orbit average
$h^\sharp(\omega):=\int h(\theta_t\omega)\,dm_{\mathbf L}(t)$
is invariant.  Ergodicity makes this average equal to $\int h\,d\eta$
for $\eta$-almost every $\omega$.  Applying this to a countable
determining family shows that $\eta=\eta_P$ for a
$\mathbf L$-periodic set $P$.  Its density equals $\rho$.
Choose a half-open fundamental domain $D_\mathbf L$ for $\mathbf L$ and
put $P_0:=P\cap D_\mathbf L$.  For $\xi\in\mathbf L^*$, set
\[
A(\xi):=\sum_{p\in P_0}e^{-2\pi i\langle\xi,p\rangle}.
\]
For $f\in C_c(\mathbb R^d)$, expand $t\mapsto S_f(\delta_{P+t})$ in Fourier series on $\mathbb R^d/\mathbf L$.
After subtracting the zero-frequency term, the coefficient at
$\xi\in\mathbf L^*\setminus\{0\}$ has absolute value
$|A(\xi)\widehat f(\xi)|/\operatorname{covol}(\mathbf L)$.
Orthogonality of the characters and Parseval's identity give
\[
\sigma_{\eta_P}
=
\frac1{\operatorname{covol}(\mathbf L)^2}
\sum_{\xi\in\mathbf L^*\setminus\{0\}}
|A(\xi)|^2\,\delta_\xi .
\]

Parseval's identity on \(\mathbb R^d/\mathbf L\) gives
\[
V_{\eta_P}(R)
=
\frac1{\operatorname{covol}(\mathbf L)^2}
\sum_{\xi\in\mathbf L^*\setminus\{0\}}
|A(\xi)|^2\,|\widehat{\mathbf1}_{B_R}(\xi)|^2.
\]
Since \(|A(\xi)|\leq\rho\operatorname{covol}(\mathbf L)\), for
\(R\geq\lambda_1(\mathbf L^*)^{-1}\) the Bessel-function estimate
\[
|\widehat{\mathbf1}_{B_R}(\xi)|
\leq
C_dR^{(d-1)/2}|\xi|^{-(d+1)/2}
\qquad(R|\xi|\geq1)
\]
(see \cite[Sec.~IV and Appendix~A]{TS03}) yields
\[
V_{\eta_P}(R)
\leq
C_{P,d}R^{d-1}
\sum_{\xi\in\mathbf L^*\setminus\{0\}}|\xi|^{-(d+1)}
=O_P(R^{d-1}).
\]
Thus every \(\eta\in\PerPP_\rho\) is hyperuniform.
\end{proof}

\emph{Proof of \cref{thm:main-hyperuniformity}.}
The residuality of \(\IHyp_\rho\) and the meagreness of
\(\Hyp_\rho\cap\EPP_{\rho,2}\) were proved at the end of
\cref{sec:intermittent-hyperuniformity}; \cref{prop:periodic-density}
provides the remaining density assertion.

\subsection{Pairwise negative correlation}
\label{subsec:pnc-structure}

The spectral representation of a reduced covariance measure belongs
to the Bartlett-spectrum framework; see \cite[Sec.~8.2]{DVJ03}.
For pairwise negative correlation we show that its off-diagonal part
is a finite negative measure.  This yields a bounded continuous
spectral density and the variance limit used below.

A process \(\eta\in\PP_{\rho,2}\) is \emph{pairwise negatively correlated}
if
\[
\Cov_\eta(N_A,N_B)\leq0
\]
for every pair of disjoint bounded Borel sets \(A,B\subset\mathbb R^d\).
Let \(\PNC_\rho\) denote this class.  This condition is weaker than
negative association \cite{TPBG25}.

\begin{remark}\label{rem:KLLY-pnc}
The fair-partition construction of \cite[Example~4.7]{KLLY25} belongs to
\(\PNC_\rho\).  Conditionally on the partition \((C_i)_i\), set
\[
p_i(A):=\rho\Vol_d(A\cap C_i).
\]
For disjoint bounded Borel sets \(A,B\), conditional independence between
cells and the fact that each cell contributes at most one point give
\[
\Cov(N_A,N_B\mid(C_i)_i)
=-\sum_i p_i(A)p_i(B).
\]
Since
\(\mathbb E[N_A\mid(C_i)_i]=\smash{\sum_i} p_i(A)=\rho\Vol_d(A)\) is deterministic,
\[
\Cov(N_A,N_B)
=-\mathbb E\sum_i p_i(A)p_i(B)\leq0.
\]
\end{remark}

Recall the definition of a F\o lner sequence from
\cref{subsec:intro-van-hove}.  A process \(\eta\in\PP_{\rho,2}\) is
\emph{F\o lner-hyperuniform} if
\[
\frac{\Var_\eta(N_{F_n})}{\Vol_d(F_n)}
\longrightarrow0\qquad(n\to\infty)
\]
for every F\o lner sequence \((F_n)\).  Denote this class by
\(\FHyp_\rho\).

\begin{proposition}\label{prop:pnc-structure}
The class \(\PNC_\rho\) is closed in \(\PP_{\rho,2}\).  For every
\(\eta\in\PNC_\rho\) there is a finite symmetric positive measure
\(\tau_\eta\) on \(\mathbb R^d\), with
\[
\tau_\eta(\mathbb R^d)\leq\rho,
\]
such that, for real bounded Borel functions $f,g$ with compact support,
\[
\Cov_\eta(S_f,S_g)
=
\rho\int f(u)g(u)\,du
-
\int_{\mathbb R^d}\int_{\mathbb R^d}
f(u)g(u+t)\,du\,d\tau_\eta(t).
\]
Equivalently,
\[
\Cov_\eta(S_f,S_g)
=
\int_{\mathbb R^d}
\widehat f(\xi)\overline{\widehat g(\xi)}
s_\eta(\xi)\,d\xi,
\qquad
s_\eta(\xi):=\rho-\widehat{\tau_\eta}(\xi),
\]
so
\[
\sigma_\eta(d\xi)=s_\eta(\xi)\,d\xi,
\qquad 0\leq s_\eta\leq2\rho,
\]
with \(s_\eta\) bounded and continuous.  Moreover, for every F\o lner
sequence \((F_n)\),
\[
\frac{\Var_\eta(N_{F_n})}{\Vol_d(F_n)}
\longrightarrow
\rho-\tau_\eta(\mathbb R^d)=s_\eta(0)
\qquad(n\to\infty).
\]
Consequently, for \(\eta\in\PNC_\rho\), the following are equivalent:
\[
\eta\in\Hyp_\rho,
\qquad
\eta\in\FHyp_\rho,
\qquad
s_\eta(0)=0.
\]
Finally, for real \(f,g\in C_c(\mathbb R^d)\),
\[
\Cov_\eta(S_f,S_g\circ\theta_t)\longrightarrow0
\qquad(|t|\to\infty).
\]
\end{proposition}

\begin{proof}
Let \(\mathcal R\) be the countable ring generated by half-open rectangles
with rational endpoints.  Translation invariance and \cref{lem:null-window} imply
that every \(A\in\mathcal R\) is a continuity set for every
\(\eta\in\PP_{\rho,2}\).  Therefore, for disjoint \(A,B\in\mathcal R\),
\[
\eta\longmapsto\Cov_\eta(N_A,N_B)
\]
is continuous in the local \(L^2\) topology.  Let
\(\smash{\alpha_\eta^{(2)}}\) be the second factorial moment measure,
\[
\int H(x,y)\,d\alpha_\eta^{(2)}(x,y)
:=
\int
\sum_{\substack{p,q\in\operatorname{supp}\omega\\p\neq q}}
H(p,q)\,d\eta(\omega),
\qquad H:\mathbb R^{2d}\to[0,\infty]\text{ Borel}.
\]
For disjoint \(A,B\),
\[
\Cov_\eta(N_A,N_B)\leq0
\quad\Longleftrightarrow\quad
\alpha_\eta^{(2)}(A\times B)
\leq
\rho^2\Vol_d(A)\Vol_d(B).
\]
Thus the inequalities for disjoint \(A,B\in\mathcal R\) hold on a
closed subset of \(\PP_{\rho,2}\).  To prove that they imply the
inequalities for all disjoint bounded Borel sets, let
\(H\in C_c(\mathbb R^d\times\mathbb R^d)\) be nonnegative with support
disjoint from the diagonal.  Partition a large box containing its support
into sufficiently small sets from \(\mathcal R\).  Only products of
distinct cells intersect the support of \(H\), so step-function approximation
and the rectangle inequalities give
\[
\int H\,d\alpha_\eta^{(2)}
\leq
\rho^2\int H(x,y)\,dx\,dy.
\]
The inequality holds for every such $H$.  Regularity of the
measures therefore gives
\[
\alpha_\eta^{(2)}\leq\rho^2\,dx\,dy
\qquad\text{off the diagonal}.
\]
In particular, the covariance inequality holds for all disjoint bounded
Borel sets, so \(\PNC_\rho\) is closed.

Since $\eta$ is simple, the factorial moment measure has no mass
on the diagonal.  Thus the off-diagonal inequality also gives
$\alpha_\eta^{(2)}\leq\rho^2dx\,dy$ on all of
$\mathbb R^{2d}$.  Define the positive measure
\[
\Theta_\eta:=\rho^2\,dx\,dy-\alpha_\eta^{(2)}.
\]
It is invariant under simultaneous translation of $x$ and $y$.
Let $\Theta_\eta^\Delta$ be the push-forward of $\Theta_\eta$
under $(x,y)\mapsto(x,y-x)$.  For each bounded Borel
$C\subset\mathbb R^d\setminus\{0\}$, the measure
\[
A\longmapsto\Theta_\eta^\Delta(A\times C)
\]
is locally finite and translation-invariant.  It is therefore a multiple
of Lebesgue measure.  Given a bounded Borel set $U$ of positive volume,
define
\[
\tau_\eta(C):=\frac{\Theta_\eta^\Delta(U\times C)}{\Vol_d(U)}.
\]
Translation invariance makes this independent of $U$.
It defines a locally finite positive measure on
$\mathbb R^d\setminus\{0\}$. Set $\tau_\eta(\{0\})=0$. Then
\[
\Theta_\eta^\Delta(dx,dt)=dx\,d\tau_\eta(t),
\qquad
 d\Theta_\eta(x,y)=dx\,d\tau_\eta(y-x).
\]
Symmetry of \(\smash{\alpha_\eta^{(2)}}\) makes \(\tau_\eta\) symmetric.

For a bounded Borel set $F$ with $0<\Vol_d(F)<\infty$, the factorial-moment identity gives
\[
\Var_\eta(N_F)
=
\rho\Vol_d(F)+\alpha_\eta^{(2)}(F\times F)
-\rho^2\Vol_d(F)^2
=
\rho\Vol_d(F)-\Theta_\eta(F\times F).
\]
Using $d\Theta_\eta(x,y)=dx\,d\tau_\eta(y-x)$, we obtain
\[
\Theta_\eta(F\times F)
=
\int_{\mathbb R^d}\Vol_d(F\cap(F-t))\,d\tau_\eta(t).
\]
Therefore
\begin{equation}\label{eq:pairwise-negative-variance}
\frac{\Var_\eta(N_F)}{\Vol_d(F)}
=
\rho-
\int_{\mathbb R^d}
c_F(t)\,d\tau_\eta(t),
\qquad
c_F(t):=
\frac{\Vol_d(F\cap(F-t))}{\Vol_d(F)}.
\end{equation}
For expanding cubes $Q_R$, the functions $c_{Q_R}$ satisfy
$0\le c_{Q_R}\le1$ and $c_{Q_R}(t)\to1$ for every $t$.
Since $\Var_\eta(N_{Q_R})\geq0$, \eqref{eq:pairwise-negative-variance} gives
$\int c_{Q_R}\,d\tau_\eta\leq\rho$.
Fatou's lemma now yields
\[
\tau_\eta(\mathbb R^d)\leq\rho.
\]

For bounded Borel functions $f,g$ with compact support, the factorial moment identity is
\[
\int S_fS_g\,d\eta
=\rho\int f(u)g(u)\,du
+\int f(x)g(y)\,d\alpha_\eta^{(2)}(x,y).
\]
Substituting
$\smash{\alpha_\eta^{(2)}}=\rho^2dx\,dy-\Theta_\eta$ and subtracting the
product of the first moments gives
\[
\Cov_\eta(S_f,S_g)
=
\rho\int fg
-
\int\!\!\int f(u)g(u+t)\,du\,d\tau_\eta(t).
\]
For real $f,g\in C_c^\infty$, Fourier inversion and Plancherel give the stated
density $s_\eta=\rho-\widehat{\tau_\eta}$.  The corresponding identity for
complex-valued smooth tests follows by sesquilinearity.
For bounded Borel functions $f,g$ with compact support, approximate them
in $L^2(dx)$ by smooth compactly supported functions
with supports in one fixed compact set. As $n\to\infty$,
Cauchy--Schwarz gives, uniformly in $t$,
\[
\left|\int f(u)g(u+t)\,du-\int f_n(u)g_n(u+t)\,du\right|
\leq
\|f-f_n\|_2\|g\|_2+\|f_n\|_2\|g-g_n\|_2
\longrightarrow0.
\]
Since \(\tau_\eta\) is finite, the covariance formula passes
to the limit and Plancherel gives the corresponding Fourier-side identity;
since \(\tau_\eta\) is also symmetric, \(s_\eta\) is bounded, continuous,
and real-valued.  The
covariance form is nonnegative when \(f=g\), so \(s_\eta\geq0\) almost
everywhere; continuity gives \(s_\eta\geq0\) everywhere.  In addition,
\[
s_\eta(\xi)
\leq
\rho+|\widehat{\tau_\eta}(\xi)|
\leq2\rho.
\]
For $f,g\in C_c^\infty$, this is the Bartlett identity with spectral
measure $s_\eta\,d\xi$.  Approximation in the Schwartz topology,
using translation boundedness, and uniqueness in
\cref{prop:bartlett-spectrum} give $\sigma_\eta=s_\eta\,d\xi$.

If $(F_n)$ is a F\o lner sequence, then
\[
c_{F_n}(t)
=
1-\frac{\Vol_d(F_n\mathbin{\triangle}(F_n-t))}
        {2\Vol_d(F_n)}
\longrightarrow1
\]
for every \(t\).  Dominated convergence with respect to the finite measure
\(\tau_\eta\) proves the F\o lner variance limit.

Since the balls form a F\o lner sequence, \(\FHyp_\rho\subseteq\Hyp_\rho\).
If \(\eta\) is hyperuniform, the F\o lner variance limit along balls gives
\(s_\eta(0)=0\); conversely, \(s_\eta(0)=0\) makes the displayed limit
zero along every F\o lner sequence.  This proves the three-way equivalence.

For real \(f,g\in C_c(\mathbb R^d)\), the product
\(s_\eta\widehat f\,\overline{\widehat g}\) belongs to \(L^1(\mathbb R^d)\).
The Riemann--Lebesgue lemma therefore gives
\[
\Cov_\eta(S_f,S_g\circ\theta_t)\longrightarrow0
\qquad(|t|\to\infty).
\]
\end{proof}

\begin{corollary}\label{cor:pnc-nowhere-dense}
The class \(\PNC_\rho\) is nowhere dense in \(\PP_{\rho,2}\).
\end{corollary}

\begin{proof}
By \cref{prop:pnc-structure}, \(\PNC_\rho\) is closed.  Let \(P\) be a periodic set with \(\operatorname{dens}(P)=\rho\).  Its orbit
measure \(\eta_P\) is not a point mass, so
\(\Var_{\eta_P}(S_f)>0\) for some \(f\in C_c(\mathbb R^d)\).  If
\(0\neq\lambda\in\operatorname{Stab}(P)\), then
\[
\Cov_{\eta_P}(S_f,S_f\circ\theta_{n\lambda})
=
\Var_{\eta_P}(S_f)>0
\qquad(n\geq1),
\]
contradicting the final assertion of \cref{prop:pnc-structure}.  Hence
\(\PerPP_\rho\cap\PNC_\rho=\varnothing\).  Since
\(\PerPP_\rho\) is dense by \cref{prop:periodic-density}, the closed
set \(\PNC_\rho\) has empty interior.
\end{proof}

\section{Revival and near-quadratic van Hove fluctuations}
\label{sec:revival}

\subsection{Revival sequences for point processes}\label{subsec:point-revival}

\begin{proposition}\label{prop:revival-residual}
The set \(\Rev_\rho\) is residual in \(\smash{\EPP_{\rho,2}}\).
\end{proposition}

\begin{proof}
Let \((A_j)_{j\geq1}\) enumerate the countable
ring of finite unions of bounded half-open rectangles with rational
endpoints.  For
\(m,\ell\geq1\), define
\[
\mathcal V_{m,\ell}
:=
\bigcup_{|t|>\ell}
\bigcap_{j=1}^m
\left\{
\eta\in\smash{\EPP_{\rho,2}}:
\|N_{A_j+t}-N_{A_j}\|_{L^2(\eta)}<m^{-1}
\right\}.
\]
Fix $j$ and $t$, and suppose $\eta_k\to\eta$ locally in $L^2$.
Choose $m$ with $A_j\cup(A_j+t)\subset Q_m$.
By \cref{lem:null-window}, both count functions are almost surely
continuous under $\eta$ for the vague topology.
The expectations of $S_{\varphi_m}^2$ converge by definition of the
topology. Since $S_{\varphi_m}^2$ is continuous in the vague topology,
weak convergence and convergence of its integrals give
\[
\lim_{M\to\infty}\sup_{k\geq1}
\int S_{\varphi_m}^2\mathbf1_{\{S_{\varphi_m}^2>M\}}\,d\eta_k=0.
\]
The squared count difference is almost surely continuous under $\eta$
and bounded by $4S_{\varphi_m}^2$.
Thus the expectation of
$(N_{A_j+t}-N_{A_j})^2$ is continuous in $\eta$.
Each $\mathcal V_{m,\ell}$ is therefore open.

If $P$ is periodic, then for $t\in\operatorname{Stab}(P)$ the
configuration itself is unchanged by translation by $t$.  Thus
$N_{A+t}=N_A$ $\eta_P$-almost surely.  Since
\(\operatorname{Stab}(P)\) contains vectors of arbitrarily large norm,
every periodic invariant measure belongs to every \(\mathcal V_{m,\ell}\).
By \cref{prop:periodic-density}, each \(\mathcal V_{m,\ell}\) is dense.  Thus
\[
\mathcal V:=\bigcap_{m,\ell\geq1}\mathcal V_{m,\ell}
\]
is a dense $G_\delta$ subset of \(\smash{\EPP_{\rho,2}}\).
Fix $\eta\in\mathcal V$.  For each $n$, the condition
$\eta\in\mathcal V_{n,n}$ gives a vector $t_n$ with $|t_n|>n$
and
\[
\|N_{A_j+t_n}-N_{A_j}\|_{L^2(\eta)}<n^{-1}
\qquad(1\leq j\leq n).
\]
Let $A$ be a bounded Borel set.  Choose a bounded cube $D\supset A$.
By approximation in Lebesgue measure, there are ring sets $B_k\subset D$
such that
\[
\Vol_d(A\mathbin\triangle B_k)\longrightarrow0
\qquad(k\to\infty).
\]

Translation invariance gives
\[
\int N_{A\triangle B_k}\,d\eta
=\rho\Vol_d(A\triangle B_k)\longrightarrow0\qquad(k\to\infty).
\]
Moreover $0\leq N_{A\triangle B_k}\leq N_D$ and
$N_D\in L^2(\eta)$.  Markov's inequality gives $N_{A\triangle B_k}\to0$ in probability.
Since $N_{A\triangle B_k}^2\leq N_D^2\in L^1(\eta)$,
uniform integrability upgrades convergence in probability to
\[
\|N_{A\triangle B_k}\|_{L^2(\eta)}\longrightarrow0
\qquad(k\to\infty).
\]
For any such ring set $B$ and every $n$, translation invariance and the triangle
inequality give
\[
\|N_{A+t_n}-N_A\|_2
\leq
\|N_{(A\triangle B)+t_n}\|_2
+\|N_{B+t_n}-N_B\|_2
+\|N_{A\triangle B}\|_2
\leq
2\|N_{A\triangle B}\|_2
+\|N_{B+t_n}-N_B\|_2.
\]
For each fixed $B=B_k$, the last term tends to zero as
$n\to\infty$ by the defining inequalities for the $A_j$.
Then let $k\to\infty$, using
$\|N_{A\triangle B_k}\|_2\to0$.  It follows that
$\|N_{A+t_n}-N_A\|_2\to0$ for every bounded Borel $A$.
Thus $(t_n)$ is a revival sequence and $\mathcal V\subseteq\Rev_\rho$.
\end{proof}

\subsection{Near-quadratic van Hove fluctuations}\label{subsec:van-hove-maximal}

Recall the van Hove condition and the notation \(\partial_rF\) from
\cref{subsec:intro-van-hove}.  Every van Hove sequence is F\o lner: for
fixed \(t\in\mathbb R^d\),
\[
(F+t)\mathbin{\triangle}F\subseteq\partial_{|t|+1}F,
\]
so the F\o lner ratio is bounded by the corresponding van Hove boundary
ratio.

\begin{proposition}\label{prop:van-hove-subquadratic}
For every \(\eta\in\smash{\EPP_{\rho,2}}\) and every van Hove sequence
\((F_n)\),
\[
\Var_\eta(N_{F_n})
=
o\!\left(\Vol_d(F_n)^2\right).
\]
\end{proposition}

\begin{proof}
Let $\varphi\in C_c(\mathbb R^d)$ be even and nonnegative, with
$\smash{\int}\varphi=1$ and $\operatorname{supp}\varphi\subset B_r$.
For a unitary $\mathbb R^d$-representation, the normalized averages
along a F\o lner sequence converge strongly to the invariant
projection.  To see this, note that the averages fix invariant vectors.
On vectors $H-U_sH$ they tend to zero by the F\o lner property;
the closed span of these differences is the orthogonal complement of
the invariant subspace.
Since $\eta$ is ergodic and $\int(S_\varphi-\rho)\,d\eta=0$, the
invariant projection of $S_\varphi-\rho$ vanishes.  Hence
\[
\frac{1}{\Vol_d(F_n)}
\int_{F_n}
\bigl(S_\varphi\circ\theta_t-\rho\bigr)\,dt
\longrightarrow0
\qquad\text{in }L^2(\eta)\text{ as }n\to\infty,
\]
by the F\o lner property.  By Fubini, the term
on the left is
\(Z_{\mathbf 1_{F_n}*\varphi}^\eta/\Vol_d(F_n)\).

Define
\[
g_n:=\mathbf 1_{F_n}-\mathbf 1_{F_n}*\varphi.
\]
We have $|g_n|\leq1$, and $g_n$ vanishes outside $\partial_rF_n$.
Partition $\mathbb R^d$ into half-open unit cubes $C$.  If such a
cube intersects $\partial_rF_n$, it lies inside
$\partial_{r+\sqrt d+1}F_n$ (up to its boundary), so the number of
contributing cubes is bounded by the volume of that enlarged
boundary neighbourhood.  For each unit cube $C$, translation invariance and $|g_n|\leq1$ give
\[
\|Z_{g_n\mathbf1_C}^\eta\|_2
\leq\|N_C\|_2+\rho\Vol_d(C)
=\|N_{Q_1}\|_2+\rho.
\]
Apply Minkowski's inequality to the contributing unit cubes and use
the uniform $L^2$ bound above.  There is a constant $C_\eta<\infty$
such that
\[
\|Z_{g_n}^\eta\|_{L^2(\eta)}
\leq
C_\eta\,\Vol_d(\partial_{r+\sqrt d+1}F_n).
\]
Since \(\Vol_d(\partial_{r+\sqrt d+1}F_n)/\Vol_d(F_n)\to0\),
\[
\frac{\|Z_{g_n}^\eta\|_{L^2(\eta)}}{\Vol_d(F_n)}\longrightarrow0
\qquad(n\to\infty).
\]
Since
\[
N_{F_n}-\rho\Vol_d(F_n)
=Z_{\mathbf 1_{F_n}*\varphi}^\eta+Z_{g_n}^\eta,
\]
the centered count divided by \(\Vol_d(F_n)\) converges to zero in
\(L^2(\eta)\).  Taking squared norms proves the claim.
\end{proof}

\begin{remark}[F\o lner versus van Hove sequences]
\label{rem:folner-versus-van-hove}
The conclusion of \cref{prop:van-hove-subquadratic} fails if ``van Hove''
is replaced by ``F\o lner''.  First take intensity \(1\) and let
\[
\int F\,d\eta
:=\int_{[0,1)^d}F(\delta_{\mathbb Z^d+u})\,du
\qquad(F\text{ bounded Borel}).
\]
Let $U$ denote the coordinate on $[0,1)^d$ with Lebesgue measure.
For \(\varepsilon_n:=n^{-3/2}\) and \(K_n:=n^{2d}\), define
\[
F_n
:=
[0,n)^d
\cup
\bigcup_{i=1}^{K_n}
\bigl((n+i)e_1+[0,\varepsilon_n)^d\bigr).
\]
The microscopic cubes have total volume
\(K_n\varepsilon_n^d=n^{d/2}\), so
\[
\Vol_d(F_n)=n^d+n^{d/2}\sim n^d.
\]
For every fixed \(t\in\mathbb R^d\),
\[
\Vol_d((F_n+t)\mathbin{\triangle}F_n)
\leq
C_{d,t}n^{d-1}+2n^{d/2}
=o(\Vol_d(F_n)),
\]
and hence \((F_n)\) is F\o lner.  Every microscopic cube contains a point
of \(\mathbb Z^d+U\) if and only if
\(U\in[0,\varepsilon_n)^d\), simultaneously for all \(i\), while
\([0,n)^d\) contains \(n^d\) points almost surely.  Thus
\[
N_{F_n}
=
n^d+K_n\mathbf1_{\{U\in[0,\varepsilon_n)^d\}},
\]
and therefore
\[
\Var_\eta(N_{F_n})
=
K_n^2\varepsilon_n^d(1-\varepsilon_n^d)
\sim n^{5d/2}.
\]
Consequently
\[
\frac{\Var_\eta(N_{F_n})}{\Vol_d(F_n)^2}
\sim n^{d/2}\longrightarrow\infty\qquad(n\to\infty).
\]
The divergence comes from an event of measure $\varepsilon_n^d=n^{-3d/2}$
on which the microscopic cubes contribute $K_n=n^{2d}$ points;
their contribution to the variance is of order
$K_n^2\varepsilon_n^d=n^{5d/2}$.
Dilating both the lattice and the sets by \(\rho^{-1/d}\) gives the same
conclusion at any prescribed intensity \(\rho>0\).

For $r=1/4$, the boundary neighbourhoods of the microscopic cubes
have total volume of order $n^{2d}$, whereas $\Vol_d(F_n)\sim n^d$.
Hence $(F_n)$ is not van Hove.
\end{remark}

\begin{lemma}\label{lem:revival-covariance}
Let \(\eta\in\smash{\EPP_{\rho,2}}\) and let
\(A\subset\mathbb R^d\) be bounded Borel.  Then, for \(s,t\in\mathbb R^d\),
\[
\Cov_\eta(N_{A+s},N_{A+t})
=
\Var_\eta(N_A)
-
\frac12\|N_{A+s}-N_{A+t}\|_{L^2(\eta)}^2.
\]
Consequently, if \(A+s_1,\ldots,A+s_m\) are pairwise disjoint and
\[
\|N_{A+s_i}-N_A\|_{L^2(\eta)}\leq\delta
\qquad(1\leq i\leq m),
\]
then, with \(\sigma_A:=\sqrt{\Var_\eta(N_A)}\),
\[
\sqrt{\Var_\eta\!\left(
N_{\bigsqcup_{i=1}^m(A+s_i)}
\right)}
\geq
m(\sigma_A-\delta).
\]
\end{lemma}

\begin{proof}
The random variables \(N_{A+s}\) and \(N_{A+t}\) have the same mean and
variance as \(N_A\).  Expanding
\(\|N_{A+s}-N_{A+t}\|_2^2\) gives the covariance identity.  For the second assertion, write
$Y_i:=N_{A+s_i}-\rho\Vol_d(A)$ and
$Y_0:=N_A-\rho\Vol_d(A)$.  By the reverse triangle inequality,
\[
\left\|\sum_{i=1}^mY_i\right\|_2
\geq m\|Y_0\|_2-\sum_{i=1}^m\|Y_i-Y_0\|_2
\geq m(\sigma_A-\delta).
\]
Disjointness identifies $\smash{\sum_i}Y_i$ with the centered count of the
union.
\end{proof}

\begin{proposition}\label{prop:revival-mfvh}
Let \(\eta\in\Rev_\rho\), let \(\alpha_n>0\) with \(\alpha_n\to0\), and
let \(M_n>0\) with \(M_n\to\infty\).  Then there is a van Hove sequence
\((F_n)\), depending on \(\eta\), \((\alpha_n)\), and \((M_n)\), such that
\[
\Var_\eta(N_{F_n})
\geq
M_n\Vol_d(F_n)^{2-\alpha_n}
\qquad(n\geq1).
\]
Consequently \(\Rev_\rho\subseteq\MFVH_\rho\).
\end{proposition}

\begin{proof}
Fix $\eta\in\Rev_\rho$ and a revival sequence $(t_k)$.
Choose $A$ from the ring of finite unions of bounded rational rectangles
such that
\[
\sigma^2:=\Var_\eta(N_A)>0.
\]
Indeed, if the counts on this ring were all deterministic, $\eta$
would be a Dirac mass.  Translation invariance would then make its
configuration translation-invariant, which is impossible for a locally
finite counting measure of positive intensity.

\emph{Sets with variance bounded away from zero.}
Choose $R_n\to\infty$.  For each $n$, ergodicity and the mean
ergodic theorem give
\[
\left\|\frac1{\Vol_d(Q_T)}\int_{Q_T+v}
\bigl(N_{A+s}-\rho\Vol_d(A)\bigr)\,ds\right\|_2\longrightarrow0
\qquad(T\to\infty),
\]
uniformly in $v$: the averages over $Q_T+v$ are unitary translates of
the average over $Q_T$.  Taking the inner product with
$N_{Q_{R_n}}-\rho\Vol_d(Q_{R_n})$ shows that
\[
\frac1{\Vol_d(Q_T)}\int_{Q_T+v}
\Cov_\eta(N_{Q_{R_n}},N_{A+s})\,ds\longrightarrow0
\qquad(T\to\infty)
\]
uniformly in $v$.  Choose $T$ so that this average is greater than
$-\sigma^2/4$.  Then choose $v$ so that every $A+s$, $s\in Q_T+v$,
is disjoint from $Q_{R_n}$.  Some $s_n\in Q_T+v$ satisfies
\[
(A+s_n)\cap Q_{R_n}=\varnothing,
\qquad
\Cov_\eta(N_{Q_{R_n}},N_{A+s_n})\geq-\frac{\sigma^2}{4}.
\]

Set
\[
E_n:=Q_{R_n}\sqcup(A+s_n).
\]
Then
\[
\Var_\eta(N_{E_n})
\geq
\sigma^2+2\Cov_\eta(N_{Q_{R_n}},N_{A+s_n})
\geq
\frac{\sigma^2}{2}.
\]
For every fixed $r>0$,
\[
\Vol_d(\partial_rE_n)
\leq
\Vol_d(\partial_rQ_{R_n})+\Vol_d(\partial_rA),
\]
where the second term is independent of the translation $s_n$ and
finite because $A$ is bounded.  The right-hand side is
$o(\Vol_d(Q_{R_n}))$, while $\Vol_d(E_n)\geq\Vol_d(Q_{R_n})$.

\emph{Variance amplification.}
Since $\Var_\eta(N_{E_n})>0$ and $\alpha_n>0$, choose an
integer $m_n$ large enough that
\[
 m_n^{\alpha_n}\geq
 \frac{4M_n\Vol_d(E_n)^{2-\alpha_n}}{
 \Var_\eta(N_{E_n})},
\]
or equivalently,
\[
\frac{m_n^2}{4}\Var_\eta(N_{E_n})
\geq
M_n\bigl(m_n\Vol_d(E_n)\bigr)^{2-\alpha_n}.
\]
Since $\partial E_n\subseteq\partial Q_{R_n}\cup(\partial A+s_n)$
and both sets have zero Lebesgue measure, the revival property gives,
for all sufficiently large \(k\),
\[
\|N_{E_n+t_k}-N_{E_n}\|_{L^2(\eta)}
\leq
\frac12\Var_\eta(N_{E_n})^{1/2}.
\]
Since $|t_k|\to\infty$, we may choose indices
$k_1<\cdots<k_{m_n}$, each satisfying the preceding $L^2$ bound,
with $r_{n,i}:=t_{k_i}$ and
\[
|r_{n,i}-r_{n,j}|>\operatorname{diam}(E_n)
\qquad(i\neq j).
\]
Then the translates \(E_n+r_{n,i}\) are pairwise disjoint.  Set
\[
F_n
:=
\bigsqcup_{i=1}^{m_n}(E_n+r_{n,i}).
\]
By \cref{lem:revival-covariance},
\[
\Var_\eta(N_{F_n})
\geq
\frac{m_n^2}{4}\Var_\eta(N_{E_n})
\geq
M_n\,\Vol_d(F_n)^{2-\alpha_n}.
\]

The sequence \((F_n)\) is van Hove: for every fixed \(r>0\),
\[
\partial_rF_n
\subseteq
\bigcup_{i=1}^{m_n}\bigl(\partial_rE_n+r_{n,i}\bigr),
\]
so
\[
\frac{\Vol_d(\partial_rF_n)}{\Vol_d(F_n)}
\leq
\frac{\Vol_d(\partial_rE_n)}{\Vol_d(E_n)}
\longrightarrow0\qquad(n\to\infty).
\]
Thus \(\eta\in\MFVH_\rho\).
\end{proof}

\begin{corollary}\label{cor:mfvh-residual}
The set \(\MFVH_\rho\) is residual in \(\smash{\EPP_{\rho,2}}\).
\end{corollary}

\begin{proof}
By \cref{prop:revival-residual,prop:revival-mfvh}, the residual set
\(\Rev_\rho\) is contained in \(\MFVH_\rho\).
\end{proof}

\appendix
\section{Rigidity sequences for actions}\label{subsec:generic-rigidity}

Let \(\operatorname{Aut}(X,\mu)\) be the group of measure-preserving
automorphisms modulo null sets, equipped with its weak topology.  This is a
Polish group; see \cite[Sec.~1]{Kec10}.  For an action \(a:\mathbb R^d\times X\longrightarrow X\), let
\(a_t(x):=a(t,x)\).  Identifying \(a\) with the continuous homomorphism
\(t\mapsto a_t\) into \(\operatorname{Aut}(X,\mu)\), equip
\(\Ad_d(X,\mu)\) with the topology of uniform convergence on compact subsets of
\(\mathbb R^d\).

For \(a\in\Ad_d(X,\mu)\), let \(U_t^a\) denote its Koopman operators on
\(L^2(X,\mu)\).  The action \(a\) is \emph{weakly mixing} if the diagonal
action on \((X\times X,\mu\times\mu)\) is ergodic.  A sequence
\((t_n)\subset\mathbb R^d\) is a \emph{rigidity sequence} for \(a\) if
\[
|t_n|\longrightarrow\infty
\qquad\text{and}\qquad
U_{t_n}^a\longrightarrow I
\quad\text{strongly on }L^2(X,\mu).
\]

\begin{proposition}\label{prop:generic-rigidity-actions}
The set of weakly mixing actions in \(\Ad_d(X,\mu)\) which admit a rigidity
sequence is a dense \(G_\delta\)-subset of \(\Ad_d(X,\mu)\).
\end{proposition}

For $a\in\Ad_d(X,\mu)$ and a bounded half-open rectangular box
$Q\subset\mathbb R^d$, a \emph{$Q$-tower} consists of a Borel set
$B\subset X$ and a finite measure $\nu_B$ on $B$ for which
\[
\Phi:Q\times B\longrightarrow X,\qquad (t,y)\longmapsto a_t y
\]
is injective and Borel, and
$\mu|_{\Phi(Q\times B)}=\Phi_*(\Vol_d|_Q\otimes\nu_B)$.
The measure $\nu_B$ is the \emph{transverse measure} on the base $B$;
it is not the restriction of $\mu$ to $B$.  The remainder is
$X\setminus\Phi(Q\times B)$.
The rectangular tower theorem of
\cite[Theorem~1 and Corollary~1, pp.~177--178]{Lin75}
gives the following for essentially free measure-preserving actions:
for every $\varepsilon>0$ there exists $L_0>0$ such that, whenever
$Q=\prod_{i=1}^d[0,L_i)$ with $L_i\geq L_0$ for all $i$, there is a
$Q$-tower satisfying
\[
\mu(X\setminus\Phi(Q\times B))<\varepsilon.
\]
Its transverse base measure is nonatomic.  Indeed, an atom of
$\nu_B$ would give positive $\mu$-measure to a bounded part of one
orbit, contradicting finiteness and translation invariance of $\mu$.

\begin{lemma}[Density of conjugacy classes]\label{lem:conjugacy-dense}
Every essentially free action in $\Ad_d(X,\mu)$ has a dense
conjugacy class.
\end{lemma}

\begin{proof}
Fix two essentially free actions $a,b$.  We show that a conjugate
of $b$ approximates $a$.  Fix compact $K\subset\mathbb R^d$, measurable
$A_1,\ldots,A_m\subset X$, and $\varepsilon>0$.
Choose a half-open rectangular box $Q$ with
\[
\Vol_d\{t\in Q:t+K\not\subset Q\}
<\frac{\varepsilon}{8}\Vol_d(Q).
\]
Take $Q$-towers for $a$ and $b$, with remainders of measure less
than $\varepsilon/32$ each, and bases $B_a,B_b$.
Their transverse base measures are nonatomic. Restrict the larger
base so that both have the same transverse measure, then choose a
measure-preserving isomorphism $J:B_b\to B_a$. The two restricted
towers have the same measure; their combined remainders have measure
less than $\varepsilon/8$. Match these remainders to extend
\[
T(b_t y):=a_t(J(y)),\qquad y\in B_b,\quad t\in Q,
\]
to a measure-space automorphism of $X$.
For $c_s:=Tb_sT^{-1}$ and $t,t+s\in Q$,
\[
c_s(a_t x)=a_{t+s}x=a_s(a_t x),\qquad x\in B_a.
\]
Uniformly for $s\in K$, the set on which $c_s$ and $a_s$ differ
has measure at most the remainder measure plus
$\Vol_d(\{t\in Q:t+K\not\subset Q\})\nu_{B_a}(B_a)
<\varepsilon/4$. Hence
\[
\sup_{s\in K}\mu(c_sA_i\triangle a_sA_i)<\varepsilon
\qquad(1\leq i\leq m).
\]
Thus $TbT^{-1}$ belongs to the prescribed neighbourhood of $a$.
\end{proof}

\begin{proof}
Let $(f_j)_{j\geq1}\subset L^2_0(X,\mu)$ be a countable dense
family.  Weak mixing is a $G_\delta$ condition.  Indeed, the mean
ergodic theorem applied to the diagonal action (equivalently, to the
tensor-square representation) shows that $a$ is weakly mixing precisely
when, for every $i,j,k,\ell\geq1$, there is an integer $R\geq\ell$
such that
\[
\frac1{\Vol_d(Q_R)}
\int_{Q_R}
\left|
\langle U_t^a f_i,f_j\rangle
\right|^2\,dt
<
\frac1k.
\]
For fixed $R,i,j$, compact-open convergence of actions gives uniform
convergence of the coefficients $\langle U_t^a f_i,f_j\rangle$ for
$t\in Q_R$.  Their averaged squared moduli therefore depend
continuously on $a$.

Choose a sequence $(h_j)_{j\geq1}$ dense in the unit ball of
$L^2(X,\mu)$.  The existence of a rigidity sequence is equivalent to
\[
\forall\,m,\ell\geq1\quad
\exists\,t\in\mathbb R^d,\ |t|>\ell:
\quad
\max_{1\leq j\leq m}
\|U_t^a h_j-h_j\|_2<\frac1m.
\]
For each fixed $t$, the finitely many inequalities define an open set in
the action topology.  Taking the union over $|t|>\ell$ gives an open set.  Conversely, applying it with $m=\ell=n$ and selecting
one $t_n$ for each $n$ gives $|t_n|\to\infty$ and strong convergence to
the identity on a dense family, hence on all of $L^2(X,\mu)$ by
unitarity.

Danilenko and Ryzhikov construct a rank-one weakly mixing rigid flow
\((V_t)_{t\in\mathbb R}\) with a rigidity sequence \(s_n\) \cite{DR12,DR13}.
Such a flow is essentially free: by ergodicity its stabilizer is almost
surely a fixed closed subgroup of \(\mathbb R\); nonatomicity excludes
\(\mathbb R\), while a stabilizer \(c\mathbb Z\), \(c>0\), would make the
flow factor through \(\mathbb R/c\mathbb Z\) and hence give a nonconstant
eigenfunction, contradicting weak mixing.  The coordinatewise action
\[
(t_1,\ldots,t_d).(x_1,\ldots,x_d)
:=
(V_{t_1}x_1,\ldots,V_{t_d}x_d)
\]
on the \(d\)-fold product is therefore essentially free and weakly mixing,
and \((s_n,\ldots,s_n)\) is a rigidity sequence because \(V_{s_n}\to I\)
strongly.  Transporting the product action to $(X,\mu)$ gives an element of
$\Ad_d(X,\mu)$ with both properties.  Both properties are invariant under conjugacy. By
\cref{lem:conjugacy-dense}, their intersection is dense.
\end{proof}

\begin{proposition}\label{prop:revival-passes-to-factors}
Let \(a\in\Ad_d(X,\mu)\) and
\(\kappa\in\smash{\RR^a_{\rho,2}(X,\mu)}\).  If \(a\) is weakly mixing, then the point process \(\kappa_*\mu\) is weakly mixing.  If \((t_n)\) is a rigidity sequence for \(a\),
then \((t_n)\) is a revival sequence for \(\kappa_*\mu\).
\end{proposition}

\begin{proof}
With $\eta:=\kappa_*\mu$, weak mixing passes to factors.  If $(t_n)$
is a rigidity sequence for $a$, fix a bounded Borel continuity set $A$.
The function $N_A\circ\kappa$ belongs to $L^2(X,\mu)$, and equivariance
identifies its Koopman translates with the translated window counts.
Thus
\[
\|N_{A+t_n}-N_A\|_{L^2(\eta)}
=
\|N_{A+t_n}\circ\kappa-N_A\circ\kappa\|_{L^2(\mu)}
\longrightarrow0\qquad(n\to\infty),
\]
by strong convergence of $U_{t_n}^a$ to the identity.
\end{proof}

\begin{proposition}[Revival and Koopman rigidity]
\label{prop:revival-koopman}
Let $\eta\in\EPP_{\rho,2}$ and let $(t_n)\subset\mathbb R^d$ satisfy
$|t_n|\to\infty$. Then $(t_n)$ is a revival sequence for $\eta$
if and only if
\[
F\circ\theta_{t_n}\longrightarrow F\quad\text{in }L^2(\eta)
\qquad(F\in L^2(\eta)).
\]
\end{proposition}

\begin{proof}
Strong Koopman convergence implies revival by applying it to $N_A$
for each bounded Borel set $A$ with $\Vol_d(\partial A)=0$.
Conversely, suppose $(t_n)$ is a revival sequence. Let $\mathscr Q$
be the countable family of bounded half-open rectangles with rational
endpoints. For each $Q\in\mathscr Q$,
\[
\|N_Q\circ\theta_{t_n}-N_Q\|_2
=\|N_{Q+t_n}-N_Q\|_2\longrightarrow0\qquad(n\to\infty).
\]
These counts are integer-valued. Hence, for every finite collection
$Q_1,\ldots,Q_m\in\mathscr Q$, the probability that any of their
counts changes under $\theta_{t_n}$ tends to zero. Every bounded
function of $(N_{Q_1},\ldots,N_{Q_m})$ therefore converges to its
Koopman translate in $L^2(\eta)$. Such functions are dense in
$L^2(\eta)$ because the counts on $\mathscr Q$ generate the Borel
$\sigma$-algebra of $\operatorname{Conf}(\mathbb R^d)$. The Koopman
operators are unitary, so convergence extends to every $F\in L^2(\eta)$.
\end{proof}

\begin{proposition}[Rigidity sequences and linear maximal rigidity]
\label{prop:rigidity-implies-linear}
Let $\eta\in\PP_{\rho,2}$. Suppose that there are $t_n\in\mathbb R^d$
with $|t_n|\to\infty$ such that
\[
F\circ\theta_{t_n}\longrightarrow F\quad\text{in }L^2(\eta)
\qquad(F\in L^2(\eta)).
\]
Then $\eta\in\LMRig_\rho$.
\end{proposition}

\begin{proof}
Fix a bounded Borel set $A$ and a bounded Borel function $h$ vanishing
outside $A$. Put $h_n(u):=h(u-t_n)$. For all sufficiently large $n$,
$h_n$ vanishes on $A$. Translation invariance and the definition of
$\theta_t$ give
\[
Z_{h_n}^\eta=Z_h^\eta\circ\theta_{t_n}
\longrightarrow Z_h^\eta\quad\text{in }L^2(\eta)\text{ as }n\to\infty.
\]
Thus $Z_h^\eta$ belongs to the closed linear span of centered exterior
statistics.
\end{proof}

In particular, \cref{prop:revival-koopman} gives
$\Rev_\rho\subseteq\LMRig_\rho$.
If $a\in\Ad_d(X,\mu)$ admits a rigidity sequence and
$\kappa\in\RR^a_{\rho,2}(X,\mu)$, then
$\kappa_*\mu\in\LMRig_\rho$: the Koopman representation of the
factor is a subrepresentation of that of $a$.

\section{Explicit intermittency and linear rigidity from irrational rotations}
\label{app:rotation-examples}

Let $\mathbb T:=\mathbb R/\mathbb Z$ and $I:=[0,1/2)\subset\mathbb T$.
For an irrational $\alpha\in(0,1)$ and $(\vartheta,u)\in\mathbb T\times[0,1)$,
set
\[
\Lambda_{\alpha,\vartheta,u}
:=\{n+u:n\in\mathbb Z,\ \vartheta+n\alpha\in I\}.
\]
Define the probability measure $\eta_\alpha$ on
$\operatorname{Conf}(\mathbb R)$ by
\[
\int F\,d\eta_\alpha
:=\int_0^1\int_{\mathbb T}
F(\delta_{\Lambda_{\alpha,\vartheta,u}})\,d\vartheta\,du
\qquad(F\text{ bounded Borel}).
\]
This is the cut-and-project process of \cite{BH24} for the lattice
\[
\Gamma_\alpha:=\{(n,n\alpha):n\in\mathbb Z\}
<\mathbb R\times\mathbb T
\]
and window $I$.  Within the parameter interval $(0,1)$, the
hyperuniform parameters have full Lebesgue measure but form a meagre
set, whereas intermittent hyperuniformity holds for a dense
$G_\delta$-set.  The genericity results in the main part of the paper
refer instead to Baire category on $\EPP_{\rho,2}$, which carries no
specified probability measure.

Write the regular continued fraction of $\alpha$ as
$\alpha=[0;a_1,a_2,\ldots]$.  For $j\geq1$, let $p_j/q_j$ be
its $j$th convergent, in lowest terms, and put $q_0:=1$ and
$q_{-1}:=0$.  Then
$q_{j+1}=a_{j+1}q_j+q_{j-1}$ and $q_j\to\infty$ as $j\to\infty$.

\begin{proposition}\label{prop:rotation-intermittency}
For every irrational $\alpha\in(0,1)$, the measure $\eta_\alpha$ is
translation-invariant and ergodic, has intensity $1/2$, and belongs to
$\LMRig_{1/2}$.  It is supported on uniformly discrete and relatively
dense point configurations.  Moreover,
\[
V_{\eta_\alpha}(q_j/2)\leq4\qquad(j\geq1).
\]
There is a dense $G_\delta$-set $\mathscr G\subset(0,1)$ of Liouville
numbers such that
\[
\eta_\alpha\in\IHyp_{1/2}\cap\LMRig_{1/2}
\qquad(\alpha\in\mathscr G).
\]
One explicit parameter in $\mathscr G$ is
\[
\alpha_*:=\sum_{n=1}^{\infty}3^{-n!}.
\]
\end{proposition}

For $\beta\in\mathbb R$ and an integer $N\geq1$, set
\[
D_N^\beta(\vartheta)
:=\sum_{j=0}^{N-1}
\bigl(\mathbf1_I(\vartheta+j\beta)-1/2\bigr).
\]

\begin{lemma}[Rational shadowing]\label{lem:rotation-shadowing}
Let $p,q$ be relatively prime integers with $q$ positive and odd.
If $M\geq1$ is an integer and
\[
\left|\beta-\frac pq\right|<\frac{1}{100Mq^2},
\]
then
\[
\int_{\mathbb T}|D_{Mq}^\beta(\vartheta)|^2\,d\vartheta
\geq\frac6{25}M^2.
\]
\end{lemma}

\begin{proof}
For $0\leq j<Mq$, the circle distance between
$\vartheta+j\beta$ and $\vartheta+jp/q$ is less than $1/(100q)$.
Let $B$ be the set of phases for which one of the $q$ points
$\vartheta+jp/q$, $0\leq j<q$, lies within $1/(100q)$ of either endpoint
of $I$.  The union bound gives
\[
|B|\leq q\cdot2\cdot\frac{2}{100q}=\frac1{25}.
\]
Outside $B$, the two itineraries agree for all $Mq$ steps.
Since $q$ is odd, a rational period contains either $(q-1)/2$ or
$(q+1)/2$ visits to $I$.  Its centered discrepancy is therefore
$1/2$ or $-1/2$, with the same sign in every period.  Thus
$|D_{Mq}^\beta|=M/2$ outside $B$, and
\[
\int_{\mathbb T}|D_{Mq}^\beta|^2
\geq\frac{24}{25}\frac{M^2}{4}
=\frac6{25}M^2.
\]
\end{proof}

\begin{proof}[Proof of \cref{prop:rotation-intermittency}]
\emph{Translation invariance and geometry.}
For $t\in\mathbb R$, put $k=\lfloor u-t\rfloor$.  Translation of a
configuration by $-t$ corresponds to the transformation
\[
(\vartheta,u)\longmapsto
(\vartheta-k\alpha,\ u-t-k).
\]
This action preserves the product of Haar measure and Lebesgue measure.
Its time-one map is $(\vartheta,u)\mapsto(\vartheta+\alpha,u)$.
A function invariant under the whole action is therefore independent
of $\vartheta$, by irrationality of $\alpha$, and invariance under
the remaining times makes it constant in $u$.  Hence the action,
and its point-process factor, are ergodic.

For every bounded Borel set $A\subset\mathbb R$,
\[
\int N_A\,d\eta_\alpha
=\frac12\int_0^1\sum_{n\in\mathbb Z}\mathbf1_A(n+u)\,du
=\frac12|A|.
\]
All local counts are bounded by those of $\mathbb Z+u$, and distinct
points have distance at least one.  Minimality of the irrational
rotation, applied to $\operatorname{int}I$, gives a finite bound on
return gaps, uniform in the phase.  Thus the configurations are
relatively dense.

\smallskip\noindent
\emph{Linear maximal rigidity.}
Let $A\subset\mathbb R$ be bounded and let $h$ be a bounded Borel
function vanishing outside $A$.  Choose an integer $K$ bounding the
number of sites of every translate $\mathbb Z+u$ in a bounded interval
containing $A$.  Put
\[
h_j(v):=h(v-q_j),\qquad
\varepsilon_j:=\operatorname{dist}(q_j\alpha,\mathbb Z).
\]
The supports of $h_j$ are disjoint from $A$ for all sufficiently large $j$,
and $\varepsilon_j\to0$.  Reindexing the sums gives
\[
\begin{split}
S_{h_j}-S_h
=\sum_{n\in\mathbb Z}h(n+u)
\bigl(&\mathbf1_I(\vartheta+n\alpha+q_j\alpha)\\
&-\mathbf1_I(\vartheta+n\alpha)\bigr).
\end{split}
\]
There are at most $K$ nonzero summands.  For each summand the two
indicators disagree on a set of phases of measure $2\varepsilon_j$.
Cauchy--Schwarz and integration therefore give
\[
\|S_{h_j}-S_h\|_{L^2(\eta_\alpha)}^2
\leq2K^2\|h\|_\infty^2\varepsilon_j\longrightarrow0\qquad(j\to\infty).
\]
Since $\int h_j=\int h$, the same convergence holds for the centered
statistics.  Hence $Z_h^{\eta_\alpha}$ lies in the closed linear span of exterior
centered statistics, proving $\eta_\alpha\in\LMRig_{1/2}$.

\smallskip\noindent
\emph{Bounded variance.}
For an integer $N\geq1$, the interval $[0,N)$ contains exactly the
sites $u,u+1,\ldots,u+N-1$.  Translation invariance therefore gives
\begin{equation}\label{eq:rotation-variance-discrepancy}
V_{\eta_\alpha}(N/2)
=\int_{\mathbb T}|D_N^\alpha(\vartheta)|^2\,d\vartheta.
\end{equation}
The function $\mathbf1_I-1/2$ has mean zero and total variation two.
The Denjoy--Koksma inequality \cite[Chap.~VI]{Her79} yields
$|D_{q_j}^\alpha|\leq2$, and hence
$V_{\eta_\alpha}(q_j/2)\leq4$.
Since $q_j\to\infty$, for every $\beta>0$ we have
\[
\lim_{j\to\infty}
\frac{V_{\eta_\alpha}(q_j/2)}{(q_j/2)^\beta}=0.
\]

\smallskip\noindent
\emph{Large variance.}
For $k\geq6$, set
\[
\mathscr O_k:=
\bigcup_{\substack{q\geq k\ {\rm odd},\ p\in\mathbb Z\\(p,q)=1}}
\left(\frac pq-q^{-k},\frac pq+q^{-k}\right)\cap(0,1),
\qquad
\mathscr G:=((0,1)\setminus\mathbb Q)\cap
\bigcap_{k\geq6}\mathscr O_k.
\]
Reduced rationals with arbitrarily large odd denominators are dense.
Thus every $\mathscr O_k$ is open and dense, and $\mathscr G$ is a
dense $G_\delta$-set consisting of Liouville numbers.

For $\alpha\in\mathscr G$ and each $k\geq6$, choose $p,q$ as in the
definition of $\mathscr O_k$ and put
\[
M:=\lfloor q^{k-2}/100\rfloor,\qquad N:=Mq.
\]
Then $M\geq q^{k-2}/200$ and
$|\alpha-p/q|<1/(100Mq^2)$.  The rational-shadowing lemma gives
\[
V_{\eta_\alpha}(N/2)\geq\frac6{25}M^2.
\]
Consequently, for every fixed $\gamma<2$,
\[
\frac{V_{\eta_\alpha}(N/2)}{N^\gamma}
\geq
\frac6{25}\,200^{-(2-\gamma)}
q^{(k-2)(2-\gamma)-\gamma}
\longrightarrow+\infty.
\]
Here $q\geq k$, so $N\to\infty$.  Since $R=N/2\to\infty$ as $k\to\infty$, this gives
$\limsup_{R\to\infty}V_{\eta_\alpha}(R)/R^\gamma=\infty$
for every $\gamma<2$.  Together with the bounded-variance subsequence it gives
$\eta_\alpha\in\IHyp_{1/2}$.

\smallskip\noindent
\emph{Explicit parameter.}
In the base-$3$ expansion of $\alpha_*=\sum_{n\geq1}3^{-n!}$,
the digit at position $n!$ is $1$, and all other digits are $0$.
The increasing gaps between these positions show that the base-$3$
digit sequence is not eventually periodic, so $\alpha_*$ is irrational.
The $n$th partial sum of the displayed series is $p_n/q_n$, where
$q_n=3^{n!}$ and $p_n\equiv1\pmod3$.
Thus the fraction is reduced and $q_n$ is odd.  Moreover,
\[
0<\alpha_*-\frac{p_n}{q_n}
<2\,3^{-(n+1)!}
=2q_n^{-(n+1)}.
\]
For every $k\geq6$, choosing $n$ sufficiently large gives
$q_n\geq k$ and $|\alpha_*-p_n/q_n|<q_n^{-k}$.
Hence $\alpha_*\in\mathscr G$.
\end{proof}

\begin{corollary}[Category and measure in the rotation family]
\label{cor:rotation-category-measure}
The set of irrational $\alpha\in(0,1)$ for which $\eta_\alpha$ is
hyperuniform has full Lebesgue measure but is meagre in $(0,1)$.
The intermittently hyperuniform parameters contain the dense
$G_\delta$-set $\mathscr G$ from \cref{prop:rotation-intermittency}.
Every irrational parameter gives a linearly maximally rigid process.
\end{corollary}

\begin{proof}
The last two assertions follow from \cref{prop:rotation-intermittency}.
Since intermittent hyperuniformity excludes hyperuniformity, the set of
hyperuniform parameters is contained in the meagre complement of
$\mathscr G$.

Fix $0<\varepsilon<1/2$. For each integer $q\geq1$, the set
\[
\{\alpha\in(0,1):\operatorname{dist}(q\alpha,\mathbb Z)
<q^{-1-\varepsilon}\}
\]
has Lebesgue measure at most $2q^{-1-\varepsilon}$.
By Borel--Cantelli, for almost every $\alpha$ there is $q_0$ such that
$\operatorname{dist}(q\alpha,\mathbb Z)\geq q^{-1-\varepsilon}$
for $q\geq q_0$. If $p_j/q_j$ is a continued-fraction convergent of
such an irrational $\alpha$, then
\[
q_j^{-1-\varepsilon}
\leq |q_j\alpha-p_j|<q_{j+1}^{-1}
\]
for all sufficiently large $j$. Consequently
$a_{j+1}\leq C_\alpha q_j^\varepsilon$, where $a_{j+1}$ is the
next continued-fraction coefficient.

Fix an integer $N\geq1$.  Choose $J$ maximal with $q_J\leq N$
and subtract successively the largest possible multiples of
$q_J,q_{J-1},\ldots,q_0$.  The recurrence for $q_j$ gives
\[
N=\sum_{j=0}^J b_jq_j,\qquad 0\leq b_j\leq a_{j+1}.
\]
Since $q_{j+2}\geq2q_j$, we have $J=O_\alpha(\log(2N))$.
Partition the first $N$ iterates into consecutive blocks of lengths
$q_j$, with $b_j$ blocks of length $q_j$.  Applying Denjoy--Koksma
at the initial phase of each block gives
\[
\sup_{\vartheta\in\mathbb T}|D_N^\alpha(\vartheta)|
\leq 2\sum_{j=0}^J b_j
\leq C_\alpha N^\varepsilon\log(2N).
\]
By \eqref{eq:rotation-variance-discrepancy},
\[
V_{\eta_\alpha}(N/2)
\leq C_\alpha N^{2\varepsilon}\log^2(2N)=o(N).
\]
If $2R\in[N,N+1)$, the counts in $[0,2R)$ and $[0,N)$ differ by
at most one. Their centered counts therefore differ in $L^2$ by
at most two. Translation invariance identifies their variances with those of
balls of radii $R$ and $N/2$. Hence $V_{\eta_\alpha}(R)=o(R)$,
which is hyperuniformity in dimension one.
\end{proof}

\section*{Funding}

This work was supported by the Swedish Research Council under grant
VR~11253322.

\section*{Use of AI tools}

During the preparation of this manuscript, the author used
ChatGPT (OpenAI) and Astra for mathematical exploration, proof analysis,
and drafting; Aristotle (Harmonic) for substantial formalization and
verification of selected arguments in Lean; and Claude Sonnet
(Anthropic) for an independent adversarial mathematical review.
The author reviewed and revised the AI-assisted contributions and takes
full responsibility for the mathematical results and the final manuscript.

\end{document}